\documentclass[10pt,reqno,twoside]{amsart}

\usepackage{amsfonts}
\usepackage{paralist}
  \usepackage{graphics}
  \usepackage{epsfig}
 \usepackage[colorlinks=true]{hyperref}
\usepackage{hyperref}

\numberwithin{equation}{section}
\usepackage{latexsym}
\usepackage{amsmath}
\usepackage{amssymb}
\usepackage{mathrsfs}
\usepackage{graphicx}
\usepackage{ifthen}
\usepackage[T1]{fontenc}
\usepackage[latin1]{inputenc}
\usepackage{color}
\newtheorem{Satz}{Theorem}[section]

\newtheorem{Def}[Satz]{Definition}

\newcommand{\oo}{\overline \Omega}

\newcommand{\po}{\partial\Omega}

\newcommand{\diag}{\text {diag}}

\newcommand{\supp}{\text {supp}}

\thanks{$^*$This research is supported by
the Student Research Training of Jiangsu Province
under Grant
No. 202110307029Y, and
the National Natural Science Foundation of China
under Grant Nos. 11601232,   11671354 and  11775116.
}

\begin{document}

\title[Bubbling solutions]
{Bubbling solutions for a planar exponential nonlinear elliptic equation with a singular source}

\maketitle

%\centerline{\scshape Haitao Yang}
%\medskip
%{\footnotesize
% \centerline{School of Mathematics, Zhejiang University, Hangzhou 310027, China}
% \smallskip
% \centerline{Email:\,\,\,\,htyang@css.zju.edu.cn}
%}
%
%\bigskip

\centerline{\scshape 
Jingyi Dong,   Jiamei Hu,  Yibin Zhang\footnote{Corresponding author: yibin10201029@njau.edu.cn} 
}
\medskip
{\footnotesize
 \centerline{College of Sciences, Nanjing Agricultural University, Nanjing 210095, China}
%  \smallskip
% \centerline{Email:\,\,\,\,yibin10201029@njau.edu.cn}
}

%\author[Yibin Zhang]{Yibin Zhang}
%\address{College of Sciences, Nanjing Agricultural University, Nanjing 210095, China}
%\email{yibin10201029@njau.edu.cn}

\begin{abstract}
Let $\Omega$ be a bounded domain in $\mathbb{R}^2$ with smooth boundary,
we study the following  elliptic Dirichlet problem
$$
\begin{cases}
-\Delta\upsilon=
e^{\upsilon}-s\phi_1-4\pi\alpha\delta_p-h(x)\,\,\,\,
\,\textrm{in}\,\,\,\,\,\Omega,\\[2mm]
\upsilon=0
\quad\quad\quad\quad\quad\quad
\qquad\qquad\quad\quad\,\,\,\,
\textrm{on}\,\ \,\partial\Omega,
\end{cases}
$$
where $s>0$ is a large  parameter,
$h\in C^{0,\gamma}(\overline{\Omega})$,
$p\in\Omega$, $\alpha\in(-1,+\infty)\setminus\mathbb{N}$,
$\delta_p$ denotes the Dirac measure supported at point $p$
and $\phi_1$ is a positive first
eigenfunction of
the problem
$-\Delta\phi=\lambda\phi$
under Dirichlet boundary condition in $\Omega$.
If $p$ is a strict local maximum point of $\phi_1$,
we show that such a problem has a family of solutions $\upsilon_s$
with arbitrary $m$ bubbles accumulating to $p$,
and the quantity $\int_{\Omega}e^{\upsilon_s}\rightarrow8\pi(m+1+\alpha)\phi_1(p)$
as $s\rightarrow+\infty$.
\end{abstract}

\newpage

\section{Introduction}
Let $\Omega$ be a bounded domain in $\mathbb{R}^2$ with smooth boundary.
This paper deals with the analysis of solutions in the distributional
sense for the following problem  involving a  singular source
\begin{equation}\label{1.1}
\aligned
\left\{\aligned
&-\Delta\upsilon=
e^{\upsilon}-s\phi_1-4\pi\alpha\delta_p-h(x)\,\,\,\,
\,\textrm{in}\,\,\,\,\,\Omega,\\[2mm]
&\upsilon=0
\quad\quad\quad\quad\quad\quad
\quad\qquad\qquad\quad\quad
\textrm{on}\,\ \,\partial\Omega,
\endaligned\right.
\endaligned
\end{equation}
where
$s>0$ is a large  parameter, $p\in\Omega$, $\alpha\in(-1,+\infty)\setminus\mathbb{N}$,
$\delta_p$ denotes the Dirac measure supported at point $p$,
$h\in C^{0,\gamma}(\overline{\Omega})$ is given,
$\phi_1>0$ is an eigenfunction of
$-\Delta$ with Dirichlet boundary condition corresponding to the
first eigenvalue $\lambda_1$.
Clearly, if we set $\rho(x)=(-\Delta)^{-1}h$ in $H_0^1(\Omega)$ and
let $G(x,y)$ be the Green's function associated to $-\Delta$ with Dirichlet
boundary condition, namely
\begin{equation}\label{1.2}
\left\{\aligned
&-\Delta_xG(x,y)=8\pi\delta_{y}(x),
\,\,\,\,\,\,\,\,
x\in\Omega,\\[2mm]
&G(x,y)=0,
\,\qquad\qquad\quad\,\,\quad\,
x\in\partial\Omega,
\endaligned\right.
\end{equation}
and $H(x,y)$ be its regular part defined as
\begin{equation}\label{1.3}
\aligned
H(x,y)=G(x,y)-4\log\frac1{|x-y|},
\endaligned
\end{equation}
then equation (\ref{1.1}) is equivalent to solving for
$u=\upsilon+\frac{s}{\lambda_1}\phi_1+\frac{\alpha}{2}G(\,\cdot,\,p)+\rho$, the problem
\begin{equation}\label{1.4}
\aligned
\left\{\aligned
&-\Delta u=|x-p|^{2\alpha}
k(x)e^{-t\phi_1}e^u\,\,\,\,\,
\textrm{in}\,\,\,\,\,\Omega,\\[2mm]
&u=0\,\,\qquad\,\,\qquad\,\,
\quad\quad\qquad\qquad\,\textrm{on}
\,\,\,\partial\Omega,
\endaligned\right.
\endaligned
\end{equation}
where $k(x)=e^{-\rho(x)-\frac{\alpha}{2}H(x,p)}$ and
$t=s/\lambda_1$.
We are interested in the existence of solutions of problem (\ref{1.4})
(or (\ref{1.1}))
which exhibit the {\it concentration phenomenon}  when
the  parameter $t\rightarrow+\infty$.

This work is  directly  motivated by
the study of the regular case $\alpha=0$ in
equation (\ref{1.1}), namely
the following elliptic equation of
{\it Ambrosetti-Prodi type} \cite{AP}:
\begin{equation}\label{1.5}
\left\{\aligned
&-\Delta \upsilon=e^\upsilon-s\phi_1-h(x)\,\,\,\,
\,\textrm{in}\,\,\,\,\,\Omega,\\[2mm]
&\upsilon=0\,\,\quad\quad\,\,\,\,\quad\quad\quad\,
\quad\quad\,\,\,\,\,\,\textrm{on}\,\,\,\,\po,
\endaligned\right.
\end{equation}
or its equivalent form
\begin{equation}\label{1.6}
\aligned
\left\{\aligned
&-\Delta
u=k(x)e^{-t\phi_1}e^u\,\,\,\,\,
\textrm{in}\,\,\,\,
\Omega,\\[2mm]
&u=0\quad\,\,\quad
\quad\quad\quad\quad\quad\,
\textrm{on}\,\,\,\po,
\endaligned\right.
\endaligned
\end{equation}
where $\Omega$ is a bounded smooth domain in $\mathbb{R}^N$ ($N\geq2$).
In the early 1980s, Lazer and McKenna conjectured that (\ref{1.5}) has
an unbounded number of solutions as $s\rightarrow+\infty$ (see \cite{LM}).
When $N=2$, del Pino and Mu\~{n}oz \cite{DM} proved the Lazer-McKenna conjecture
for problem (\ref{1.5})
by constructing  {\it non-simple}  bubbling
solutions of  (\ref{1.6}) with the following
properties
\begin{equation}\label{1.7}
\aligned
k(x)e^{-t\phi_1}e^{u_t}\rightharpoonup8\pi\sum_{i=1}^m m_i\delta_{\xi_i}
\,\,\,\quad\,\,\,
\textrm{and}
\,\,\,\quad\,\,\,
u_t=\sum\limits_{i=1}^{m}m_iG(x,\xi_i)+o(1),
\endaligned
\end{equation}
where $m_i>1$ and  $\xi_i$'s are distinct maxima of $\phi_1$.
Surprising enough,  this multiple bubbling phenomenon  is in strong
opposition to a slightly modified but widely studied  version of equation (\ref{1.6}), namely the
Liouville-type  equation or sometimes referred to as the Gelfand equation
\begin{equation}\label{1.8}
\aligned
\left\{\aligned
&-\Delta
u=\varepsilon^2k(x)e^u\,\,\,\,\,
\textrm{in}\,\,\,\,
\Omega,\\[2mm]
&u=0\quad\quad
\quad\quad\quad\quad
\textrm{on}\,\,\,\po,
\endaligned\right.
\endaligned
\end{equation}
with $\varepsilon\rightarrow0$,
where
$\Omega\subset\mathbb{R}^2$ is a  bounded smooth domain
and $k(x)\in\mathcal{C}^2(\overline{\Omega})$ is
a non-negative, not identically zero function.
Indeed, the asymptotic analysis in
\cite{BM,LS,MW,NS} shows  that if $u_\varepsilon$ is an
unbounded family of
solutions of equation (\ref{1.8}) for which
$\varepsilon^2\int_{\Omega}k(x)e^{u_\varepsilon}$ is uniformly bounded,
then, up to a subsequence, there exists an integer $m\geq1$ such that
$u_\varepsilon$ makes
$m$  distinct  points  {\it simple} blow-up on
$\mathcal{S}=\{\xi_1,\ldots,\xi_m\}\subset\widetilde{\Omega}$ with
$\widetilde{\Omega}=\{x\in\Omega|\,k(x)>0\}$,
more precisely
\begin{equation}\label{1.9}
\aligned
\varepsilon^2k(x)e^{u_\varepsilon}\rightharpoonup8\pi\sum_{i=1}^m \delta_{\xi_i}
\,\,\ \,\quad\ \,\ \,\,
\textrm{and}
\,\,\,\quad\ \,\,\,
u_\varepsilon=\sum\limits_{i=1}^{m}G(x,\xi_i)+o(1).
\endaligned
\end{equation}
Also the location of $m$-tuple of these bubbling points  can be viewed  as a critical point
of a functional in terms of the Green's function and its regular part.
Conversely, the existence of solutions for equation
(\ref{1.8}) with these bubbling behaviors  has been founded in
\cite{BP,CL,DKM,EGP,W}.
In particular,
the construction of solutions with arbitrary  $m$ distinct bubbling points
is achieved  in  some special  cases:  for any $m\geq1$ if
$\Omega$ is not simply connected (\cite{DKM}),
and for any  $m\in\{1,\ldots,h\}$ provided that $\Omega$ is an
$h$-dumbell with thin handles (\cite{EGP}).
Finally, we mention that  in the recent paper \cite{YZ},
it has been proven that
the Lazer-McKenna conjecture also holds true
for problem (\ref{1.5}) in dimension $N\geq3$
with some  symmetries,
by constructing non-simple
bubbling solutions to a two-dimensional
anisotropic version of  (\ref{1.6}).

Our  motivation also directly comes from the study of
a slightly modified version of
equation (\ref{1.1}) (or (\ref{1.4})), namely
the singular Liouville equation
\begin{equation}\label{1.10}
\left\{\aligned
&-\Delta \upsilon=\varepsilon^2e^\upsilon-4\pi\alpha\delta_p-h(x)\,\,\,\,
\,\textrm{in}\,\,\,\,\,\Omega,\\[2mm]
&\upsilon=0\,\,\quad\quad\,\,\,\,\quad\quad\quad\,
\qquad\qquad\,\,\,\,\,\,\textrm{on}\,\,\,\,\po,
\endaligned\right.
\end{equation}
or its equivalent form
\begin{equation}\label{1.11}
\aligned
\left\{\aligned
&-\Delta
u=\varepsilon^2k(x)|x-p|^{2\alpha}e^u\,\,\,\,\,
\textrm{in}\,\,\,\,
\Omega,\\[2mm]
&u=0\quad\,\,\quad\,\,
\quad\qquad\qquad\quad\quad\,
\textrm{on}\,\,\,\po,
\endaligned\right.
\endaligned
\end{equation}
with $\Omega\subset\mathbb{R}^2$,
$-1<\alpha\neq0$
and $\varepsilon\rightarrow0$.
%This type of singular equation appears in the
%several superconductivity theories of the self-duel regime,
%such as the Abelian Maxwell-Higgs and Chern-Simons-Higgs theories
%\cite{NT,T}. It also arises in the  Tur-Yanovsky vortex pattern
%of planar stationary Euler equations for an incompressible and homogeneous
%fluid \cite{TY,DEM},  and in the construction of  planar conformal metrics  with conical
%singularity of order $\alpha$ \cite{LT}.
This type of singular equation arises in the  Tur-Yanovsky vortex pattern
of planar stationary Euler equations for an incompressible and homogeneous
fluid \cite{TY,DEM},   the construction of  planar conformal metrics  with conical
singularity of order $\alpha$ \cite{LT},
and  the
several superconductivity theories of the self-duel regime,
such as the Abelian Maxwell-Higgs and Chern-Simons-Higgs theories
\cite{NT,T}.

For equation (\ref{1.10}) involving $\alpha>0$,
solutions with $m$ distinct bubbling points away from the singular source
$p$ have been founded  first in \cite{DKM} provided that $m<1+\alpha$.
Later in \cite{D1}, this result has been extended to the case of multiple singular sources,
and more specifically,
it is shown that, under suitable restrictions on the weights, if several sources
exist, then the more involved topology should generate a larger number of bubbling solutions than
the singleton case considered in \cite{DKM}.  However, the problem  of
finding solutions of $(\ref{1.10})|_{\alpha>0}$
with additional bubbles around the singular source $p$ is of different nature. Indeed, the asymptotic analysis
in \cite{E,T1} shows that if these solutions exist, then the bubble at the singularity provides an
additional contribution of $8\pi(1+\alpha)\delta_p$ in the limit of (\ref{1.9}). More precisely,
if $u_{\varepsilon}$ is an unbounded family of solutions of $(\ref{1.11})|_{\alpha>0}$ for which
$\varepsilon^2\int_{\Omega}k(x)|x-p|^{2\alpha}e^{u_\varepsilon}$ is uniformly bounded
and $u_{\varepsilon}$ is unbounded in any neighborhood of $p$,
then, up to a subsequence, there exists an integer $m\geq0$ such that
$u_\varepsilon$ makes
$m+1$ distinct  points  {\it simple} blow-up on
$\mathcal{S}=\{p,\,\xi_1,\ldots,\xi_m\}\subset\Omega$,
namely
\begin{equation}\label{1.12}
\aligned
\varepsilon^2k(x)|x-p|^{2\alpha}e^{u_\varepsilon}\rightharpoonup
8\pi(1+\alpha)\delta_p+8\pi\sum_{i=1}^m \delta_{\xi_i}
\,\,\ \,\quad\
\textrm{and}
\,\,\,\quad\
u_\varepsilon=(1+\alpha)G(x,p)+\sum\limits_{i=1}^{m}G(x,\xi_i)+o(1).
\endaligned
\end{equation}
Moreover, the location of the $m$ distinct points
$\xi_1,\ldots,\xi_m\in\Omega\setminus\{p\}$
can be characterized as a critical point of some
certain functional in terms of
the Green's function and its regular part.
Reciprocally, the construction of solutions of equation
$(\ref{1.10})|_{\alpha>0}$  with  bubbles around
 $p$ has been carried out first in \cite{E} for the
 case $\alpha\in(0,+\infty)\setminus\mathbb{N}$,
later in \cite{DEM} for the case $\alpha\in\mathbb{N}$
where  given any positive integer $\alpha$ and any
sufficiently small complex number $a$, it is proven  that
there exists a solution of equation (\ref{1.10}) with $h(x)\equiv0$
and $\delta_p$ replaced by $\delta_{p_{a,\varepsilon}}$ for a suitable
$p_{a,\varepsilon}\in\Omega$ with $\alpha+1$ bubbling points at the vertices of
a sufficiently tiny regular polygon centered in point $p_{a,\varepsilon}$;
moreover  $p_{a,\varepsilon}$ lies close to a zero point of a vector field
explicitly built upon derivatives of order $\alpha+1$ of the regular part of
Green's function of the domain. Recently,
for equation (\ref{1.11})
with $\alpha\in\mathbb{N}$ and the potential
$k(x)$ replaced by $a(x)e^{-\frac12(\alpha-1)H(x,p)}$, it has been proven
in  \cite{D2} that if the local potential $a(x)$ and the geometry of the domain
satisfy some conditions at the singular source $p$, then there exists a solution $u_{\varepsilon}$
bubbling only at $p$  and  satisfying
$\varepsilon^2\int_{\Omega}a(x)e^{-\frac12(\alpha-1)H(x,p)}|x-p|^{2\alpha}e^{u_\varepsilon}\rightarrow8\pi(1+\alpha)$
as $\varepsilon\rightarrow0$.

In the present paper, we consider the singular case
of problem (\ref{1.1}) (or (\ref{1.4}))
involving $\alpha\in(-1,+\infty)\setminus\mathbb{N}$ and
try to prove the existence of its
non-simple bubbling solutions in a constructive way. We  find that if
the singular source $p$ is a strict local maximum point of $\phi_1$ in the domain,
then problem (\ref{1.4}) (or (\ref{1.1})) has a
family of solutions with the accumulation of arbitrarily many bubbles at source $p$.
This can be stated as following:

\vspace{1mm}
\vspace{1mm}
\vspace{1mm}
\vspace{1mm}

\noindent{\bf Theorem 1.1.} {\it
Let $\alpha\in(-1,+\infty)\setminus\mathbb{N}$ and assume that $p$ is a strict local maximum
point of $\phi_1(x)$ in $\Omega$.
Then for any integer $m\geq1$, there exists $t_m>0$
such that for any $t>t_m$,  problem {\upshape(\ref{1.4})} has
a family of solutions $u_t$
satisfying
%$$
%\aligned
%\lim_{t\rightarrow+\infty}\int_{\Omega}|x-p|^{2\alpha}
%k(x)e^{-t\phi_1}e^{u_t}=8\pi(m+1+\alpha)\phi_1(p).
%\endaligned
%$$
%More precisely,
$$
\aligned
u_t(x)=\left[\log
\frac{1}{(\varepsilon_{0,t}^2\mu_{0,t}^2+|x-p|^{2(1+\alpha)})^2}
+(1+\alpha)H(x,p)\right]
+
\sum\limits_{i=1}^{m}\left[\,\log
\frac1{(\varepsilon_{i,t}^2\mu_{i,t}^2+|x-\xi_{i,t}|^2)^2}
+H(x,\xi_{i,t})
\,\right]+o(1),
\endaligned
$$
where $o(1)\rightarrow0$, as $t\rightarrow+\infty$, uniformly on
each compact subset of $\overline{\Omega}\setminus\{p,\,\xi_{1,t},\ldots,\xi_{m,t}\}$,
the parameters $\varepsilon_{0,t}$, $\varepsilon_{i,t}$,  $\mu_{0,t}$  and
$\mu_{i,t}$ satisfy
$$
\aligned
\varepsilon_{0,t}=e^{-\frac12t\phi_1(p)},
\,\,\ \quad\ \,\,
\varepsilon_{i,t}=e^{-\frac12t\phi_1(\xi_{i,t})},
\,\,\ \quad\ \,\,
\frac1C\leq\mu_{0,t}\leq Ct^{2m\beta},
\,\,\ \quad\ \,\,
\frac1C\leq\mu_{i,t}\leq Ct^{(2m+\alpha)\beta},
\endaligned
$$
for  some   $C>0$,
and $(\xi_{1,t},\ldots,\xi_{m,t})\in\Omega^m$ satisfies
$$
\aligned
\xi_{i,t}\rightarrow p\,\ \ \,\,\textrm{for all}\,\,\,i,\,\,
\ \ \ \,
\textrm{and}\,\,\ \ \,\,\,|\xi_{i,t}-\xi_{j,t}|>t^{-\beta}\,\,\ \,\,\forall
\,\,\,
i\neq j,
\endaligned
$$
with $\beta=\frac{1}{2}(m+1)(m+1+\alpha)$.
}

\vspace{1mm}
\vspace{1mm}
\vspace{1mm}
\vspace{1mm}

The equivalent result for problem (\ref{1.1}) can be stated in the following form.

\vspace{1mm}
\vspace{1mm}
\vspace{1mm}
\vspace{1mm}

\noindent{\bf Theorem 1.2.} {\it
Let $\alpha\in(-1,+\infty)\setminus\mathbb{N}$ and assume that $p$ is a strict local maximum
point of $\phi_1(x)$ in $\Omega$.  Then for any integer $m\geq1$ and any
$s$ large enough,  there exists a  family of solutions $\upsilon_s$
of  problem {\upshape(\ref{1.1})} with $m$ distinct bubbles  accumulating to $p$. Moreover,
$$
\aligned
\lim_{s\rightarrow+\infty}\int_{\Omega}e^{\upsilon_s}=8\pi(m+1+\alpha)\phi_1(p).
\endaligned
$$
}

\vspace{1mm}

Moreover, for the case $m=0$,  we have
the corresponding results for problems
(\ref{1.1}) and (\ref{1.4}), respectively.

\vspace{1mm}
\vspace{1mm}
\vspace{1mm}
\vspace{1mm}

\noindent{\bf Theorem 1.3.} {\it
Let $\alpha\in(-1,+\infty)\setminus\mathbb{N}$.
Then   there exists $t_0>0$
such that for any $t>t_0$,  problem {\upshape(\ref{1.4})} has
a family of solutions $u_t$  such that as $t$ tends to $+\infty$,
$$
\aligned
u_t(x)=\left[\log
\frac{1}{(\mu_{0}^2e^{-t\phi_1(p)}+|x-p|^{2(1+\alpha)})^2}
+(1+\alpha)H(x,p)\right]
+o(1),
\endaligned
$$
 uniformly on
each compact subset of $\overline{\Omega}\setminus\{p\}$,
where  the parameter   $\mu_{0}$   satisfies $1/C\leq\mu_0\leq C$
for some $C>0$.
}

\vspace{1mm}
\vspace{1mm}
\vspace{1mm}
\vspace{1mm}

\noindent{\bf Theorem 1.4.} {\it
Let $\alpha\in(-1,+\infty)\setminus\mathbb{N}$.  Then for any
$s$ large enough,  there always exists a  family of solutions $\upsilon_s$
of  problem {\upshape(\ref{1.1})} such that
$$
\aligned
\lim_{s\rightarrow+\infty}\int_{\Omega}e^{\upsilon_s}=8\pi(1+\alpha)\phi_1(p).
\endaligned
$$
}

According to Theorems 1.1 and 1.2,  it follows  that if  the singular source
$p$ is an isolated local maximum point of $\phi_1$,
then for any  integer $m\geq1$
there exists a family of solutions  of problem (\ref{1.4}) which exhibits
the phenomenon of $m+1$-bubbling at $p$, namely,
$|x-p|^{2\alpha}k(x)e^{-t\phi_1}e^{u_t}\rightharpoonup 8\pi(m+1+\alpha)\delta_{p}$
and $u_t=(m+1+\alpha)G(x,p)+o(1)$. While
for the case $m=0$, by arguing exactly along the
sketch of the proof of Theorem 1.1 we can prove
the corresponding results in Theorems 1.3 and 1.4,
and further find that problem (\ref{1.4}) should always admit a family of solutions
blowing up at the singular source $p$ whether
$p$ is an isolated local maximum point of $\phi_1$
or not.

The strategy for proving  our main results  relies on a very
well-known  Lyapunov-Schmidt reduction procedure. In Section $2$ we exactly describe
the ansatz for the solution of problem (\ref{1.4}) and  rewrite  problem (\ref{1.4}) in terms of a linearized operator for which
a solvability theory, subject to suitable orthogonality conditions,
is performed  through solving a linearized  problem in Section $3$. In Section $4$ we
solve  a   nonlinear projected problem.
In Section $5$ we set up a maximization problem. In the last section we show that the solution to the
maximization problem indeed yields  a solution of problem (\ref{1.4})
with the qualitative properties  as predicted in
Theorem 1.1.

Throughout this paper, the symbol $C$ will always denote
a generic positive  constant independent of $t$,
it could be changed from one line to another.

\vspace{1mm}
\vspace{1mm}
\vspace{1mm}
\vspace{1mm}

\section{Ansatz for the solution}
In this section we will provide an ansatz for
solutions of   problem (\ref{1.4}).
For the sake of convenience we always fix the point
$p$ as an isolated local maximum point of $\phi_1$ in $\Omega$, and further
assume
\begin{equation}\label{2.1}
\aligned
\phi_1(p)=1.
\endaligned
\end{equation}
The configuration space for $m$ concentration points
$\xi=(\xi_1,\ldots,\xi_m)$ we try to seek is the following
\begin{eqnarray}\label{2.2}
\mathcal{O}_t:=\left\{\,\xi=(\xi_1,\ldots,\xi_m)\in\big(B_d(p)\big)^m\left|\,\,
|\xi_i-p|\geq\frac1{t^\beta},
\,\,\,\,\,
|\xi_i-\xi_j|\geq\frac1{t^\beta},
\,\,\,\,\,1-\phi_1(\xi_i)\leq\frac{1}{\sqrt{t}},
\right.\right.
&&\nonumber\\[0.5mm]
i,j=1,\ldots,m,\,\,\,i\neq j
\big\},
\qquad\qquad\qquad\qquad\qquad\qquad\,\,
\qquad\qquad\qquad\qquad\qquad\qquad\,
&&
\end{eqnarray}
where $d>0$ is a sufficiently small but fixed number, independent of $t$, and
$\beta$ is given by
\begin{equation}\label{2.3}
\aligned
\beta=\frac{\,(m+1)(m+1+\alpha)\,}{2}.
\endaligned
\end{equation}
Let us fix $\xi\in\mathcal{O}_t$. For numbers $\mu_0>0$ and
$\mu_i>0$, $i=1,\ldots,m$, yet to
be determined, we define
\begin{equation}\label{2.4}
\aligned
u_0(x)=\log
\frac{8\mu_0^2(1+\alpha)^2}{k(p)(\varepsilon_{0}^2\mu_{0}^2+|x-p|^{2(1+\alpha)})^2},
\qquad\qquad
u_i(x)=\log
\frac{8\mu_i^2}{k(\xi_i)|\xi_i-p|^{2\alpha}(\varepsilon_{i}^2\mu_{i}^2+|x-\xi_i|^2)^2},
\endaligned
\end{equation}
which satisfy in entire $\mathbb{R}^2$
\begin{equation}\label{2.5}
\aligned
-\Delta
u_0=\varepsilon_0^2k(p)|x-p|^{2\alpha}e^{u_0},
\,\qquad\quad\quad\,
-\Delta
u_i=\varepsilon_i^2k(\xi_i)|\xi_i-p|^{2\alpha}e^{u_i},
\endaligned
\end{equation}
having the properties
\begin{equation}\label{2.6}
\aligned
\int_{\mathbb{R}^2}\varepsilon_0^2k(p)|x-p|^{2\alpha}e^{u_0}=8\pi(1+\alpha),
\,\quad\quad\quad\quad\,
\int_{\mathbb{R}^2}\varepsilon_i^2k(\xi_i)|\xi_i-p|^{2\alpha}e^{u_i}=8\pi,
\endaligned
\end{equation}
where
\begin{equation}\label{2.7}
\aligned
\varepsilon_0=\varepsilon_0(t)\equiv e^{-\frac12t},
\,\qquad\qquad\,
\varepsilon_i=\varepsilon_i(t)\equiv e^{-\frac12t\phi_1(\xi_i)}.
\endaligned
\end{equation}

Our ansatz is then
\begin{equation}\label{2.8}
\aligned
U(x):=\sum\limits_{i=0}^{m}\,U_i(x)
=\sum\limits_{i=0}^{m}\,
\big[u_i(x)+H_i(x)\big],
\endaligned
\end{equation}
where $H_i(x)$ is a correction term defined as the solution of
\begin{equation}\label{2.9}
\left\{\aligned
&\Delta
H_i=0
\qquad
\textrm {in}\,\,\,\,\,\Omega,\\[1mm]
&H_i=-u_i
\,\,\,\quad
\textrm{on}\,\,\,\po.
\endaligned\right.
\end{equation}

\vspace{1mm}

\noindent{\bf Lemma 2.1.}\,\,{\it For any points  $\xi=(\xi_1,\ldots,\xi_m)\in\mathcal{O}_t$ and any $t$ large enough, then we have
\begin{equation}\label{2.10}
\aligned
H_0(x)=(1+\alpha)H(x,p)-\log
\frac{8\mu_0^2(1+\alpha)^2}{k(p)}+O\left(\varepsilon_0^2\mu_0^2\right),
\endaligned
\end{equation}
\begin{equation}\label{2.11}
\aligned
H_i(x)=H(x,\xi_i)-\log
\frac{8\mu_i^2}{k(\xi_i)|\xi_i-p|^{2\alpha}}+O\left(\varepsilon_i^2\mu_i^2\right),
\,\quad\,i=1,\ldots,m,
\endaligned
\end{equation}
uniformly in $\overline{\Omega}$, where $H$ is the regular part of Green's function defined in
{\upshape(\ref{1.3})}.
}

\vspace{1mm}

\begin{proof}
If we set $z(x)=H_0(x)-(1+\alpha)H(x,p)+\log\frac{8\mu_0^2(1+\alpha)^2}{k(p)}$,
then $z(x)$ is a harmonic function. Hence by (\ref{1.3}), (\ref{2.4}), (\ref{2.9}) and the maximum principle,
$$
\aligned
\max_{x\in\oo}\big|z(x)\big|&=\max_{x\in\po}\left|-u_0(x)-4(1+\alpha)\log|x-p|+\log\frac{8\mu_0^2(1+\alpha)^2}{k(p)}\right|\\
&=\max_{x\in\po}\left|\log\frac1{\,|x-p|^{4(1+\alpha)}\,}-\log
\frac{1}{(\varepsilon_0^2\mu_0^2+|x-p|^{2(\alpha+1)})^2}\right|=O\left(\varepsilon_0^2\mu_0^2\right),
\endaligned
$$
uniformly in $\oo$, as $s\rightarrow+\infty$, which implies that expansion (\ref{2.10}) holds.
Furthermore,  expansion (\ref{2.11}) can be also  obtained along these analogous arguments of (\ref{2.10}).
\end{proof}

\vspace{1mm}
\vspace{1mm}
\vspace{1mm}

Observe that $u_0$ and $u_i$, $i=1,\ldots,m$ are good approximations for a solution of problem
(\ref{1.4}) near points $p$  and $\xi_i$, $i=1,\ldots,m$, respectively.
We expect that the ansatz in (\ref{2.8}) is more accurate near $p$ and each $\xi_i$,
namely the remainders
$U-u_0=H_0+\sum_{j\neq 0}(u_j+H_j)$
and
$U-u_i=H_i+\sum_{j\neq i}(u_j+H_j)$ vanish at main order
near $p$ or $\xi_i$, respectively. This can be achieved through the following
precise choices of  the concentration parameters $\mu_0$ and $\mu_i$:
\begin{equation}\label{2.12}
\aligned
\log
\frac{8\mu_0^2(1+\alpha)^2}{k(p)}=
(1+\alpha)H(p,p)+\sum_{j=1}^m G(p,\xi_j),
\endaligned
\end{equation}
\begin{equation}\label{2.13}
\aligned
\log
\frac{8\mu_i^2}{k(\xi_i)|\xi_i-p|^{2\alpha}}=
H(\xi_i,\xi_i)+(1+\alpha)G(\xi_i,p)
+\sum_{j=1,\,j\neq
i}^m G(\xi_i,\xi_j),
\,\quad\,i=1,\ldots,m.
\endaligned
\end{equation}
We thus fix $\mu_0$ and $\mu_i$ {\it a priori} as  functions of $\xi$ in $\mathcal{O}_t$ and
write $\mu_0=\mu_0(\xi)$ and $\mu_i=\mu_i(\xi)$ for all $i=1,\ldots,m$.
Since $\xi=(\xi_1,\ldots,\xi_m)\in\mathcal{O}_t$,
there exists a constant $C>0$ independent of $t$ such that
\begin{equation}\label{2.14}
\aligned
\frac1C\leq\mu_0\leq Ct^{2m\beta}
\quad\quad\,\,\textrm{and}\quad\quad\,\,
\big|\partial_{\xi_{kl}}\log\mu_0\big|\leq Ct^{\beta},\,
\,\quad\,\forall\,\,k=1,\ldots,m,\,\,l=1,2,
\endaligned
\end{equation}
and
\begin{equation}\label{2.15}
\aligned
\frac1C\leq\mu_i\leq Ct^{(2m+\alpha)\beta}
\quad\quad\,\,\textrm{and}\quad\quad\,\,
\big|\partial_{\xi_{kl}}\log\mu_i\big|\leq Ct^{\beta},\,
\,\quad\,\forall\,\,i, k=1,\ldots,m,\,\,l=1,2.
\endaligned
\end{equation}

Consider now the change of variables
\begin{equation}\label{2.16}
\aligned
\omega(y)=u(\varepsilon_0 y)-2t\,\ \ \ \,\,\,\forall\,\,y\in\overline{\Omega}_t,
\endaligned
\end{equation}
with
$$
\aligned
\Omega_t=\varepsilon_0^{-1}\Omega=e^{\frac12t}\Omega,
\endaligned
$$
then $u(x)$ solves  equation (\ref{1.4}) if and only if $\omega(y)$ satisfies
\begin{equation}\label{2.17}
\aligned
\left\{\aligned
&-\Delta \omega=|\varepsilon_0y-p|^{2\alpha}q(y,t)e^\omega
\quad
\textrm{in}\,\,\,\,\,\Omega_t,\\[2mm]
&\omega=-2t
\qquad\qquad\qquad\quad\,\,\,
\quad\quad\,\textrm{on}
\,\,\,\partial\Omega_t,
\endaligned\right.
\endaligned
\end{equation}
where
\begin{equation}\label{2.18}
\aligned
q(y,t)\equiv k(\varepsilon_0 y)\exp
\left\{-t\big[\phi_1(\varepsilon_0 y)-1\big]\right\}.
\endaligned
\end{equation}
Let us write $p'=p/\varepsilon_0$ and $\xi'_i=\xi_i/\varepsilon_0$, $i=1,\ldots,m$, and define the initial approximate solution of (\ref{2.17}) as
\begin{equation}\label{2.19}
\aligned
V(y)=U(\varepsilon_0 y)-2t,
\endaligned
\end{equation}
where  $U$ is defined in (\ref{2.8}). Moreover, set
\begin{equation}\label{2.20}
\aligned
W(y)=|\varepsilon_0y-p|^{2\alpha}q(y,t)e^{V(y)},
\endaligned
\end{equation}
and the ``error term''
\begin{equation}\label{2.21}
\aligned
E(y)=\Delta V(y)+|\varepsilon_0y-p|^{2\alpha}q(y,t)e^{V(y)}.
\endaligned
\end{equation}

Let us see how well $-\Delta V(y)$ match with $W(y)$ through $V(y)$
so that the ``error term''
$E(y)$ is sufficiently small for any $y\in\overline{\Omega}_t$.
A simple computation shows that
\begin{eqnarray}\label{2.22}
-\Delta V(y)=-\varepsilon_0^2\Delta_{x}U(x)
=-\varepsilon_0^2\sum\limits_{i=0}^{m}\Delta_{x}U_i(x)
=-\varepsilon_0^2\sum\limits_{i=0}^{m}
\Delta_{x}\big[u_i(x)+H_i(x)\big]
&&\nonumber\\
=\varepsilon_0^4k(p)|x-p|^{2\alpha}e^{u_0}
+\varepsilon_0^2\sum\limits_{i=1}^{m}\varepsilon_i^2k(\xi_i)|\xi_i-p|^{2\alpha}e^{u_i}
\,\qquad\qquad\qquad\,\,&&\nonumber\\[0.5mm]
=\frac{8\varepsilon_0^4\mu_0^2(1+\alpha)^2|x-p|^{2\alpha}}{\,(\varepsilon_{0}^2\mu_{0}^2+|x-p|^{2(1+\alpha)})^2\,}
+\sum\limits_{i=1}^{m}
\frac{8\varepsilon_0^2\varepsilon_i^2\mu_i^2}{\,(\varepsilon_{i}^2\mu_{i}^2+|x-\xi_i|^2)^2\,}
\,\qquad\qquad\,\,\,\,
&&\nonumber\\[1mm]
=\left(\frac{\varepsilon_0}{\rho_0v_0}\right)^2\frac{8(1+\alpha)^2\big|\frac{\varepsilon_0 y-p}{\rho_0v_0}\big|^{2\alpha}}{\,\big(1+\big|\frac{\varepsilon_0 y-p}{\rho_0v_0}\big|^{2(1+\alpha)}\big)^2\,}+\sum\limits_{i=1}^{m}
\frac{1}{\gamma_{i}^2}\frac{8}{\big(1+\big|\frac{y-\xi'_i}{\gamma_{i}}\big|^2\big)^2},
\,\quad\quad\,
&&
\end{eqnarray}
where
\begin{equation}\label{2.23}
\aligned
\rho_0:=\varepsilon_0^{\frac{1}{1+\alpha}}=\exp\left\{
-\frac1{2(1+\alpha)}t\right\},
\,\,\,\quad\,\,\,\,
v_0:=\mu_0^{\frac{1}{1+\alpha}},
\,\,\,\quad\,\,\,
\gamma_i:=\frac{1}{\varepsilon_0}\varepsilon_i\mu_i=\mu_i\exp\left\{
-\frac12t\big[\phi_1(\xi_i)-1\big]\right\}.
\endaligned
\end{equation}
Then if $|y-\xi'_i|\leq1/(\varepsilon_0 t^{2\beta})$
for some index $i\in\{1,\ldots,m\}$,
\begin{equation}\label{2.24}
\aligned
-\Delta V(y)=\frac{1}{\gamma_{i}^2}\frac{8}{\big(1+\big|\frac{y-\xi'_i}{\gamma_{i}}\big|^2\big)^2}
+O\left(
\varepsilon_0^4\mu_0^2t^{(4+2\alpha)\beta}
\right)
+\sum\limits_{j=1,\,j\neq i}^{m}O\left(
\varepsilon_0^2\varepsilon_j^2\mu_j^2t^{4\beta}
\right),
\endaligned
\end{equation}
and if $|y-p'|\leq1/(\varepsilon_0 t^{2\beta})$,
\begin{equation}\label{2.25}
\aligned
-\Delta V(y)=\left(\frac{\varepsilon_0}{\rho_0v_0}\right)^2\frac{8(1+\alpha)^2\big|\frac{\varepsilon_0 y-p}{\rho_0v_0}\big|^{2\alpha}}{\,\big(1+\big|\frac{\varepsilon_0 y-p}{\rho_0v_0}\big|^{2(1+\alpha)}\big)^2\,}
+\sum\limits_{j=1}^{m}O\left(
\varepsilon_0^2\varepsilon_j^2\mu_j^2t^{4\beta}
\right),
\endaligned
\end{equation}
while if $|y-p'|>1/(\varepsilon_0 t^{2\beta})$ and
$|y-\xi'_i|>1/(\varepsilon_0 t^{2\beta})$
for all $i=1,\ldots,m$,
\begin{equation}\label{2.26}
\aligned
-\Delta V(y)=
O\left(
\varepsilon_0^4\mu_0^2t^{(8+4\alpha)\beta}
\right)
+\sum\limits_{j=1}^{m}O\left(
\varepsilon_0^2\varepsilon_j^2\mu_j^2t^{8\beta}
\right).
\endaligned
\end{equation}
On the other hand, let us  first fix the index
$i\in\{1,\ldots,m\}$ and the region $|y-\xi'_i|\leq1/(\varepsilon_0 t^{2\beta})$.
Then we have
\begin{eqnarray}\label{2.27}
W(y)=\varepsilon_0^4|\varepsilon_0y-p|^{2\alpha}q(y, t)\exp
\left\{
\sum\limits_{j=0}^{m}\,
\big[u_j(\varepsilon_0 y)+H_j(\varepsilon_0 y)\big]\right\}
\quad\qquad\qquad
\qquad\qquad\qquad
\quad\quad\,\,\,\,
&&\nonumber\\
=\varepsilon_0^4k(\varepsilon_0 y)|\varepsilon_0y-p|^{2\alpha}
\exp\left\{-t\big[\phi_1(\varepsilon_0 y)-1\big]
+\log\frac{8\mu_i^2}{k(\xi_i)|\xi_i-p|^{2\alpha}(\varepsilon_{i}^2\mu_{i}^2+\varepsilon_0^2|y-\xi'_i|^2)^2}
\right.
&&\nonumber\\
\left.
+H_i(\varepsilon_0 y)
+
\sum\limits_{j=0,\,j\neq i}^{m}\,
\big[u_j(\varepsilon_0 y)
+H_j(\varepsilon_0 y)\big]\right\}
\qquad\qquad\qquad\qquad\quad
\qquad\qquad\qquad\qquad\,\,\,\,
&&\nonumber\\
=\frac{1}{\gamma_{i}^2}\frac{8}{\big(1+\big|\frac{y-\xi'_i}{\gamma_{i}}\big|^2\big)^2}
\times\frac{k(\varepsilon_0 y)}{k(\xi_i)}
\times\frac{|\varepsilon_0 y-p|^{2\alpha}}{|\xi_i-p|^{2\alpha}}
\times\exp
\left\{-t\big[\phi_1(\varepsilon_0 y)-\phi_1(\xi_i)\big]
\right\}
\qquad\qquad\,\,\,\,\,
&&\nonumber\\
\times\exp\left\{H_i(\varepsilon_0 y)
+
\sum\limits_{j=0,\,j\neq i}^{m}\,
\big[u_j(\varepsilon_0 y)
+H_j(\varepsilon_0 y)\big]\right\}.
\qquad\qquad\qquad\quad\,\,
\qquad\qquad\qquad\quad\,
&&
\end{eqnarray}
%\begin{eqnarray}\label{}
%W(y)=\varepsilon_0^4|\varepsilon_0y-p|^{2\alpha}q(y, t)\exp
%\left\{
%\sum\limits_{j=0}^{m}\,
%\big[u_j(\varepsilon_0 y)+H_j(\varepsilon_0 y)\big]\right\}
%\quad\qquad\qquad\quad\quad\quad
%\qquad\qquad\qquad\quad\quad\quad
%&&\nonumber\\
%=\varepsilon_0^4k(\varepsilon_0 y)|\varepsilon_0y-p|^{2\alpha}
%\exp\left\{-t\big[\phi_1(\varepsilon_0 y)-1\big]
%+\left[\log
%\frac{8\mu_0^2(1+\alpha)^2}{k(p)(\varepsilon_{0}^2\mu_{0}^2+|\varepsilon_0y-p|^{2(\alpha+1)})^2}
%+H_0(\varepsilon_0y)\right]
%\right.
%&&\nonumber\\
%\left.+
%\sum\limits_{j=1}^{m}\,
%\left[\log\frac{8\mu_j^2}{k(\xi_j)|\xi_j-p|^{2\alpha}(\varepsilon_{j}^2\mu_{j}^2+\varepsilon_0^2|y-\xi'_j|^2)^2}
%+H_j(\varepsilon_0 y)\right]\right\}
%\qquad\qquad\qquad\qquad\qquad\qquad\quad
%&&\nonumber\\
%=\frac{1}{\gamma_{i}^2}\frac{8}{\big(1+\big|\frac{y-\xi'_i}{\gamma_{i}}\big|^2\big)^2}
%\times\frac{k(\varepsilon_0 y)}{k(\xi_i)}
%\times\frac{|\varepsilon_0 y-p|^{2\alpha}}{|\xi_i-p|^{2\alpha}}
%\times\exp
%\left\{-t\big[\phi_1(\varepsilon_0 y)-\phi_1(\xi_i)\big]
%\right\}
%\qquad\qquad\qquad\qquad\,
%&&\nonumber\\
%\times\exp\left\{H_i(\varepsilon_0 y)+\left[\log
%\frac{8\mu_0^2(1+\alpha)^2}{k(p)(\varepsilon_{0}^2\mu_{0}^2+|\varepsilon_0y-p|^{2(\alpha+1)})^2}
%+H_0(\varepsilon_0 y)\right]
%\right.
%\qquad\qquad\qquad\qquad\qquad\quad\,
%&&\nonumber\\
%\left.+
%\sum\limits_{j=1,\,j\neq i}^{m}\,
%\left[\log\frac{8\mu_j^2}{k(\xi_j)|\xi_j-p|^{2\alpha}(\varepsilon_{j}^2\mu_{j}^2+\varepsilon_0^2|y-\xi'_j|^2)^2}+H_j(\varepsilon_0 y)\right]\right\},
%\qquad\qquad\qquad\qquad\qquad\quad&&
%\end{eqnarray}
From (\ref{2.4}), (\ref{2.10}), (\ref{2.11}) and the fact that $H(\cdot,x)$ is $C^1(\overline{\Omega})$ for any $x\in\Omega$,
we have that for $|y-\xi'_i|\leq1/(\varepsilon_0 t^{2\beta})$,
$$
\aligned
H_i&(\varepsilon_0 y)
+
\sum\limits_{j=0,\,j\neq i}^{m}\,
\big[u_j(\varepsilon_0 y)
+H_j(\varepsilon_0 y)\big]\\
&=H(\xi_i,\xi_i)-\log
\frac{8\mu_i^2}{k(\xi_i)|\xi_i-p|^{2\alpha}}+
(1+\alpha)\left[\log\frac{1}{|\xi_i-p|^{4}}+H(\xi_i,p)
+O\left(
\frac{|\varepsilon_0 y-\xi_i|}{|\xi_i-p|}
+\frac{\varepsilon_0^2\mu_0^2}{|\xi_i-p|^{2(1+\alpha)}}
\right)
\right]\\
&\quad\,+
\sum\limits_{j=1,\,j\neq i}^{m}\,
\left[\log\frac{1}{|\xi_i-\xi_j|^4}+H(\xi_i,\xi_j)
+O\left(
\frac{|\varepsilon_0 y-\xi_i|}{|\xi_i-\xi_j|}
+
\frac{\varepsilon_j^2\mu_j^2}{|\xi_i-\xi_j|^2}
\right)
\right]+O\left(\varepsilon_0|y-\xi'_i|\right)+\sum_{j=0}^mO\left(\varepsilon_j^2\mu_j^2\right)\\
&=H(\xi_i,\xi_i)-\log
\frac{8\mu_i^2}{k(\xi_i)|\xi_i-p|^{2\alpha}}
+(1+\alpha)G(\xi_i,p)+
\sum\limits_{j=1,\,j\neq i}^{m}\,
G(\xi_i,\xi_j)
+
O\left(\varepsilon_0t^{\beta}|y-\xi'_i|\right)\\
&\quad\,
+
O\left(\varepsilon_i^2\mu_i^2\right)
+
O\left(\varepsilon_0^2\mu_0^2t^{2\beta(1+\alpha)}\right)
+\sum_{j=1,\,j\neq i}^mO\left(\varepsilon_j^2\mu_j^2t^{2\beta}\right)\\
&=O\left(\varepsilon_0t^{\beta}|y-\xi'_i|\right)
+
O\left(\varepsilon_i^2\mu_i^2\right)
+
O\left(\varepsilon_0^2\mu_0^2t^{2\beta(1+\alpha)}\right)
+
\sum_{j=1,\,j\neq i}^mO\left(\varepsilon_j^2\mu_j^2t^{2\beta}\right),
\endaligned
$$
where the last equality is due to the choice of $\mu_i$
in (\ref{2.13}). Thus if
$|y-\xi'_i|\leq1/(\varepsilon_0 t^{2\beta})$ for some $i\in\{1,\ldots,m\}$,
\begin{eqnarray}\label{2.28}
W(y)=\frac{1}{\gamma_{i}^2}\frac{8}{\big(1+\big|\frac{y-\xi'_i}{\gamma_{i}}\big|^2\big)^2}
\left\{1+O\big(\varepsilon_0t^{\beta}|y-\xi'_i|\big)
+
O\big(\varepsilon_i^2\mu_i^2\big)
+
O\big(\varepsilon_0^2\mu_0^2t^{2\beta(1+\alpha)}\big)
+\sum_{j=1,\,j\neq i}^m O\big(\varepsilon_j^2\mu_j^2t^{2\beta}\big)\right\},
\end{eqnarray}
and  by (\ref{2.24}),
\begin{eqnarray}\label{2.29}
E(y)=\frac{1}{\gamma_{i}^2}\frac{8}{\big(1+\big|\frac{y-\xi'_i}{\gamma_{i}}\big|^2\big)^2}
\left\{O\left(\varepsilon_0t^{\beta}|y-\xi'_i|\right)
+
O\left(\varepsilon_i^2\mu_i^2\right)
+
O\left(\varepsilon_0^2\mu_0^2t^{2\beta(1+\alpha)}\right)
+\sum_{j=1,\,j\neq i}^mO\left(\varepsilon_j^2\mu_j^2t^{2\beta}\right)
\right\}
&&\nonumber\\
+\,O\left(
\varepsilon_0^4\mu_0^2t^{(4+2\alpha)\beta}
\right)
+\sum\limits_{j=1,\,j\neq i}^{m}O\left(
\varepsilon_0^2\varepsilon_j^2\mu_j^2t^{4\beta}
\right).
\qquad\qquad\qquad\qquad\qquad\qquad
\qquad\qquad\qquad\qquad\,
&&
\end{eqnarray}
Similarly,  if   $|y-p'|\leq1/(\varepsilon_0 t^{2\beta})$, by  (\ref{2.1}), (\ref{2.4}), (\ref{2.10}), (\ref{2.11}) and (\ref{2.12}) we can compute
\begin{equation}\label{2.30}
\aligned
W(y)
=\left(\frac{\varepsilon_0}{\rho_0v_0}\right)^2\frac{8(1+\alpha)^2\big|\frac{\varepsilon_0 y-p}{\rho_0v_0}\big|^{2\alpha}}{\,\big(1+\big|\frac{\varepsilon_0 y-p}{\rho_0v_0}\big|^{2(1+\alpha)}\big)^2\,}
\left[1+O\left(\varepsilon_0t^{\beta}|y-p'|\right)
+
O\left(\varepsilon_0^2\mu_0^2\right)
+\sum_{j=1}^mO\left(\varepsilon_j^2\mu_j^2t^{2\beta}\right)\right],
\endaligned
\end{equation}
and  by  (\ref{2.25}),
\begin{eqnarray}\label{2.31}
E(y)=\left(\frac{\varepsilon_0}{\rho_0v_0}\right)^2\frac{8(1+\alpha)^2\big|\frac{\varepsilon_0 y-p}{\rho_0v_0}\big|^{2\alpha}}{\,\big(1+\big|\frac{\varepsilon_0 y-p}{\rho_0v_0}\big|^{2(1+\alpha)}\big)^2\,}
\left[O\left(\varepsilon_0t^{\beta}|y-p'|\right)
+
O\left(\varepsilon_0^2\mu_0^2\right)
+\sum_{j=1}^mO\left(\varepsilon_j^2\mu_j^2t^{2\beta}\right)\right]
&&\nonumber\\
+\sum\limits_{j=1}^{m}O\left(
\varepsilon_0^2\varepsilon_j^2\mu_j^2t^{4\beta}
\right).
\qquad\qquad\qquad\qquad
\qquad\qquad\qquad\qquad\qquad
\qquad\qquad\qquad\quad\,\,\,
&&
\end{eqnarray}
While if  $|y-p'|>1/(\varepsilon_0 t^{2\beta})$ and
$|y-\xi'_i|>1/(\varepsilon_0 t^{2\beta})$
for all $i=1,\ldots,m$,
by (\ref{2.4}), (\ref{2.10}) and (\ref{2.11}) we obtain
\begin{equation}\label{2.32}
\aligned
W(y)
=O\left(\frac{\varepsilon_0^2e^{-t\phi_1(\varepsilon_0 y)}}
{\,|\varepsilon_0 y-p|^{4+2\alpha}\,}\prod_{i=1}^m\frac{1}{\,|\varepsilon_0 y-\xi_i|^4\,}
\right).
\endaligned
\end{equation}
Then by (\ref{2.26}),
\begin{equation}\label{2.33}
\aligned
E(y)=O\left(\frac{\varepsilon_0^2e^{-t\phi_1(\varepsilon_0 y)}}{\,|\varepsilon_0 y-p|^{4+2\alpha}\,}\prod_{i=1}^m\frac{1}{\,|\varepsilon_0 y-\xi_i|^4\,}
\right)+
O\left(
\varepsilon_0^4\mu_0^2t^{(8+4\alpha)\beta}
\right)
+\sum\limits_{j=1}^{m}O\left(
\varepsilon_0^2\varepsilon_j^2\mu_j^2t^{8\beta}
\right).
\endaligned
\end{equation}

In what remains of this paper we will seek
solutions  of problem (\ref{2.17}) in the form
$\omega=V+\phi$,  where $\phi$ will represent a lower order correction.
In terms of $\phi$, problem (\ref{2.17}) becomes
\begin{equation}\label{2.34}
\left\{\aligned
&\mathcal{L}(\phi):=-\Delta\phi-W\phi=E+N(\phi)
\,\,\,\,
\textrm{in}\,\,\,\,\,\Omega_t,\\[2mm]
&\phi=0
\quad\qquad\qquad\qquad\,
\qquad\qquad\quad\,\,
\textrm{on}\,\,\,
\partial\Omega_t,
\endaligned\right.
\end{equation}
where the ``nonlinear term''  is given by
\begin{equation}\label{2.35}
\aligned
N(\phi)=W\big(e^\phi-1-\phi\big).
\endaligned
\end{equation}

\vspace{1mm}
\vspace{1mm}
\vspace{1mm}

\section{Solvability of a linear problem}
In this section we consider  the solvability of
the following linear problems: given $h\in L^{\infty}(\Omega_t)$
and points $\xi=(\xi_1,\ldots,\xi_m)\in\mathcal{O}_t$,
we find a function $\phi$ such that for certain  scalars $c_{ij}$,
$i=1,\ldots,m$, $j=1,2$, one has
\begin{equation}\label{3.1}
\aligned
\left\{
\aligned
&\mathcal{L}(\phi)=-\Delta\phi-W\phi=h+\sum\limits_{i=1}^m\sum\limits_{j=1}^2c_{ij}\chi_iZ_{ij}
\quad\,\,\,\,
\textrm{in}\,\,\,\,\,\Omega_t,\\
&\phi=0
\quad\quad\quad\quad\quad
\,\quad\quad\,\,\,\,\,
\qquad\qquad\qquad\quad
\quad\quad\,\,
\textrm{on}
\,\,\,\po_t,\\[1mm]
&\int_{\Omega_t}\chi_iZ_{ij}\phi=0
\qquad\qquad\qquad
\forall\,\,i=1,\ldots,m,
\,\,j=1,2,
\endaligned
\right.
\endaligned
\end{equation}
where $W=|\varepsilon_0y-p|^{2\alpha}q(y,t)e^{V}$
satisfies (\ref{2.28}), (\ref{2.30}) and (\ref{2.32}), and $Z_{ij}$, $\chi_i$ are defined as follows. Let
$\mathcal{Z}_p$, $\mathcal{Z}_0$, $\mathcal{Z}_1$ and $\mathcal{Z}_2$ be
\begin{equation}\label{3.2}
\aligned
\mathcal{Z}_{p}(z)=\frac{|z|^{2(1+\alpha)}-1}{|z|^{2(1+\alpha)}+1},
\,\,\quad\quad\quad\,\,
\mathcal{Z}_{0}(z)=\frac{|z|^2-1}{|z|^2+1},
\,\,\quad\quad\quad\,\,
\mathcal{Z}_{j}(z)=\frac{4z_j}{|z|^2+1},\,\,\,\,j=1,\,2.
\endaligned
\end{equation}
It is well known  that
\begin{itemize}
  \item any bounded
solution to
\begin{equation}\label{3.3}
\aligned
\Delta
\phi+\frac{8(1+\alpha)^2|z|^{2\alpha}}{(1+|z|^{2(1+\alpha)})^2}\phi=0
\,\quad\textrm{in}\,\,\,\mathbb{R}^2,
\endaligned
\end{equation}
where $-1<\alpha\not\in\mathbb{N}$,
is proportional to $\mathcal{Z}_p$ (see \cite{E,CY,EPW});
  \item any bounded
solution to
\begin{equation}\label{3.4}
\aligned
\Delta
\phi+\frac{8}{(1+|z|^2)^2}\phi=0
\,\quad\textrm{in}\,\,\,\mathbb{R}^2,
\endaligned
\end{equation}
is a linear combination of $\mathcal{Z}_j$, $j=0,1,2$ (see \cite{BP,CL}).
\end{itemize}
Then we define
\begin{equation}\label{3.5}
\aligned
Z_{p}(y)=\frac{\varepsilon_0}{\rho_0 v_0}\mathcal{Z}_p\left(\frac{\varepsilon_0y-p}{\rho_0 v_0}\right)
\,\,\quad\,\,\,\textrm{and}
\,\,\quad\,\,\,
Z_{ij}(y)=\frac{1}{\gamma_i}\mathcal{Z}_j\left(\frac{y-\xi_i'}{\gamma_i}\right),
\,\,\,\,\,\,i=1,\ldots,m,\,\,j=0,1,2.
\endaligned
\end{equation}

Next, we consider  a large but fixed positive number $R_0$ and set
a radial, smooth non-increasing  cut-off function $\chi(r)$ with
$0\leq\chi(r)\leq1$, $\chi(r)=1$ for $r\leq R_0$ and $\chi(r)=0$ for $r\geq R_0+1$.
Let
\begin{equation}\label{3.6}
\aligned
\chi_p(y)=\chi\left(\frac{\big|\varepsilon_0y-p\big|}{\rho_0 v_0}\right)
\,\,\,\quad\,\,\,
\textrm{and}
\,\,\,\quad\,\,\,
\chi_i(y)=\chi\left(\frac{\big|y-\xi_i'\big|}{\gamma_i}\right),\,
\quad\,i=1,\ldots,m.
\endaligned
\end{equation}

Equation (\ref{3.1}) will be solved for
$h\in L^\infty(\Omega_t)$, but we need to estimate the size of
the solution in terms of the following $L^\infty$-weighted norm:
\begin{eqnarray}\label{3.7}
\quad
\|h\|_{*}:=
\left\|\left[
\varepsilon_0^2+
\left(\frac{\varepsilon_0}{\rho_0v_0}\right)^2\frac{\big|\frac{\varepsilon_0 y-p}{\rho_0v_0}\big|^{2\alpha}}
{\,\big(1+\big|\frac{\varepsilon_0 y-p}{\rho_0v_0}\big|\big)^{4+2\hat{\alpha}+2\alpha}\,}
+
\sum\limits_{i=1}^m\frac{1}{\gamma_i^2}\frac{1}{\big(1+\big|\frac{y-\xi'_i}{\gamma_i}\big|\big)^{4+2\hat{\alpha}}}
\right]^{-1}h(y)
\right\|_{L^\infty(\Omega_t)},
\end{eqnarray}
%\begin{equation}\label{3.7}
%\aligned
%\|h\|_{*}=
%\sup_{y\in\Omega_t}\left|\left[
%\varepsilon_0^2+
%\left(\frac{\varepsilon_0}{\rho_0v_0}\right)^2\frac{\big|\frac{\varepsilon_0 y-p}{\rho_0v_0}\big|^{2\alpha}}
%{\,\big(1+\big|\frac{\varepsilon_0 y-p}{\rho_0v_0}\big|\big)^{4+2\hat{\alpha}+2\alpha}\,}
%+
%\sum\limits_{i=1}^m\frac{1}{\gamma_i^2}\frac{1}{\big(1+\big|\frac{y-\xi'_i}{\gamma_i}\big|\big)^3}
%\right]^{-1}h(y)
%\right|,
%\endaligned
%\end{equation}
where $\hat{\alpha}+1$ is a sufficiently small but fixed positive number, independent of $t$,
such that $-1<\hat{\alpha}<\min\big\{\alpha,\,-2/3\big\}$.

\vspace{1mm}
\vspace{1mm}
\vspace{1mm}
\vspace{1mm}

\noindent{\bf Proposition 3.1.}\,\,{\it
Let $m$ be a positive  integer. Then
there exist constants $t_m>1$ and $C>0$ such
that for any $t>t_m$,
any points $\xi=(\xi_1,\ldots,\xi_m)\in\mathcal{O}_t$ and any
$h\in L^{\infty}(\Omega_t)$, there is a unique solution
$\phi\in L^{\infty}(\Omega_t)$,
$c_{ij}\in\mathbb{R}$,
$i=1,\ldots,m$, $j=1,2$
to problem {\upshape(\ref{3.1})}. Moreover
\begin{equation}\label{3.8}
\aligned
\|\phi\|_{L^{\infty}(\Omega_t)}\leq Ct\|h\|_{*}.
\endaligned
\end{equation}

}
\begin{proof}
The proof of this result consists of  five steps which we state and prove next.

\vspace{1mm}

{\bf Step 1:} The operator $\mathcal{L}$ satisfies the maximum principle in
$\widetilde{\Omega}_t:=
\Omega_t\setminus\big[\bigcup_{i=1}^mB_{R_1\gamma_i}(\xi'_i)
\cup B_{R_1\rho_0 v_0/\varepsilon_0}(p')\cup B^c_{2d/\varepsilon_0}(p')\big]$
for $R_1$ large but $d$ small independent of $p$,
namely if $\psi$ satisfies $\mathcal{L}(\psi)=-
\Delta\psi-W\psi\geq0$ in $\widetilde{\Omega}_t$,
$\psi\geq0$ on $\partial\widetilde{\Omega}_t$,
then $\psi\geq0$ in $\widetilde{\Omega}_t$.
In order to prove it,  we shall  first  find a  function $\mathcal{Z}$
 such that $\mathcal{L}(\mathcal{Z})>0$ and $\mathcal{Z}>0$  in $\widetilde{\Omega}_t$. Indeed, let
\begin{equation}\label{3.9}
\aligned
\mathcal{Z}(y)=
\Psi_0(\varepsilon_0y)
+
\mathcal{\widehat{Z}}_p\left(
a\left|\frac{\varepsilon_0 y-p}{\rho_0v_0}\right|
\right)
+
\sum_{i=1}^m\mathcal{Z}_0\left(
a\left|\frac{y-\xi'_i}{\gamma_i}\right|
\right),
\endaligned
\end{equation}
where $a>0$, $\Psi_0$ satisfies
$-\Delta\Psi_0=1$ in $\Omega$, $\Psi_0=2$ on $\partial\Omega$,
$\mathcal{Z}_0$ is defined in (\ref{3.2}) and
$$
\aligned
\mathcal{\widehat{Z}}_p(z)=
\frac{|z|^{2(1+\hat{\alpha})}-1}{|z|^{2(1+\hat{\alpha})}+1}.
\endaligned
$$
Observe that
$$
\aligned
-\Delta \mathcal{Z}=\varepsilon_0^2+\left(\frac{\varepsilon_0}{\rho_0v_0}\right)^2\frac{8a^2(1+\hat{\alpha})^2\big|a\frac{\varepsilon_0 y-p}{\rho_0v_0}\big|^{2\hat{\alpha}}}
{\big(1+\big|a\frac{\varepsilon_0 y-p}{\rho_0v_0}
\big|^{2(1+\hat{\alpha})}\big)^2}\mathcal{\widehat{Z}}_p\left(a\left|
\frac{\varepsilon_0 y-p}{\rho_0v_0}
\right|\right)
+\sum_{i=1}^m\frac{1}{\gamma_i^2}\frac{8a^2}{\big(1+\big|
a\frac{y-\xi'_i}{\gamma_i}
\big|^2\big)^2}\mathcal{Z}_0\left(
a\left|\frac{y-\xi'_i}{\gamma_i}\right|
\right).
\endaligned
$$
Then if
$3^{1/(2+2\hat{\alpha})}\rho_0v_0/(a\varepsilon_0)\leq\big|y-p'\big|\leq1/(\varepsilon_0 t^{2\beta})$,
by (\ref{2.30}),
$$
\aligned
\mathcal{L}(\mathcal{Z})&\geq
\left(\frac{\varepsilon_0}{\rho_0v_0}\right)^2a^2\frac{4(1+\hat{\alpha})^2\big|a\frac{\varepsilon_0 y-p}{\rho_0v_0}\big|^{2\hat{\alpha}}}
{\big(1+\big|a\frac{\varepsilon_0 y-p}{\rho_0v_0}
\big|^{2(1+\hat{\alpha})}\big)^2}
-W\|\mathcal{Z}\|_{\infty}\\[1mm]
&\geq
\left(\frac{\varepsilon_0}{\rho_0v_0}\right)^2\frac{(1+\hat{\alpha})^2}
{a^{2(1+\hat{\alpha})}\big|\frac{\varepsilon_0 y-p}{\rho_0v_0}
\big|^{4+2\hat{\alpha}}}
-
\left(\frac{\varepsilon_0}{\rho_0v_0}\right)^2\frac{8C(1+\alpha)^2\big|\frac{\varepsilon_0 y-p}{\rho_0v_0}\big|^{2\alpha}}
{\,\big(1+\big|\frac{\varepsilon_0 y-p}{\rho_0v_0}\big|^{2(1+\alpha)}\big)^2\,}\\[2mm]
%&\geq
%\left(\frac{\varepsilon_0}{\rho_0v_0}\right)^2\frac{(1+\hat{\alpha})^2}
%{a^{2(1+\hat{\alpha})}\big|\frac{\varepsilon_0 y-p}{\rho_0v_0}
%\big|^{4+2\hat{\alpha}}}
%-
%\left(\frac{\varepsilon_0}{\rho_0v_0}\right)^2\frac{8C(m+1)(1+\alpha)^2
%}{\,\big|\frac{\varepsilon_0 y-p}{\rho_0v_0}\big|^{4+2\alpha}\,}\\
&\geq
\left(\frac{\varepsilon_0}{\rho_0v_0}\right)^2\frac{1}
{\,\big|\frac{\varepsilon_0 y-p}{\rho_0v_0}
\big|^{4+2\hat{\alpha}}\,}\left[
\frac{(1+\hat{\alpha})^2}
{a^{2(1+\hat{\alpha})}}
-
\frac{8C(1+\alpha)^2
}{\,\big|\frac{\varepsilon_0 y-p}{\rho_0v_0}\big|^{2(\alpha-\hat{\alpha})}\,}
\right].
\endaligned
$$
Similarly,  if
$3^{1/2}\gamma_i/a\leq\big|y-\xi'_i\big|\leq1/(\varepsilon_0 t^{2\beta})$
for some index $i\in\{1,\ldots,m\}$, by (\ref{2.28}),
$$
\aligned
\mathcal{L}(\mathcal{Z})\geq
\frac{1}{\gamma_i^2}\frac{1}{\,\big|
\frac{y-\xi'_i}{\gamma_i}
\big|^4\,}\left(\,\frac{1}{\,a^2\,}
-
8C
\right).
\endaligned
$$
While if  $1/(\varepsilon_0 t^{2\beta})<|y-p'|<2d/\varepsilon_0$ and
$|y-\xi'_i|>1/(\varepsilon_0 t^{2\beta})$
for all $i=1,\ldots,m$,
by  (\ref{2.32}),
$$
\aligned
\mathcal{L}(\mathcal{Z})
&\geq
\varepsilon_0^2+
\left(\frac{\varepsilon_0}{\rho_0v_0}\right)^2\frac{(1+\hat{\alpha})^2}
{a^{2(1+\hat{\alpha})}\big|\frac{\varepsilon_0 y-p}{\rho_0v_0}
\big|^{4+2\hat{\alpha}}}
+\sum_{i=1}^m\frac{1}{\gamma_i^2}\frac{1}{a^2\big|
\frac{y-\xi'_i}{\gamma_i}
\big|^4}-
C\frac{\varepsilon_0^2e^{-t\phi_1(\varepsilon_0 y)}}{|\varepsilon_0 y-p|^{4+2\alpha}}\prod_{i=1}^m\frac{1}{|\varepsilon_0 y-\xi_i|^4}.
\endaligned
$$
Hence if $a$  is taken sufficiently small but fixed, and $R_1$ is chosen sufficiently large depending
on this $a$, then by (\ref{2.23}) we can easily conclude that
$\mathcal{L}(\mathcal{Z})>0$  in
$\widetilde{\Omega}_t$.

Next, we suppose that the operator $\mathcal{L}$ does not satisfy the maximum principle in
$\widetilde{\Omega}_t$. Since $\mathcal{Z}>0$ in $\widetilde{\Omega}_t$, it follows that
the function $\psi/\mathcal{Z}$ has a  negative minimum  point $y_0$
in $\widetilde{\Omega}_t$. A direct computation gives
$$
\aligned
-\Delta\left(\frac{\psi}{\mathcal{Z}}\right)=\frac{1}{\mathcal{Z}^2}\big[
\mathcal{Z}\mathcal{L}(\psi)
-\psi\mathcal{L}(\mathcal{Z})
\big]+\frac{2}{\mathcal{Z}}\nabla \mathcal{Z}\nabla\left(\frac{\psi}{\mathcal{Z}}\right).
\endaligned
$$
Then $-\Delta\big(\psi/\mathcal{Z}\big)(y_0)>0$, which contradicts to the fact that
$y_0$ is a minimum  point of $\psi/\mathcal{Z}$ in $\widetilde{\Omega}_t$.

\vspace{1mm}
\vspace{1mm}
\vspace{1mm}

{\bf Step 2:} Let $R_1$  be as before.
Since $\xi=(\xi_1,\ldots,\xi_m)\in\mathcal{O}_t$,
$\rho_0 v_0=o(1/t^{2\beta})$ and
$\varepsilon_0\gamma_i=o(1/t^{2\beta})$ for $t$ large enough,
we find $B_{R_1\rho_0 v_0/\varepsilon_0}(p')$ and
$B_{R_1\gamma_i}(\xi'_i)$, $i=1,\ldots,m$, disjointed and included in $\Omega_t$.
Let us  consider  the following norm
\begin{equation}\label{3.10}
\aligned
\|\phi\|_{**}=\sup_{y\in\bigcup_{i=1}^mB_{R_1\gamma_i}(\xi'_i)\cup B_{R_1\rho_0 v_0/\varepsilon_0}(p')
\cup\big(\Omega_t\setminus B_{2d/\varepsilon_0}(p')\big)}
|\phi|(y).
\endaligned
\end{equation}
We claim that there is a constant $C>0$ independent of  $t$ such that,
if $\phi$ is the solution of the linear equation
\begin{equation}\label{3.11}
\aligned
\left\{\aligned
&\mathcal{L}(\phi)=-\Delta\phi-W\phi=h
\,\,\quad\,\textrm{in}\,\,\,\,\,\Omega_t,\\[2mm]
&\phi=0
\,\,\,\,\quad\,
\,\,\,\,\quad\,
\,\quad\,\,\,\,
\,\,\,\qquad\,
\quad\,\textrm{on}
\,\,\,\po_t,
\endaligned\right.
\endaligned
\end{equation}
then
\begin{equation}\label{3.12}
\aligned
\|\phi\|_{L^{\infty}(\Omega_t)}\leq C\left(\|\phi\|_{**}+
\|h\|_{*}\right),
\endaligned
\end{equation}
for any $h\in L^\infty(\Omega_t)$ and any points $\xi=(\xi_1,\ldots,\xi_m)\in\mathcal{O}_t$.
We will establish this estimate with the use of suitable barriers.
Let $M$ be a large number such that
$\Omega\subset B(p,M)$ and $\Omega\subset B(\xi_i,M)$ for all $i=1,\ldots,m$. Consider
$\psi_0$ and $\psi_i$, $i=1,\ldots,m$, respectively, as the solutions of the
problems
$$
\begin{aligned}
\left\{\begin{aligned}
&-\Delta\psi_0=\left(\frac{\varepsilon_0}{\rho_0v_0}\right)^2
\frac{4}{\big|\frac{\varepsilon_0 y-p}{\rho_0v_0}\big|^{4+2\hat{\alpha}}}+4\varepsilon_0^2,
\quad\quad
R_1<\left|\frac{\varepsilon_0 y-p}{\rho_0v_0}\right|
<\frac{M}{\rho_0v_0},\\[1mm]
&\psi_0(y)=0,\ \ \ \ \ \ \ \ \ \,\textrm{for}\quad\left|\frac{\varepsilon_0 y-p}{\rho_0v_0}\right|=R_1,\,\,
\ \ \ \,
\qquad\quad
\left|\frac{\varepsilon_0 y-p}{\rho_0v_0}\right|=\frac{M}{\rho_0v_0}.
\end{aligned}\right.
\end{aligned}
$$
and
$$
\left\{\begin{aligned}
&-\Delta\psi_i=\frac{1}{\gamma_i^2}
\frac{4}{\big|\frac{y-\xi'_i}{\gamma_i}\big|^{4+2\hat{\alpha}}}+4\varepsilon_0^2,
\,\,\quad\quad\,\,
R_1<\left|\frac{y-\xi'_i}{\gamma_i}\right|
<\frac{M}{\varepsilon_0\gamma_i},\\[1mm]
&\psi_i(y)=0,\,\ \ \ \ \
\ \,\textrm{for}\quad\left|\frac{y-\xi'_i}{\gamma_i}\right|=R_1,
\,\,\ \ \ \ \,\,
\quad
\left|\frac{y-\xi'_i}{\gamma_i}\right|=\frac{M}{\varepsilon_0\gamma_i}.
\end{aligned}\right.
$$
Then the solutions $\psi_0$ and $\psi_i$, $i=1,\ldots,m$, are
the positive functions, respectively  given by
$$
\aligned
\psi_0(y)=&-\left[\frac1{(1+\hat{\alpha})^2} \frac1{\big|\frac{\varepsilon_0 y-p}{\rho_0v_0}\big|^{2(1+\hat{\alpha})}}
+\rho^2_0v^2_0\left|\frac{\varepsilon_0 y-p}{\rho_0v_0}\right|^2\right]
+\left[
\frac1{(1+\hat{\alpha})^2} \frac1{R_1^{2
(1+\hat{\alpha})}}+R_1^2\rho^2_0v^2_0\right]\\[1mm]
&-
\frac1{(1+\hat{\alpha})^2}\left(\frac1{R_1^{2
(1+\hat{\alpha})}}-\frac{(\rho_0v_0)^{2(1+\hat{\alpha})}}{M^{2(1+\hat{\alpha})}}\right)
\frac{\,\log\big|\frac{\varepsilon_0 y-p}{R_1\rho_0v_0}\big|\,}{\,\log\frac{M}{R_1\rho_0v_0}\,}
+\big(M^2-R_1^2\rho^2_0v^2_0\big)
\frac{\,\log\big|\frac{\varepsilon_0 y-p}{R_1\rho_0v_0}\big|\,}{\,\log\frac{M}{R_1\rho_0v_0}\,},
\endaligned
$$
$$
\aligned
\psi_i(y)=&-\left[\frac1{(1+\hat{\alpha})^2} \frac1{\big|\frac{y-\xi'_i}{\gamma_i}\big|^{2(1+\hat{\alpha})}}
+\varepsilon^2_0\gamma^2_i\left|\frac{y-\xi'_i}{\gamma_i}\right|^2\right]
+\left[
\frac1{(1+\hat{\alpha})^2} \frac1{R_1^{2
(1+\hat{\alpha})}}+R_1^2\varepsilon^2_0\gamma^2_i\right]\\[1mm]
&-
\frac1{(1+\hat{\alpha})^2}\left(\frac1{R_1^{2
(1+\hat{\alpha})}}-\frac{(\varepsilon_0\gamma_i)^{2(1+\hat{\alpha})}}{M^{2(1+\hat{\alpha})}}\right)
\frac{\,\log\left|\frac{y-\xi'_i}{R_1\gamma_i}\right|\,}{\,\log\frac{M}{\,R_1\varepsilon_0\gamma_i}\,}
+\big(M^2-R_1^2\varepsilon^2_0\gamma^2_i\big)
\frac{\,\log\left|\frac{y-\xi'_i}{R_1\gamma_i}\right|\,}{\,\log\frac{M}{\,R_1\varepsilon_0\gamma_i}\,},
\quad
i=1,\ldots,m.
\endaligned
$$
%$$
%\aligned
%\psi_0(y)=\varphi_0\left(\left|\frac{\varepsilon_0 y-p}{\rho_0v_0}\right|\right)
%-
%\varphi_0\left(\frac{M}{\rho_0v_0}\right)
%\frac{\,\log\big|\frac{\varepsilon_0 y-p}{R_1\rho_0v_0}\big|\,}{\,\log\frac{M}{R_1\rho_0v_0}\,},
%\endaligned
%$$
%$$
%\aligned
%\psi_i(y)=\varphi_i\left(\left|\frac{y-\xi'_i}{\gamma_i}\right|\right)-\varphi_i\left(\frac{M}{\varepsilon_0\gamma_i}\right)
%\frac{\,\log\left|\frac{y-\xi'_i}{R_1\gamma_i}\right|\,}{\,\log\frac{M}{\,R_1\varepsilon_0\gamma_i}\,},
%\,\qquad\,
%i=1,\ldots,m,
%\endaligned
%$$
%where
%$$
%\aligned
%\varphi_0(r)=\frac1{(1+\hat{\alpha})^2} \left(\frac1{R_1^{2
%(1+\hat{\alpha})}}-\frac1{r^{2(1+\hat{\alpha})}}\right)+\rho^2_0v^2_0(R_1^2-r^2),
%\endaligned
%$$
%$$
%\aligned
%\varphi_i(r)=\frac1{(1+\hat{\alpha})^2} \left(\frac1{R_1^{2
%(1+\hat{\alpha})}}-\frac1{r^{2(1+\hat{\alpha})}}\right)+\varepsilon^2_0\gamma^2_i(R_1^2-r^2),
%\,\qquad\,
%i=1,\ldots,m.
%\endaligned
%$$
Clearly,  the functions $\psi_0$ and $\psi_i$, $i=1,\ldots,m$, are
uniformly bounded from above by a constant independent of $t$. Let us  consider the function $\mathcal{Z}(y)$ defined in the previous step.
We  take the barrier
\begin{equation}\label{3.13}
\aligned
\widetilde{\phi}(y)=2\|\phi\|_{**}\mathcal{Z}(y)+
\|h\|_{*}\sum_{i=0}^m\psi_i(y).
\endaligned
\end{equation}
Choosing $R_1$ larger if necessary, we have that
for $y\in\bigcup_{i=1}^m\partial B_{R_1\gamma_i}(\xi'_i)
\cup\partial B_{R_1\rho_0 v_0/\varepsilon_0}(p')\cup\partial B_{2d/\varepsilon_0}(p')$,  by
  (\ref{3.9}),
$$
\aligned
\big(\widetilde{\phi}\pm\phi\big)(y)\geq2\|\phi\|_{**}\mathcal{Z}(y)
\pm\phi(y)\geq\|\phi\|_{**}
\pm\phi(y)\geq|\phi(y)|
\pm\phi(y)
\geq0,
\endaligned
$$
and for
$y\in\Omega_t\setminus\big[\bigcup_{i=1}^mB_{R_1\gamma_i}(\xi'_i)\cup B_{R_1\rho_0 v_0/\varepsilon_0}(p')\cup B^c_{2d/\varepsilon_0}(p')\big]$,
by  (\ref{2.28}), (\ref{2.30}), (\ref{2.32}) and (\ref{3.7}),
$$
\aligned
\mathcal{L}\big(\widetilde{\phi}\pm\phi\big)(y)
\geq&
\|h\|_{*}\sum_{i=0}^m\mathcal{L}\big(\psi_i\big)(y)\pm
h(y)
=\|h\|_{*}\left\{-\sum_{i=0}^m\Delta\psi_i(y)
-W\sum_{i=0}^m\psi_i(y)\right\}\pm
h(y)\\[0.5mm]
\geq&\|h\|_{*}\left\{\left(\frac{\varepsilon_0}{\rho_0v_0}\right)^2
\frac{4}{\big|\frac{\varepsilon_0 y-p}{\rho_0v_0}\big|^{4+2\hat{\alpha}}}
+\sum_{i=1}^m\frac{1}{\gamma_i^2}
\frac{4}{\big|\frac{y-\xi'_i}{\gamma_i}\big|^{4+2\hat{\alpha}}}+4(m+1)\varepsilon_0^2
\right.\\[1mm]
&-C\left.\left[
\left(\frac{\varepsilon_0}{\rho_0v_0}\right)^2\frac{8(1+\alpha)^2\big|\frac{\varepsilon_0 y-p}{\rho_0v_0}\big|^{2\alpha}}{\,\big(1+\big|\frac{\varepsilon_0 y-p}{\rho_0v_0}\big|^{2(1+\alpha)}\big)^2\,}
+\sum_{i=1}^m\frac{1}{\gamma_{i}^2}\frac{8}{\big(1+\big|\frac{y-\xi'_i}{\gamma_{i}}\big|^2\big)^2}
+\varepsilon_0^2t^{4\beta(2m+2+\alpha)}e^{-t\phi_1(\varepsilon_0 y)}\right]
\right\}\pm
h(y)\\[1.5mm]
\geq&
|h(y)|\pm h(y)\geq0.
\endaligned
$$
By the  maximum
principle in the previous step  we obtain  that $-\tilde{\phi}\leq\phi\leq\tilde{\phi}$
in
$\widetilde{\Omega}_t=
\Omega_t\setminus\big[\bigcup_{i=1}^mB_{R_1\gamma_i}(\xi'_i)
\cup B_{R_1\rho_0 v_0/\varepsilon_0}(p')\cup B^c_{2d/\varepsilon_0}(p')\big]$,
which combined with (\ref{3.13}) gives estimate (\ref{3.12}).

\vspace{1mm}
\vspace{1mm}
\vspace{1mm}

{\bf Step 3:} Take $R_0=2R_1$, $R_1$ being the constant in the previous
two steps.
We prove uniform a priori estimates for solutions $\phi$ of
equation (\ref{3.11}), when $h\in L^{\infty}(\Omega_t)$ and $\phi$
satisfies more orthogonality conditions than those of (\ref{3.1})
in the following way
\begin{equation}\label{3.14}
\aligned
\int_{\Omega_t}\chi_pZ_p\phi=0
\,\qquad\,
\textrm{and}
\,\qquad\,\int_{\Omega_t}\chi_iZ_{ij}\phi=0,\,\,\,\,\,
\,\,\,\,\,i=1,\ldots,m,\,\,\,j=0,1,2.
\endaligned
\end{equation}
Namely, we prove that there exists a constant $C>0$ independent of $t$
such that for any $h\in L^{\infty}(\Omega_t)$ and any points $\xi=(\xi_1,\ldots,\xi_m)\in\mathcal{O}_t$,
\begin{equation}\label{3.15}
\aligned
\|\phi\|_{L^{\infty}(\Omega_t)}\leq C
\|h\|_{*},
\endaligned
\end{equation}
for $t$ large enough.
By contradiction,  assume that there are  sequences of
parameters $t_n\rightarrow+\infty$,
points $\xi^n=(\xi_1^n,\ldots,\xi_m^n)\in\mathcal{O}_{t_n}$,
functions $h_n$, $W_n$  and associated solutions $\phi_n$ of
equation (\ref{3.11}) with orthogonality conditions (\ref{3.14})
such that
\begin{equation}\label{3.16}
\aligned
\|\phi_n\|_{L^{\infty}(\Omega_{t_n})}=1
\,\,\,\,\quad\,\,\,\,
\textrm{and}
\,\,\,\,\quad\,\,\,\,
\|h_n\|_{*}\rightarrow0,
\,\,\quad\,\,\textrm{as}\,\,\,\,n\rightarrow+\infty.
\endaligned
\end{equation}
Let us  set
$$
\aligned
\widehat{\phi}^n_{p^c}(x)=\phi_n\big(x/\varepsilon_0^n\big),
\,\,\qquad\,\,
\,\,\quad\,\,
\widehat{h}^n_{p^c}(x)=h_n\big(x/\varepsilon_0^n\big)
\,\,\qquad\,\,\textrm{for all}\,\,\,x\in\Omega\setminus B_{2d}(p),
\endaligned
$$
and
$$
\aligned
\widehat{\phi}^n_p(z)=\phi_n\big((\rho_0^nv_0^nz+p)/\varepsilon_0^n\big),
\,\,\quad\,\,
\qquad
\,\,\quad\,\,
\widehat{h}^n_p(z)=h_n\big((\rho_0^nv_0^nz+p)/\varepsilon_0^n\big),
\endaligned
$$
and for all $i=1,\ldots,m$,
$$
\aligned
\widehat{\phi}^n_i(z)=\phi_n\big(\gamma_{i}^nz+(\xi^n_i)'\big),
\,\,\quad\,\,
\,\,\,\,\qquad\,\,\,\,
\,\,\quad\,\,
\widehat{h}^n_i(z)=h_n\big(\gamma_{i}^nz+(\xi^n_i)'\big),
\endaligned
$$
where
$\mu^n=\big(\mu^n_0,\mu^n_1,\ldots,\mu_m^n\big)$,
$\varepsilon_0^n=\exp\left\{-\frac12t_n\right\}$,
$\varepsilon^n_i=\exp\left\{-\frac12t_n\phi_1(\xi_i^n)\right\}$,
$\rho^n_0=(\varepsilon_0^n)^{\frac{1}{1+\alpha}}=\exp\left\{
-\frac1{2(1+\alpha)}t_n\right\}$,
$v^n_0=(\mu_0^n)^{\frac{1}{1+\alpha}}$ and
$\gamma^n_i=\frac1{\varepsilon_0^n}\varepsilon_i^n\mu^n_i
=\mu^n_i\exp\left\{
-\frac12t_n\big[\phi_1(\xi^n_i)-1\big]\right\}$.
First, using the expansion of $W_n$ in (\ref{2.32}),
we have that $\widehat{\phi}^n_{p^c}(x)$ satisfies
$$
\aligned
\left\{\aligned
&-\Delta\widehat{\phi}^n_{p^c}(x)+
O\left(\frac{e^{-t_n\phi_1(x)}}
{\,|x-p|^{4+2\alpha}\,}\prod_{i=1}^m\frac{1}{\,|x-\xi_i^n|^4\,}
\right)
\widehat{\phi}^n_{p^c}(x)
=\left(\frac{1}{\varepsilon^n_0}\right)^2\widehat{h}^n_{p^c}(x)
\,\,\quad\textrm{in}\,\quad
\Omega\setminus B_{2d}(p),\\[1mm]
&\widehat{\phi}^n_{p^c}(x)=0
\qquad\qquad\qquad\qquad\qquad\qquad
\qquad\qquad\qquad\qquad\qquad\qquad
\qquad\quad
\textrm{on}\,\qquad
\partial\Omega.
\endaligned
\right.
\endaligned
$$
By  the definition of the $\|\cdot\|_{*}$-norm in (\ref{3.7})
we find that
$\big(\frac{1}{\varepsilon^n_0}\big)^2\big|\widehat{h}^n_{p^c}(x)\big|\leq C\|h_n\|_{*}\rightarrow 0$
uniformly in $\Omega\setminus B_{2d}(p)$. Obviously, elliptic regularity theory  implies that
$\widehat{\phi}^n_{p^c}$
converges uniformly in $\Omega\setminus B_{2d}(p)$ to a trivial  solution $\widehat{\phi}^{\infty}_{p^c}$,
namely $\widehat{\phi}^{\infty}_{p^c}\equiv0$ in $\Omega\setminus B_{2d}(p)$.

Next, using the expansion of $W_n$ in (\ref{2.30}),
we find that $\widehat{\phi}^n_p(z)$ satisfies
$$
\aligned
-\Delta\widehat{\phi}^n_p(z)-\frac{8(1+\alpha)^2\big|z\big|^{2\alpha}}{\,\big(1+\big|z\big|^{2(1+\alpha)}\big)^2\,}
\left[1+O\left(\rho^n_0v^n_0t_n^{\beta}|z|\right)
+
o\left(1\right)\right]
\widehat{\phi}^n_p(z)=\left(\frac{\rho^n_0v^n_0}{\varepsilon^n_0}\right)^2\widehat{h}^n_p(z)
\endaligned
$$
for any $z\in B_{R_0+2}(0)$.
Thanks to the definition of the $\|\cdot\|_{*}$-norm in (\ref{3.7}), we have that
 for any $q\in\big(1, -1/\hat{\alpha}\big)$,
$\big(\frac{\rho^n_0v^n_0}{\varepsilon^n_0}\big)^2\widehat{h}^n_p\rightarrow 0$ in
$L^{q}\big(B_{R_0+2}(0)\big)$.
Since
$\frac{8(1+\alpha)^2|z|^{2\alpha}}{(1+|z|^{2(1+\alpha)})^2}$
is bounded in $L^{q}\big(B_{R_0+2}(0)\big)$,
elliptic regularity theory readily implies that
$\widehat{\phi}^n_p$
converges uniformly over
compact subsets near the origin to a bounded solution $\widehat{\phi}^{\infty}_p$ of equation
$(\ref{3.3})$, which satisfies
\begin{equation}\label{3.17}
\aligned
\int_{\mathbb{R}^2}\chi \mathcal{Z}_p\widehat{\phi}_p^{\infty}=0.
\endaligned
\end{equation}
Then $\widehat{\phi}^{\infty}_p$
is proportional to $\mathcal{Z}_p$.
Since $\int_{\mathbb{R}^2}\chi \mathcal{Z}_p^2>0$, by
(\ref{3.17}) we deduce that $\widehat{\phi}^{\infty}_p\equiv0$ in $B_{R_1}(0)$.

Finally,
using the expansion of $W_n$ in (\ref{2.28})
and elliptic regularity, we can derive  that for each $i\in\{1,\ldots,m\}$,
$\widehat{\phi}^n_i$
converges uniformly over
compact subsets near the origin to a bounded solution $\widehat{\phi}^{\infty}_i$ of equation
$(\ref{3.4})$, which satisfies
\begin{equation}\label{3.18}
\aligned
\int_{\mathbb{R}^2}\chi \mathcal{Z}_j\widehat{\phi}_i^{\infty}=0
\quad\,\,\,\textrm{for}\,\,\,\,j=0,\,1,\,2.
\endaligned
\end{equation}
Then $\widehat{\phi}^{\infty}_i$ is
a linear combination of $\mathcal{Z}_j$, $j=0,1,2$.
Notice that $\int_{\mathbb{R}^2}\chi \mathcal{Z}_j\mathcal{Z}_l=0$ for $j\neq l$
and $\int_{\mathbb{R}^2}\chi \mathcal{Z}_j^2>0$.
Hence (\ref{3.18}) implies  $\widehat{\phi}^{\infty}_i\equiv0$ in $B_{R_1}(0)$.
As a consequence, by definition (\ref{3.10}) we find
$\lim_{n\rightarrow+\infty}\|\phi_n\|_{**}=0$. But
(\ref{3.12}) and (\ref{3.16}) tell us
$\liminf_{n\rightarrow+\infty}\|\phi_n\|_{**}>0$,
which is a contradiction.

\vspace{1mm}
\vspace{1mm}
\vspace{1mm}

{\bf Step 4:} We establish uniform an a priori estimate for solutions $\phi$ to
equation (\ref{3.11}), when $h\in L^{\infty}(\Omega_t)$ and $\phi$  only satisfies the
 orthogonality conditions in (\ref{3.1})
\begin{equation}\label{3.19}
\aligned
\int_{\Omega_t}\chi_iZ_{ij}\phi=0,\,\,\,\,\,
\,\,\,\,\,\,i=1,\ldots,m,\,\,j=1,2.
\endaligned
\end{equation}
More precisely, we prove that there exists a constant $C>0$ independent of $t$
such that for any $h\in L^{\infty}(\Omega_t)$ and any points $\xi=(\xi_1,\ldots,\xi_m)\in\mathcal{O}_t$,
\begin{equation}\label{3.20}
\aligned
\|\phi\|_{L^{\infty}(\Omega_t)}\leq Ct
\|h\|_{*},
\endaligned
\end{equation}
for $t$ large enough.

Let $R>R_0+1$ be a large but fixed number. Set
\begin{equation}\label{3.21}
\aligned
\widehat{Z}_{p}(y)=Z_{p}(y)-\frac{\varepsilon_0}{\rho_0v_0}
+a_{p}G(\varepsilon_0 y,p),
\,\qquad\qquad\,
\widehat{Z}_{i0}(y)=Z_{i0}(y)-\frac1{\gamma_i}
+a_{i0}G(\varepsilon_0 y,\xi_i),
\endaligned
\end{equation}
where
\begin{equation}\label{3.22}
\aligned
a_{p}=\frac{\varepsilon_0}{\rho_0v_0\big[H(p,p)-4\log(\rho_0v_0 R)\big]},
\,\qquad\quad\qquad\,
a_{i0}=\frac1{\gamma_i\big[H(\xi_i,\xi_i)-4\log(\varepsilon_0\gamma_i R)\big]}.
\endaligned
\end{equation}
Note that by estimates (\ref{2.14})-(\ref{2.15}), and definitions (\ref{2.2}), (\ref{2.7})  and (\ref{2.23}),
\begin{equation}\label{3.23}
\aligned
C_1|\log\varepsilon_0|\leq-\log(\rho_0v_0 R)
\leq C_2|\log\varepsilon_0|,
\qquad\qquad
C_1|\log\varepsilon_i|\leq-\log(\varepsilon_0\gamma_i R)
\leq C_2|\log\varepsilon_i|,
\endaligned
\end{equation}
and
\begin{equation}\label{3.24}
\aligned
\widehat{Z}_{p}(y)=O\left(
\frac{\varepsilon_0G(\varepsilon_0 y,p)}{\rho_0v_0|\log\varepsilon_0|}
\right),
\,\,\,\qquad\quad\qquad\,
\widehat{Z}_{i0}(y)=O\left(
\frac{G(\varepsilon_0 y,\xi_i)}{\gamma_i|\log\varepsilon_i|}
\right).
\endaligned
\end{equation}
Let   $\eta_1$ and $\eta_2$ be  radial smooth cut-off functions in $\mathbb{R}^2$ such that
$$
\aligned
&0\leq\eta_1\leq1;\,\,\,\,\,\,\,|\nabla\eta_1|\leq C\,\,\textrm{in}\,\mathbb{R}^2;
\,\,\,\,\,\,\,\eta_1\equiv1\,\,\textrm{in}\,B_R(0);\,\,\,\,\,\,\,
\,\,\,
\eta_1\equiv0\,\,\textrm{in}\,\mathbb{R}^2\setminus B_{R+1}(0);\\[1mm]
&0\leq\eta_2\leq1;\,\,\,\,\,\,\,|\nabla\eta_2|\leq C\,\,\textrm{in}\,\mathbb{R}^2;
\,\,\,\,\,\,\,\eta_2\equiv1\,\,\textrm{in}\,B_{3d}(0);\,\,\,
\,\,\,\,\,\,\,\eta_2\equiv0\,\,\textrm{in}\,\mathbb{R}^2\setminus B_{6d}(0),
\endaligned
$$
where $d>0$ can be chosen as a sufficiently small but fixed number independent of $t$ such that
$B_{9d}(p)\subset\Omega$.
Set
\begin{equation}\label{3.25}
\aligned
\eta_{p1}(y)=
\eta_1\left(\frac{\big|\varepsilon_0y-p\big|}{\rho_0 v_0}\right),
\,\,\qquad\qquad\,\,
\eta_{i1}(y)=
\eta_1\left(\frac{\big|y-\xi_i'\big|}{\gamma_i}\right),
\endaligned
\end{equation}
and
\begin{equation}\label{3.26}
\aligned
\eta_{p2}(y)=
\eta_2\left(\varepsilon_0\big|y-p'\big|\right),
\,\,\,\qquad\,\,\qquad\,\,\,
\eta_{i2}(y)=
\eta_2\left(\varepsilon_0\big|y-\xi'_i\big|\right).
\endaligned
\end{equation}
We  define the two test functions
\begin{equation}\label{3.27}
\aligned
\widetilde{Z}_{p}=\eta_{p1}Z_{p}+(1-\eta_{p1})\eta_{p2}\widehat{Z}_{p},
\,\,\,\qquad\,\,\qquad\,\,\,
\widetilde{Z}_{i0}=\eta_{i1}Z_{i0}+(1-\eta_{i1})\eta_{i2}\widehat{Z}_{i0}.
\endaligned
\end{equation}
Given $\phi$ satisfying (\ref{3.11}) and (\ref{3.19}), let
\begin{equation}\label{3.28}
\aligned
\widetilde{\phi}=\phi+d_p\widetilde{Z}_{p}+\sum\limits_{i=1}^{m}d_i\widetilde{Z}_{i0}+\sum_{i=1}^m\sum\limits_{j=1}^{2}e_{ij}\chi_iZ_{ij}.
\endaligned
\end{equation}
We will first prove the existence of $d_p$,
$d_i$ and $e_{ij}$  such
that $\widetilde{\phi}$ satisfies the orthogonality conditions in (\ref{3.14}).
Remark  that $\widetilde{Z}_{i0}$
coincides with $Z_{i0}$ in $B_{R\gamma_i}(\xi'_i)$ and hence
$\widetilde{Z}_{i0}$ is still orthogonal
to $\chi_iZ_{ij}$ for $j=1, 2$.
Testing  (\ref{3.28}) against $\chi_iZ_{ij}$ and using
the orthogonality conditions in (\ref{3.14}) and (\ref{3.19}) for $j=1,2$ and the fact that
$\chi_i\chi_k\equiv0$ if $i\neq k$, we can write
\begin{equation}\label{3.29}
\aligned
e_{ij}=\left(-d_p\int_{\Omega_t}\chi_iZ_{ij}\widetilde{Z}_{p}
-\sum_{k\neq i}^md_k\int_{\Omega_t}\chi_iZ_{ij}\widetilde{Z}_{k0}
\right)\left/\int_{\Omega_t}\chi^2_iZ^2_{ij},
\,\,\quad\,
i=1,\ldots,m,\,\,j=1,2.
\right.
\endaligned
\end{equation}
Notice that $\int_{\Omega_t}\chi^2_iZ^2_{ij}=c>0$ for all $i$, $j$, and
by  (\ref{3.24}) and (\ref{3.27}),
$$
\aligned
\int_{\Omega_t}\chi_iZ_{ij}\widetilde{Z}_{p}=
O\left(\frac{\varepsilon_0\gamma_i\log t}{\rho_0v_0|\log\varepsilon_0|}
\right),
\,\,\,\qquad\,\,\,
\int_{\Omega_t}\chi_iZ_{ij}\widetilde{Z}_{k0}=
O\left(\frac{\gamma_i\log t}{\gamma_k|\log\varepsilon_k|}
\right),\,
\quad\,\,\,k\neq i.
\endaligned
$$
Then
\begin{equation}\label{3.30}
\aligned
|e_{ij}|\leq C\left(|d_p|\frac{\varepsilon_0\gamma_i\log t}{\rho_0v_0|\log\varepsilon_0|}
+\sum_{k\neq i}^m|d_k|\frac{\gamma_i\log t}{\gamma_k|\log\varepsilon_k|}\right).
\endaligned
\end{equation}
We need just to consider $d_p$  and $d_i$. Testing  (\ref{3.28}) against $\chi_pZ_p$ and $\chi_kZ_{k0}$,
respectively, and using
the orthogonality conditions in (\ref{3.14})
for $p$ and $j=0$, we get a system of $(d_p,d_1,\ldots,d_m)$,
\begin{equation}\label{3.31}
\aligned
d_p\int_{\Omega_t}\chi_pZ_p\widetilde{Z}_{p}
+\sum_{i=1}^md_i\int_{\Omega_t}\chi_pZ_p\widetilde{Z}_{i0}=&
-\int_{\Omega_t}\chi_pZ_p\phi,
\\[1mm]
d_p\int_{\Omega_t}\chi_kZ_{k0}\widetilde{Z}_{p}
+\sum_{i=1}^md_i\int_{\Omega_t}\chi_kZ_{k0}\widetilde{Z}_{i0}=&
-\int_{\Omega_t}\chi_kZ_{k0}\phi,
\,\quad\quad\,k=1,\ldots,m.
\endaligned
\end{equation}
But
$$
\aligned
\int_{\Omega_t}\chi_pZ_p\widetilde{Z}_{p}=
\int_{\Omega_t}\chi_pZ_p^2=C_1>0,
\,\,\,\qquad\qquad\,\,\,
\int_{\Omega_t}\chi_pZ_p\widetilde{Z}_{i0}=
O\left(\frac{\rho_0v_0\log t}{\varepsilon_0\gamma_i|\log\varepsilon_i|}
\right),
\endaligned
$$
and
$$
\aligned
\int_{\Omega_t}\chi_kZ_{k0}\widetilde{Z}_{p}=
O\left(
\frac{\varepsilon_0\gamma_k\log t}{\rho_0v_0|\log\varepsilon_0|}
\right),
\qquad
\int_{\Omega_t}\chi_kZ_{k0}\widetilde{Z}_{k0}=C_2>0,
\qquad
\int_{\Omega_t}\chi_kZ_{k0}\widetilde{Z}_{i0}=
O\left(\frac{\gamma_k\log t}{\gamma_i|\log\varepsilon_i|}
\right),
\quad
i\neq k.
\endaligned
$$
Let us denote $\mathcal{M}$ the coefficient matrix of system (\ref{3.31}).
From the
above estimates it follows that
$P^{-1}\mathcal{M}P$ is diagonally dominant and then invertible, where
$P=\diag\left(\rho_0 v_0\big/\varepsilon_0,\,\gamma_1,\,\ldots,\,\gamma_m\right)$. Hence $\mathcal{M}$ is also invertible and
$(d_p,d_1,\ldots,d_m)$ is well defined.

Estimate (\ref{3.20}) is a  direct consequence of the following two claims.

\vspace{1mm}
\vspace{1mm}
\vspace{1mm}

\noindent{\bf Claim 1.}\,\,{\it
\begin{equation}\label{3.32}
\aligned
\big\|\mathcal{L}(\widetilde{Z}_{p})\big\|_{*}\leq
\frac{C\varepsilon_0\log t}{\rho_0v_0|\log\varepsilon_0|},
\endaligned
\end{equation}
and
\begin{equation}\label{3.33}
\aligned
\big\|\mathcal{L}(\chi_iZ_{ij})\big\|_{*}\leq\frac{C}{\gamma_i},
\,\quad\qquad\qquad\,\,\,\quad\,
\big\|\mathcal{L}(\widetilde{Z}_{i0})\big\|_{*}\leq
\frac{C\log t}{\gamma_i|\log\varepsilon_i|}.
\endaligned
\end{equation}
}

\noindent{\bf Claim 2.}\,\,{\it
\begin{equation}\label{3.34}
\aligned
|d_p|\leq C\frac{\rho_0v_0|\log\varepsilon_0|}{\varepsilon_0}\|h\|_{*},
\,\quad\quad\,\,\quad\,\,\quad
|d_i|\leq C\gamma_i|\log\varepsilon_i|\|h\|_{*},
\,\quad\quad\,\,\quad\,\,\quad
|e_{ij}|\leq C\gamma_i\log t\,\|h\|_{*}.
\endaligned
\end{equation}
}

In fact,  by  the definition of $\widetilde{\phi}$ in (\ref{3.28}) we get
\begin{equation}\label{3.35}
\aligned
\mathcal{L}(\widetilde{\phi})=h+d_p\mathcal{L}(\widetilde{Z}_{p})+\sum\limits_{i=1}^{m}d_i\mathcal{L}(\widetilde{Z}_{i0})
+\sum_{i=1}^m\sum_{j=1}^2e_{ij}\mathcal{L}(\chi_iZ_{ij})
\,\quad\,\textrm{in}\,\,\,\,\,\Omega_t.
\endaligned
\end{equation}
Since   (\ref{3.14}) hold, by estimate (\ref{3.15})  we conclude
\begin{equation}\label{3.36}
\aligned
\|\widetilde{\phi}\|_{L^{\infty}(\Omega_t)}\leq
C\left[\|h\|_{*}+|d_p|\big\|\mathcal{L}(\widetilde{Z}_{p})\big\|_{*}
+\sum\limits_{i=1}^{m}|d_i|\big\|\mathcal{L}(\widetilde{Z}_{i0})\big\|_{*}
+\sum_{i=1}^m\sum_{j=1}^2|e_{ij}|\big\|\mathcal{L}(\chi_i Z_{ij})\big\|_{*}
\right].
\endaligned
\end{equation}
Using the definition of $\widetilde{\phi}$ again and the fact that
\begin{equation}\label{3.37}
\aligned
\big\|\widetilde{Z}_{p}\big\|_{L^{\infty}(\Omega_t)}\leq\frac{C\varepsilon_0}{\rho_0v_0},
\quad\qquad\quad\quad
\big\|\widetilde{Z}_{i0}\big\|_{L^{\infty}(\Omega_t)}\leq\frac{C}{\gamma_i},
\quad\qquad\quad\quad
\big\|\chi_iZ_{ij}\big\|_{L^{\infty}(\Omega_t)}\leq\frac{C}{\gamma_i},
\endaligned
\end{equation}
estimate (\ref{3.20}) then follows from (\ref{2.7}),  (\ref{3.36}), Claims $1$ and  $2$.

\vspace{1mm}
\vspace{1mm}
\vspace{1mm}

\noindent{\bf Proof of Claim 1.}
Let us begin with  inequality  (\ref{3.32}).
Consider four regions
$$
\aligned
\Omega_1=\left\{\,\left|\frac{\varepsilon_0y-p}{\rho_0 v_0}\right|\leq R\right\},
\quad\quad\qquad\quad&\qquad\quad
\Omega_2=\left\{R<\left|\frac{\varepsilon_0y-p}{\rho_0 v_0}\right|\leq R+1\right\},\\[1mm]
\Omega_3=\left\{R+1<\left|\frac{\varepsilon_0y-p}{\rho_0 v_0}\right|\leq\frac{3d}{\rho_0v_0}\right\},
\quad\qquad&\qquad\quad
\Omega_4=\left\{\frac{3d}{\rho_0v_0}<\left|\frac{\varepsilon_0y-p}{\rho_0 v_0}\right|\leq
\frac{6d}{\rho_0v_0}\right\}.
\quad\quad\,\,\,\,\,\,
\endaligned
$$
Observe first that, by (\ref{3.2}), (\ref{3.3}) and (\ref{3.5}),
\begin{eqnarray}\label{3.38}
\mathcal{L}(Z_{p})
=-\Delta Z_{p}-W Z_{p}
=\left[\left(\frac{\varepsilon_0}{\rho_0v_0}\right)^2\frac{8(1+\alpha)^2\big|\frac{\varepsilon_0 y-p}{\rho_0v_0}\big|^{2\alpha}}{\,\big(1+\big|\frac{\varepsilon_0 y-p}{\rho_0v_0}\big|^{2(1+\alpha)}\big)^2\,}
-W\right]Z_{p}.
\end{eqnarray}
%\begin{eqnarray}\label{3.38}
%\mathcal{L}(Z_{p})
%=\left[-\Delta Z_{p}-\left(\frac{\varepsilon_0}{\rho_0v_0}\right)^2\frac{8(1+\alpha)^2\big|\frac{\varepsilon_0 y-p}{\rho_0v_0}\big|^{2\alpha}}{\,\big(1+\big|\frac{\varepsilon_0 y-p}{\rho_0v_0}\big|^{2(1+\alpha)}\big)^2\,}Z_{p}
%\right]+\left[\left(\frac{\varepsilon_0}{\rho_0v_0}\right)^2\frac{8(1+\alpha)^2\big|\frac{\varepsilon_0 y-p}{\rho_0v_0}\big|^{2\alpha}}{\,\big(1+\big|\frac{\varepsilon_0 y-p}{\rho_0v_0}\big|^{2(1+\alpha)}\big)^2\,}
%-W\right]Z_{p}
%&&\nonumber\\[1.5mm]
%=\left(\frac{\varepsilon_0}{\rho_0v_0}\right)^3\frac{8(1+\alpha)^2\big|\frac{\varepsilon_0 y-p}{\rho_0v_0}\big|^{2\alpha}}{\,\big(1+\big|\frac{\varepsilon_0 y-p}{\rho_0v_0}\big|^{2(1+\alpha)}\big)^2\,}
%\exp
%\left\{\,O\left(\varepsilon_0t^{\beta}|y-p'|\right)
%+
%O\left(\varepsilon_0^2\mu_0^2\right)
%+\sum_{j=1}^mO\left(\varepsilon_j^2\mu_j^2t^{2\beta}\right)\right\}.
%\qquad\,&&
%\end{eqnarray}
In $\Omega_1$, by  (\ref{2.30}), (\ref{3.27}) and (\ref{3.38}),
\begin{equation}\label{3.39}
\aligned
\mathcal{L}(\widetilde{Z}_{p})=\mathcal{L}(Z_{p})=
\left(\frac{\varepsilon_0}{\rho_0v_0}\right)^3
\frac{8(1+\alpha)^2\big|\frac{\varepsilon_0 y-p}{\rho_0v_0}\big|^{2\alpha}}{\,\big(1+\big|\frac{\varepsilon_0 y-p}{\rho_0v_0}\big|^{2(1+\alpha)}\big)^2\,}\left[O\left(\varepsilon_0t^{\beta}|y-p'|\right)
+
O\left(\varepsilon_0^2\mu_0^2\right)
+\sum_{j=1}^mO\left(\varepsilon_j^2\mu_j^2t^{2\beta}\right)\right].
\endaligned
\end{equation}
In $\Omega_2$, by (\ref{1.2}), (\ref{3.21}) and (\ref{3.27}),
\begin{eqnarray}\label{3.40}
\mathcal{L}(\widetilde{Z}_{p})=\mathcal{L}(Z_{p})-(1-\eta_{p1})\mathcal{L}(Z_{p}-\widehat{Z}_{p})
-2\nabla\eta_{p1}\nabla(Z_{p}-\widehat{Z}_{p})-(Z_{p}-\widehat{Z}_{p})\Delta \eta_{p1}
\,\,\,\,\,\,\,&&\nonumber\\[1mm]
=\mathcal{L}(Z_{p})+(1-\eta_{p1})W(Z_{p}-\widehat{Z}_{p})
-2\nabla\eta_{p1}\nabla(Z_{p}-\widehat{Z}_{p})-(Z_{p}-\widehat{Z}_{p})\Delta\eta_{p1}.
\,\,\,&&
\end{eqnarray}
Notice that,  by  (\ref{3.21})-(\ref{3.22}),
\begin{equation}\label{3.41}
\aligned
Z_{p}-\widehat{Z}_{p}=\frac{\varepsilon_0}{\rho_0v_0}
-a_{p}G(\varepsilon_0 y,p)
=\frac{\varepsilon_0}{\rho_0v_0\big[H(p,p)-4\log(\rho_0v_0 R)\big]}\left[
4\log\frac{|\varepsilon_0y-p|}{R\rho_0v_0}+O\left(
\varepsilon_0|y-p'|\right)
\right],
\endaligned
\end{equation}
and then in $\Omega_2$, by (\ref{3.23}),
\begin{equation}\label{3.42}
\aligned
|Z_{p}-\widehat{Z}_{p}|=O\left(
\frac{\varepsilon_0}{R\rho_0v_0|\log\varepsilon_0|}
\right),
\,\,\quad\quad\,\quad\,\,\quad\quad\,
|\nabla\big(Z_{p}-\widehat{Z}_{p}\big)|=O\left(
\frac{\varepsilon_0^2}{R\rho_0^2v_0^2|\log\varepsilon_0|}
\right).
\endaligned
\end{equation}
Moreover, $|\nabla\eta_{p1}|=O\big(\varepsilon_0/(\rho_0v_0)\big)$ and
$|\Delta \eta_{p1}|=O\big(\varepsilon_0^2/(\rho_0^2v_0^2)\big)$.
By (\ref{2.30}), (\ref{3.38}),  (\ref{3.40})   and (\ref{3.42}) we have
that in $\Omega_2$,
\begin{equation}\label{3.43}
\aligned
\mathcal{L}(\widetilde{Z}_{p})
=O\left(\frac{\varepsilon_0^3}{R\rho_0^3v_0^3|\log\varepsilon_0|}
\right).
\endaligned
\end{equation}
In $\Omega_3$, by (\ref{1.2}),  (\ref{3.21}), (\ref{3.27}) and (\ref{3.38}),
\begin{eqnarray*}
\mathcal{L}(\widetilde{Z}_{p})=\mathcal{L}(\widehat{Z}_{p})
=\mathcal{L}(Z_{p})-\mathcal{L}(Z_{p}-\widehat{Z}_{p})
\,\quad\quad\quad\quad\quad\quad\quad\quad\quad\quad\quad
\qquad\qquad\qquad\quad\quad\quad\quad\quad\,\,\,\,
&&\nonumber\\[1.6mm]
=\left[\left(\frac{\varepsilon_0}{\rho_0v_0}\right)^2\frac{8(1+\alpha)^2\big|\frac{\varepsilon_0 y-p}{\rho_0v_0}\big|^{2\alpha}}{\,\big(1+\big|\frac{\varepsilon_0 y-p}{\rho_0v_0}\big|^{2(1+\alpha)}\big)^2\,}
-W\right]Z_{p}
+W\left[\frac{\varepsilon_0}{\rho_0v_0}
-a_{p}G(\varepsilon_0 y,p)\right]
\equiv A_1+A_2.
&&
\end{eqnarray*}
For the estimates of these two terms, we decompose $\Omega_3$  to some subregions:
$$
\aligned
&\Omega_{p}=\left\{\,R+1<\left|\frac{\varepsilon_0y-p}{\rho_0 v_0}\right|\leq
\frac{1}{\rho_0v_0 t^{2\beta}}\,\right\},\\
\Omega_{3,k}=\Omega_3\,\bigcap\,&\big\{\,|y-\xi'_k|\leq1/(\varepsilon_0 t^{2\beta})\,\big\}
\,\,\,\,\quad\,\,\textrm{and}\quad\,\,
\widetilde{\Omega}_3=\Omega_3\setminus
\left[\bigcup_{k=1}^m\Omega_{3,k}\cup \Omega_{p}\right].
\endaligned
$$
By (\ref{2.30}), (\ref{2.32}) and (\ref{3.5}),
$$
\aligned
A_1=\left\{
\aligned
&
\left(\frac{\varepsilon_0}{\rho_0v_0}\right)^3\frac{8(1+\alpha)^2\big|\frac{\varepsilon_0 y-p}{\rho_0v_0}\big|^{2\alpha}}{\,\big(1+\big|\frac{\varepsilon_0 y-p}{\rho_0v_0}\big|^{2(1+\alpha)}\big)^2\,}
\left[\,O\left(\varepsilon_0t^{\beta}|y-p'|\right)
+
O\left(\varepsilon_0^2\mu_0^2\right)
+\sum_{j=1}^mO\left(\varepsilon_j^2\mu_j^2t^{2\beta}\right)\right]
\,\,\ \,\textrm{in}\,\,\,\,\Omega_{p},\\[2mm]
&
O\left(\frac{\varepsilon_0^5}{\rho_0v_0}\mu_0^2t^{4\beta(2+\alpha)}
\right)
+O\left(\frac{\varepsilon_0^3}{\rho_0v_0}e^{-t\phi_1(\varepsilon_0 y)}t^{4\beta(2m+2+\alpha)}
\right)
\quad\qquad\qquad\qquad\quad
\qquad\qquad
\textrm{in}\,\,\,\,\,\widetilde{\Omega}_{3}.
\endaligned
\right.
\endaligned
$$
Moreover, by (\ref{3.23}) and (\ref{3.41}),
$$
\aligned
A_2=\left\{
\aligned
&\left(\frac{\varepsilon_0}{\rho_0v_0}\right)^3\frac{8(1+\alpha)^2\big|\frac{\varepsilon_0 y-p}{\rho_0v_0}\big|^{2\alpha}}{\,\big(1+\big|\frac{\varepsilon_0 y-p}{\rho_0v_0}\big|^{2(1+\alpha)}\big)^2\,}
O\left(\frac{\log |\varepsilon_0y-p|
-\log(R\rho_0v_0)+
\varepsilon_0|y-p'|}{|\log\varepsilon_0|}\right)
\,\,\,\,\,\,\textrm{in}\,\,\,\Omega_{p},\\[2mm]
&O\left(\frac{\varepsilon_0^3}{\rho_0v_0}t^{4\beta(2m+2+\alpha)}e^{-t\phi_1(\varepsilon_0 y)}
\right)
\qquad\qquad\qquad\qquad\qquad
\qquad\qquad\qquad\qquad\qquad
\textrm{in}\,\,\,\,\widetilde{\Omega}_{3}.
\endaligned
\right.
\endaligned
$$
Then in $\Omega_{p}\cup \widetilde{\Omega}_{3}$,
\begin{equation}\label{3.44}
\aligned
\mathcal{L}(\widetilde{Z}_{p})=
\mathcal{L}(\widehat{Z}_{p})=\left(\frac{\varepsilon_0}{\rho_0v_0}\right)^3\frac{8(1+\alpha)^2\big|\frac{\varepsilon_0 y-p}{\rho_0v_0}\big|^{2\alpha}}{\,\big(1+\big|\frac{\varepsilon_0 y-p}{\rho_0v_0}\big|^{2(1+\alpha)}\big)^2\,}
O\left(\frac{\log |\varepsilon_0y-p|
-\log(R\rho_0v_0)}{|\log\varepsilon_0|}\right).
\endaligned
\end{equation}
In $\Omega_{3,k}$ with all $k$, by (\ref{2.28}),  (\ref{3.24}) and (\ref{3.38}),
\begin{eqnarray}\label{3.45}
\mathcal{L}(\widetilde{Z}_{p})=\mathcal{L}(\widehat{Z}_{p})=
-\Delta Z_{p}-W\widehat{Z}_{p}
\qquad\qquad\qquad\qquad\qquad\qquad
\qquad\qquad\qquad\qquad\qquad\,\,
&&\nonumber
\\[1.5mm]
=\left(\frac{\varepsilon_0}{\rho_0v_0}\right)^2\frac{8(1+\alpha)^2\big|\frac{\varepsilon_0 y-p}{\rho_0v_0}\big|^{2\alpha}}{\,\big(1+\big|\frac{\varepsilon_0 y-p}{\rho_0v_0}\big|^{2(1+\alpha)}\big)^2\,}Z_{p}
+O\left(\frac{1}{\gamma_{k}^2}\frac{8}{\big(1+\big|\frac{y-\xi'_k}{\gamma_{k}}\big|^2\big)^2}
\cdot\frac{\varepsilon_0G(\varepsilon_0 y,p)}{\rho_0v_0|\log\varepsilon_0|}
\right)
&&\nonumber
\\[1.5mm]
=O\left(\frac{1}{\gamma_{k}^2}\frac{8}{\big(1+\big|\frac{y-\xi'_k}{\gamma_{k}}\big|^2\big)^2}
\cdot\frac{\varepsilon_0\log t}{\rho_0v_0|\log\varepsilon_0|}
\right).
\qquad\qquad\qquad\qquad\qquad\qquad
\qquad\quad\,\,\,\,&&
\end{eqnarray}
Finally in $\Omega_4$, by (\ref{3.21}) and (\ref{3.27}),
\begin{equation}\label{3.46}
\aligned
\mathcal{L}(\widetilde{Z}_{p})&=-\eta_{p2}\Delta Z_{p}-\eta_{p2}W\widehat{Z}_{p}
-2\nabla\eta_{p2}\nabla \widehat{Z}_{p}-\widehat{Z}_{p}\Delta\eta_{p2}.
\endaligned
\end{equation}
Note that from the previous choice of the number
$d$ we get that  for any $y\in\Omega_4$ and any $k=1,\ldots,m$,
$$
\aligned
|y-\xi'_k|\geq|y-p'|-|p'-\xi_k'|\geq\frac{3d}{\varepsilon_0}
-\frac{d}{\varepsilon_0}=\frac{2d}{\varepsilon_0}>\frac{1}{\varepsilon_0 t^{2\beta}}.
\endaligned
$$
This combined with (\ref{2.32}) gives
\begin{equation}\label{3.47}
\aligned
W=O\left(\varepsilon_0^2e^{-t\phi_1(\varepsilon_0 y)}\right)
\,\,\quad\,\,\textrm{in}\,\,\,\,\Omega_4.
\endaligned
\end{equation}
In addition,
$|\nabla\eta_{p2}|=O\left(\varepsilon_0\right)$,
$|\Delta \eta_{p2}|=O\left(\varepsilon_0^2\right)$
and
\begin{equation}\label{3.48}
\aligned
|\widehat{Z}_{p}|=O\left(
\frac{\varepsilon_0}{\rho_0v_0|\log\varepsilon_0|}
\right),\,\,\qquad\qquad\,
|\nabla\widehat{Z}_{p}|=O\left(
\frac{\varepsilon_0^2}{\rho_0v_0|\log\varepsilon_0|}
\right)
\,\qquad\,\textrm{in}\,\,\,\,\,\,\Omega_4.
\endaligned
\end{equation}
Hence by (\ref{3.38}), (\ref{3.46}), (\ref{3.47}) and (\ref{3.48}), we find that in $\Omega_4$,
\begin{equation}\label{3.49}
\aligned
\mathcal{L}(\widetilde{Z}_{p})
=O\left(
\frac{\varepsilon_0^3}{\rho_0v_0|\log\varepsilon_0|}
\right).
\endaligned
\end{equation}
Combining  (\ref{3.7}), (\ref{3.39}), (\ref{3.43}), (\ref{3.44}), (\ref{3.45}) and (\ref{3.49}),
we readily conclude
$$
\aligned
\big\|\mathcal{L}(\widetilde{Z}_{p})\big\|_{*}=O\left(
\frac{\varepsilon_0\log t}{\rho_0v_0|\log\varepsilon_0|}\right).
\endaligned
$$

The  inequalities in (\ref{3.33}) are easy to establish as they are very
similar to the consideration of inequality (\ref{3.32}), so we leave the detailed proof for
readers.

\vspace{1mm}
\vspace{1mm}
\vspace{1mm}

\noindent{\bf Proof of Claim 2.}
Let us prove the first two inequalities in (\ref{3.34}).
Testing  (\ref{3.35}) against $\widetilde{Z}_{p}$ and using estimates (\ref{3.36}) and (\ref{3.37}), we find
$$
\aligned
d_p
\int_{\Omega_t}&\widetilde{Z}_{p}\mathcal{L}(\widetilde{Z}_{p})
+
\sum_{k=1}^md_k
\int_{\Omega_t}\widetilde{Z}_{p}\mathcal{L}(\widetilde{Z}_{k0})
\\
=&-\int_{\Omega_t} h\widetilde{Z}_{p}
+\int_{\Omega_t}\widetilde{\phi}\mathcal{L}(\widetilde{Z}_{p})-\sum_{k=1}^m\sum_{l=1}^2e_{kl}
\int_{\Omega_t}\chi_kZ_{kl}\mathcal{L}(\widetilde{Z}_{p})
\\[1mm]
\leq&\frac{C\varepsilon_0}{\,\rho_0v_0\,}\|h\|_{*}
+C\big\|\mathcal{L}(\widetilde{Z}_{p})\big\|_{*}
\left(\|\widetilde{\phi}\|_{L^{\infty}(\Omega_t)}
+\sum_{k=1}^m\sum_{l=1}^2\frac{1}{\gamma_k}|e_{kl}|\right)
\\[1mm]
\leq&\frac{C\varepsilon_0}{\,\rho_0v_0\,}\|h\|_{*}
+
C\big\|\mathcal{L}(\widetilde{Z}_{p})\big\|_{*}
\left[\|h\|_{*}
+|d_p|\big\|\mathcal{L}(\widetilde{Z}_{p})\big\|_{*}
+
\sum\limits_{k=1}^{m}|d_k|\big\|\mathcal{L}(\widetilde{Z}_{k0})\big\|_{*}
+\sum_{k=1}^m\sum_{l=1}^2|e_{kl}|\left(\frac{1}{\gamma_k}
+\big\|\mathcal{L}(\chi_kZ_{kl})\big\|_{*}
\right)
\right],
\endaligned
$$
where we have  applied  the following  two inequalities:
$$
\aligned
\left(\frac{\varepsilon_0}{\rho_0v_0}\right)^2\int_{\Omega_t}
\frac{\big|\frac{\varepsilon_0 y-p}{\rho_0v_0}\big|^{2\alpha}}
{\,\big(1+\big|\frac{\varepsilon_0 y-p}{\rho_0v_0}\big|\big)^{4+2\hat{\alpha}+2\alpha}\,}
dy\leq C
\,\,\,\quad\,\,\textrm{and}\,\,\,\quad\,\,
\int_{\Omega_t}\frac{1}{\gamma_i^2}\frac{1}{\big(1+\big|\frac{y-\xi'_i}{\gamma_i}\big|\big)^{4+2\hat{\alpha}}}
dy\leq C,
\,\,\,\,\,\,\,i=1,\ldots,m.
\endaligned
$$
But estimate (\ref{3.30}) and Claim $1$ imply
\begin{equation}\label{3.50}
\aligned
|d_p|\left|
\int_{\Omega_t}\widetilde{Z}_{p}\mathcal{L}(\widetilde{Z}_{p})
\right|
\leq
\frac{C\varepsilon_0}{\rho_0v_0}\|h\|_{*}
+\frac{C\varepsilon_0\log^2 t}{\rho_0v_0|\log\varepsilon_0|}\left[
\frac{\varepsilon_0|d_p|}{\rho_0v_0|\log\varepsilon_0|}
+
\sum_{k=1}^m\frac{|d_k|}{\gamma_k|\log\varepsilon_k|}
\right]
+\sum_{k=1}^m\left|d_k
\int_{\Omega_t}\widetilde{Z}_{k0}\mathcal{L}(\widetilde{Z}_{p})
\right|.
\endaligned
\end{equation}
Similarly, testing  (\ref{3.35}) against $\widetilde{Z}_{i0}$ and using
(\ref{3.30}), (\ref{3.36}), (\ref{3.37}) and Claim $1$, we can derive that
\begin{eqnarray}\label{3.51}
|d_i|\left|
\int_{\Omega_t}\widetilde{Z}_{i0}\mathcal{L}(\widetilde{Z}_{i0})
\right|
\leq
\frac{C\|h\|_{*}}{\gamma_i}
+\frac{C\log^2t}{\gamma_i|\log\varepsilon_i|}
\left[
\frac{\varepsilon_0|d_p|}{\rho_0v_0|\log\varepsilon_0|}
+\sum_{k=1}^m\frac{|d_k|}{\gamma_k|\log\varepsilon_k|}
\right]
+
\left|d_p
\int_{\Omega_t}\widetilde{Z}_{i0}\mathcal{L}(\widetilde{Z}_{p})
\right|
&&\nonumber\\[1mm]
+
\sum_{k\neq i}^m\left|d_k
\int_{\Omega_t}\widetilde{Z}_{k0}\mathcal{L}(\widetilde{Z}_{i0})
\right|.
\qquad\qquad\qquad\qquad\qquad\,
\qquad\qquad\qquad\qquad\quad\,\,\,\,
&&
\end{eqnarray}
We   decompose
\begin{equation}\label{3.52}
\aligned
\int_{\Omega_t}\widetilde{Z}_{p}\mathcal{L}(\widetilde{Z}_{p})
=\sum_{l=1}^4\int_{\Omega_l}\widetilde{Z}_{p}\mathcal{L}(\widetilde{Z}_{p})
=\sum_{l=1}^4I_{l}.
\endaligned
\end{equation}
By (\ref{3.5}) and (\ref{3.39}), we get
\begin{eqnarray}\label{3.53}
I_1=\int_{\Omega_l}Z_{p}\mathcal{L}(Z_{p})
=\int_{\Omega_1}
\left(\frac{\varepsilon_0}{\rho_0v_0}\right)^4
\frac{8(1+\alpha)^2\big|\frac{\varepsilon_0 y-p}{\rho_0v_0}\big|^{2\alpha}}{\,\big(1+\big|\frac{\varepsilon_0 y-p}{\rho_0v_0}\big|^{2(1+\alpha)}\big)^2\,}\left[O\left(\varepsilon_0t^{\beta}|y-p'|\right)
+
O\left(\varepsilon_0^2\mu_0^2\right)
+\sum_{j=1}^mO\left(\varepsilon_j^2\mu_j^2t^{2\beta}\right)\right]
&&\nonumber\\[1mm]
=\left(\frac{\varepsilon_0}{\rho_0v_0}\right)^2
\left[O\left(\rho_0v_0t^{\beta}\right)
+
O\left(\varepsilon_0^2\mu_0^2\right)
+\sum_{j=1}^mO\left(\varepsilon_j^2\mu_j^2t^{2\beta}\right)\right].
\qquad\qquad\qquad\qquad\quad
\qquad\quad\qquad\qquad\qquad
\quad\,\,\,\,&&
\end{eqnarray}
By (\ref{3.24}), (\ref{3.44}) and (\ref{3.45}), we have
\begin{eqnarray}\label{3.54}
I_3=\int_{\Omega_{p}\cup \widetilde{\Omega}_{3}}\widehat{Z}_{p}\mathcal{L}(\widehat{Z}_{p})
+\sum_{k=1}^m\int_{\Omega_{3,k}}\widehat{Z}_{p}\mathcal{L}(\widehat{Z}_{p})
\qquad\qquad\qquad\qquad\quad
\qquad\qquad\qquad\qquad
\qquad\qquad\qquad\qquad
&&\nonumber\\[1mm]
=\left(\frac{\varepsilon_0}{\rho_0v_0}\right)^2
O
\left(\int_{R+1}^{3d/(\rho_0v_0)}\frac{r^{1+2\alpha}\log(r/R)}{(1+r^{2(1+\alpha)})^2}
\frac{\log(\rho_0v_0 r)}{|\log\varepsilon_0|^2}dr
+\sum_{k=1}^m\int_{0}^{1/(\varepsilon_0\gamma_kt^{2\beta})}
\frac{r}{(1+r^2)^2}
\frac{\log^2 t}{|\log\varepsilon_0|^2} dr
\right)
&&\nonumber\\[1.5mm]
=\left(\frac{\varepsilon_0}{\rho_0v_0}\right)^2
\left[O\left(\frac{1}{R^{2(1+\alpha)}|\log\varepsilon_0|}\right)
+O\left(\frac{\log^2t}{|\log\varepsilon_0|^2}\right)\right].
\qquad\qquad\qquad\qquad\,\,\,\,
\qquad\qquad\qquad\qquad\quad\quad
&&
\end{eqnarray}
By (\ref{3.48}) and (\ref{3.49}), we derive that
\begin{equation}\label{3.55}
\aligned
I_4=\int_{\Omega_4}\eta_{p2}\widehat{Z}_{p}\mathcal{L}(\widetilde{Z}_{p})
=\int_{\left\{\frac{3d}{\rho_0v_0}<\left|\frac{\varepsilon_0y-p}{\rho_0 v_0}\right|\leq
\frac{6d}{\rho_0v_0}\right\}}
O\left(
\frac{\varepsilon_0^4}{\rho_0^2v_0^2|\log\varepsilon_0|^2}
\right)
dy
=O\left(
\frac{\varepsilon_0^2}{|\log\varepsilon_0|^2}
\right).
\endaligned
\end{equation}
Regarding the expression $I_2$, by (\ref{3.40}) we  get
$$
\aligned
I_2=-\int_{\Omega_2}\widetilde{Z}_{p}(Z_{p}-\widehat{Z}_{p})\Delta\eta_{p1}
-2\int_{\Omega_2}\widetilde{Z}_{p}\nabla\eta_{p1}\nabla(Z_{p}-\widehat{Z}_{p})+\int_{\Omega_2}\widetilde{Z}_{p}\big[
\mathcal{L}(Z_{p})+(1-\eta_{p1})W(Z_{p}-\widehat{Z}_{p})
\big].
\endaligned
$$
Integrating by parts the first term and using estimates (\ref{2.30}), (\ref{3.38}) and (\ref{3.42}) for the last term, we obtain
\begin{eqnarray}\label{3.56}
I_2=
-\int_{\Omega_2}Z_{p}\nabla\eta_{p1}\nabla(Z_{p}-\widehat{Z}_{p})
+\int_{\Omega_2}(Z_{p}-\widehat{Z}_{p})\nabla\eta_{p1}\nabla(Z_{p}-\widehat{Z}_{p})
+\int_{\Omega_2}(Z_{p}-\widehat{Z}_{p})^2|\nabla\eta_{p1}|^2
&&\nonumber\\[1.5mm]
+\int_{\Omega_2}(Z_{p}-\widehat{Z}_{p})\nabla\eta_{p1}\nabla\widehat{Z}_{p}
+\left(\frac{\varepsilon_0}{\rho_0v_0}\right)^2
O\left(\frac{1}{R^{3+2\alpha}|\log\varepsilon_0|}\right)
\qquad\qquad\quad\qquad
\qquad\quad\quad\,\,
&&\nonumber\\[1.5mm]
=I_{21}+I_{22}+I_{23}+I_{24}+\left(\frac{\varepsilon_0}{\rho_0v_0}\right)^2
O\left(\frac{1}{R^{3+2\alpha}|\log\varepsilon_0|}\right).
\qquad\qquad\qquad
\qquad\qquad\qquad\quad
&&
\end{eqnarray}
By (\ref{2.14}),  (\ref{2.23}),  (\ref{3.2}), (\ref{3.5}), (\ref{3.22}), (\ref{3.25}) and (\ref{3.41}),  we conclude
\begin{eqnarray}\label{3.57}
&&I_{21}=-a_{p}\left(\frac{\varepsilon_0}{\rho_0v_0}\right)^2
\int_{\left\{R<\left|\frac{\varepsilon_0y-p}{\rho_0 v_0}\right|\leq R+1\right\}}
\frac{1}{|y-p'|}
\mathcal{Z}_p\left(\frac{\varepsilon_0y-p}{\rho_0 v_0}\right)
\eta_1'\left(\frac{\big|\varepsilon_0y-p\big|}{\rho_0 v_0}\right)
\big(4+o(1)\big)dy\nonumber\\[1.5mm]
&&\,\quad\,\,=-8\pi a_{p}\frac{\varepsilon_0}{\rho_0v_0}\int_{R}^{R+1}
\eta_1'(r)\left[
1+O\left(\frac1{r^{2(1+\alpha)}}\right)
\right]dr\nonumber\\[2mm]
&&\,\quad\,\,=\frac{2\pi(1+\alpha)}{|\log\varepsilon_0|}\left(\frac{\varepsilon_0}{\rho_0v_0}\right)^2
\left[1+O\left(\frac1{R^{2(1+\alpha)}}\right)\right].
\end{eqnarray}
By  (\ref{3.2}), (\ref{3.5}), (\ref{3.25})  and (\ref{3.42}) we find
$|\nabla\eta_{p1}|=O\big(\frac{\varepsilon_0}{\rho_0v_0}\big)$
and $|\nabla\widehat{Z}_{p}|=O\big(\frac{\varepsilon_0^2}{R^{3+2\alpha}\rho_0^2v_0^2}\big)$ in $\Omega_2$. Furthermore,
\begin{equation}\label{3.58}
\aligned
I_{22}=O\left(\frac{\varepsilon_0^2}{R\rho_0^2v_0^2|\log\varepsilon_0|^2}\right),
\quad\quad\quad
I_{23}=O\left(\frac{\varepsilon_0^2}{R\rho_0^2v_0^2|\log\varepsilon_0|^2}\right),
\quad\quad\quad
I_{24}=O\left(\frac{\varepsilon_0^2}{R^{3+2\alpha}\rho_0^2v_0^2|\log\varepsilon_0|}\right).
\endaligned
\end{equation}
Substituting estimates (\ref{3.53})-(\ref{3.58}) into (\ref{3.52}), we conclude
that for $R$ and $t$ large enough,
\begin{equation}\label{3.59}
\aligned
\int_{\Omega_t}\widetilde{Z}_{p}\mathcal{L}(\widetilde{Z}_{p})
=\frac{2\pi(1+\alpha)}{|\log\varepsilon_0|}\left(\frac{\varepsilon_0}{\rho_0v_0}\right)^2
\left[\,1+O\left(\frac1{R^{2(1+\alpha)}}\right)\right].
\endaligned
\end{equation}

According to (\ref{3.50}), we need just to calculate  $\int_{\Omega_t}\widetilde{Z}_{k0}\mathcal{L}(\widetilde{Z}_{p})$
for all $k$.
By the above estimates of
$\mathcal{L}(\widetilde{Z}_{p})$  and $\widetilde{Z}_{k0}$, we can easily prove that
$$
\aligned
\int_{\Omega_1}\widetilde{Z}_{k0}\mathcal{L}(\widetilde{Z}_{p})
=O\left(\frac{\varepsilon_0
\big(
\rho_0v_0t^{\beta}
+\sum_{j=0}^m\varepsilon_j^2\mu_j^2t^{2\beta}
\big
)
\log t
}{\rho_0v_0\gamma_k|\log\varepsilon_k|}
\right),
\,\qquad\qquad\,
\int_{\Omega_2}\widetilde{Z}_{k0}\mathcal{L}(\widetilde{Z}_{p})
=O\left(\frac{\varepsilon_0\log t}
{\rho_0v_0\gamma_k|\log\varepsilon_0||\log\varepsilon_k|}\right),
\endaligned
$$
$$
\aligned
\int_{\Omega_4}\widetilde{Z}_{k0}\mathcal{L}(\widetilde{Z}_{p})=
O\left(\frac{\varepsilon_0 }{\rho_0v_0\gamma_k|\log\varepsilon_0||\log\varepsilon_k|}\right),
\,\,\quad\qquad\qquad\,\,
\int_{\Omega_{p}\cup \widetilde{\Omega}_{3}}\widetilde{Z}_{k0}\mathcal{L}(\widetilde{Z}_{p})
=O\left(\frac{\varepsilon_0\log t}{\rho_0v_0\gamma_k|\log\varepsilon_0||\log\varepsilon_k|}\right),
\endaligned
$$
and
$$
\aligned
\int_{\Omega_{3,l}}\widetilde{Z}_{k0}\mathcal{L}(\widetilde{Z}_{p})
=O\left(\frac{\varepsilon_0\log^2 t}{\rho_0v_0\gamma_k|\log\varepsilon_0||\log\varepsilon_k|}\right)
\quad\quad\textrm{for all}\,\,\,l\neq k.
\endaligned
$$
It remains to  calculate  the integral over $\Omega_{3,k}$. From (\ref{3.27}) and  an integration by parts  we have
$$
\aligned
\int_{\Omega_{3,k}}\widetilde{Z}_{k0}\mathcal{L}(\widetilde{Z}_{p})
=\int_{\Omega_{3,k}}\widehat{Z}_{p}\mathcal{L}(\widetilde{Z}_{k0})
-\int_{\partial\Omega_{3,k}}\widehat{Z}_{k0}\frac{\partial\widehat{Z}_{p}}{\partial\nu}
+\int_{\partial\Omega_{3,k}}\widehat{Z}_{p}\frac{\partial\widehat{Z}_{k0}}{\partial\nu}.
\endaligned
$$
Observe that
$$
\aligned
\int_{\Omega_{3,k}}\widehat{Z}_{p}\mathcal{L}(\widetilde{Z}_{k0})
=
\left(\int_{\big\{|y-\xi'_k|\leq\gamma_kR\big\}}
+
\int_{\big\{\gamma_kR<|y-\xi'_k|\leq\gamma_k(R+1)\big\}}
+
\int_{\big\{\gamma_k(R+1)<|y-\xi'_k|\leq
1/(\varepsilon_0t^{2\beta})\big\}}
\right)
\widehat{Z}_{p}\mathcal{L}(\widetilde{Z}_{k0}).
\endaligned
$$
By (\ref{2.28}), (\ref{3.2}), (\ref{3.5}),  (\ref{3.21}),  (\ref{3.22}), (\ref{3.23})  and (\ref{3.27})  we
can compute that for $|y-\xi'_k|\leq\gamma_kR$,
$$
\aligned
\mathcal{L}(\widetilde{Z}_{k0})=\mathcal{L}(Z_{k0})=
O\left(\varepsilon_0t^{\beta}/\gamma_k^2\right),
\endaligned
$$
for $\gamma_kR<|y-\xi'_k|\leq\gamma_k(R+1)$,
$$
\aligned
\mathcal{L}(\widetilde{Z}_{k0})
=O\left(\frac{1}{R\gamma_k^3|\log\varepsilon_k|}
\right),
\endaligned
$$
and for $\gamma_k(R+1)<|y-\xi'_k|\leq
1/(\varepsilon_0t^{2\beta})$,
$$
\aligned
\mathcal{L}(\widetilde{Z}_{k0})=
\mathcal{L}(\widehat{Z}_{k0})=O
\left(\frac{\log|y-\xi'_k|-
\log(R\gamma_k)}{\big(1+\big|\frac{y-\xi'_k}{\gamma_{k}}\big|^2\big)^2}\cdot
\frac{1}{\gamma_k^3|\log\varepsilon_k|}
\right).
\endaligned
$$
These, together with the estimate of $\widehat{Z}_{p}$ in (\ref{3.24}),
give
$$
\aligned
\int_{\Omega_{3,k}}\widehat{Z}_{p}\mathcal{L}(\widetilde{Z}_{k0})
=O\left(\frac{\varepsilon_0\log t}{\rho_0v_0\gamma_k|\log\varepsilon_0||\log\varepsilon_k|}\right).
\endaligned
$$
As on $\partial\Omega_{3,k}$, by (\ref{2.2}) and (\ref{3.24}),
$$
\aligned
\widehat{Z}_{p}
=O\left(\frac{\varepsilon_0\log t}{\rho_0v_0|\log\varepsilon_0|}\right),
\,\qquad\,\,\,\qquad\,
|\nabla\widehat{Z}_{p}|
=O\left(\frac{\varepsilon_0^2t^{\beta}}{\rho_0v_0|\log\varepsilon_0|}\right),
\endaligned
$$
and
$$
\aligned
\widehat{Z}_{k0}
=O\left(\frac{\log t}{\gamma_k|\log\varepsilon_k|}\right),
\,\qquad\,\,\,\qquad\,
|\nabla\widehat{Z}_{k0}|
=O\left(\frac{\varepsilon_0  t^{2\beta} }{\gamma_k|\log\varepsilon_k|}\right).
\endaligned
$$
Then
$$
\aligned
\int_{\Omega_{3,k}}\widetilde{Z}_{k0}\mathcal{L}(\widetilde{Z}_{p})
=O\left(\frac{\varepsilon_0\log t}{\rho_0v_0\gamma_k|\log\varepsilon_0||\log\varepsilon_k|}\right).
\endaligned
$$
By the above estimates, we readily  have
\begin{equation}\label{3.60}
\aligned
\int_{\Omega_t}\widetilde{Z}_{k0}\mathcal{L}(\widetilde{Z}_{p})
=O\left(\frac{\varepsilon_0\log^2 t}{\rho_0v_0\gamma_k|\log\varepsilon_0||\log\varepsilon_k|}\right),
\,\quad\,k=1,\ldots,m.
\endaligned
\end{equation}
Inserting estimates (\ref{3.59}) and (\ref{3.60})  into  (\ref{3.50}), we get
\begin{equation}\label{3.61}
\aligned
\frac{\varepsilon_0|d_p|}{\rho_0v_0|\log\varepsilon_0|}
\leq
C\|h\|_{*}
+\frac{C\log^2 t}{\,|\log\varepsilon_0|\,}\left(
\frac{\varepsilon_0|d_p|}{\rho_0v_0|\log\varepsilon_0|}
+
\sum_{k=1}^m\frac{|d_k|}{\gamma_k|\log\varepsilon_k|}
\right).
\endaligned
\end{equation}
On the other hand, similar to the above arguments in (\ref{3.59})-(\ref{3.60}), we can show that
for $R$ and $t$ large enough,
\begin{equation}\label{3.62}
\aligned
\int_{\Omega_t}\widetilde{Z}_{i0}\mathcal{L}(\widetilde{Z}_{i0})
=\frac{2\pi}{\gamma_i^2|\log\varepsilon_i|}\left[1+O\left(\frac1{R^2}\right)\right],
\endaligned
\end{equation}
and
\begin{equation}\label{3.63}
\aligned
\int_{\Omega_t}\widetilde{Z}_{k0}\mathcal{L}(\widetilde{Z}_{i0})
=O\left(\frac{\log^2t}{\gamma_i\gamma_k|\log\varepsilon_i||\log\varepsilon_k|}
\right)
\quad\quad\textrm{for all}\,\,\,k\neq i.
\endaligned
\end{equation}
These, together with (\ref{3.51}) and (\ref{3.60}), imply
\begin{equation}\label{3.64}
\aligned
\frac{|d_i|}{\gamma_i|\log\varepsilon_i|}
\leq
C\|h\|_{*}
+\frac{C\log^2t}{|\log\varepsilon_i|}
\left(
\frac{\varepsilon_0|d_p|}{\rho_0v_0|\log\varepsilon_0|}
+\sum_{k=1}^m\frac{|d_k|}{\gamma_k|\log\varepsilon_k|}
\right).
\endaligned
\end{equation}
As a result,  using linear algebra arguments, by (\ref{2.7}), (\ref{3.61})  and (\ref{3.64}) we can prove  Claim 2 for
$d_p$ and $d_i$,  and then
complete the proof by (\ref{3.30}).

\vspace{1mm}
\vspace{1mm}
\vspace{1mm}

{\bf Step 5:}
Proof of Proposition 3.1. We begin by establishing the validity of the a priori
estimate (\ref{3.8}).
Using estimate (\ref{3.20}) and the fact that $\|\chi_iZ_{ij}\|_{*}=O(\gamma_i)$, we deduce
\begin{equation}\label{3.65}
\aligned
\|\phi\|_{L^{\infty}(\Omega_t)}\leq Ct\left(
\|h\|_{*}+\sum\limits_{i=1}^m\sum\limits_{j=1}^2\gamma_i|c_{ij}|
\right).
\endaligned
\end{equation}
So it suffices to estimate the values of the constants $c_{ij}$.
Let us consider the cut-off function $\eta_{i2}$ defined in
 (\ref{3.26}).   Multiplying   (\ref{3.1}) by $\eta_{i2}Z_{ij}$  and integrating
 by parts, we find
\begin{equation}\label{3.66}
\aligned
\int_{\Omega_t}
\phi \mathcal{L}(\eta_{i2}Z_{ij})=
\int_{\Omega_t}
h\eta_{i2}Z_{ij}+
\sum_{k=1}^m\sum_{l=1}^2c_{kl}
\int_{\Omega_t}\chi_kZ_{kl}\eta_{i2}Z_{ij}.
\endaligned
\end{equation}
Notice that
$$
\aligned
\mathcal{L}(\eta_{i2}Z_{ij})=&
\eta_{i2}\mathcal{L}(Z_{ij})
-Z_{ij}\Delta \eta_{i2}
-2\nabla\eta_{i2}\nabla Z_{ij}
=\left[\frac{1}{\gamma_{i}^2}\frac{8}{\big(1+\big|\frac{y-\xi'_i}{\gamma_{i}}\big|^2\big)^2}-W\right]\eta_{i2}Z_{ij}
+O\left(\varepsilon_0^3\right).
\endaligned
$$
For the estimate of the first term, we decompose  $\supp(\eta_{i2})$ to some subregions:
$$
\aligned
\widehat{\Omega}_{p}=\supp(\eta_{i2})\bigcap\big\{|y-p'|\leq
1/(\varepsilon_0 t^{2\beta})\big\},
\,\quad\quad\quad\,
\widehat{\Omega}_{k1}=\supp(\eta_{i2})\bigcap\big\{|y-\xi'_k|\leq
1/(\varepsilon_0 t^{2\beta})\big\},
\,\,\,\,k=1,\ldots,m,
\endaligned
$$
$$
\aligned
\widehat{\Omega}_{2}=\supp(\eta_{i2})\setminus\left[\,\bigcup_{k=1}^m\widehat{\Omega}_{k1}
\cup\widehat{\Omega}_{p}\right],
\endaligned
$$
where $\supp(\eta_{i2})=\left\{|y-\xi'_i|\leq
6d/\varepsilon_0\right\}$.
Notice that, by (\ref{2.2}),
\begin{equation}\label{3.67}
\aligned
|y-\xi'_i|\geq |\xi_i'-p'|-|y-p'|\geq
|\xi_i'-p'|-\frac{1}{\varepsilon_0 t^{2\beta}}
\geq
\frac{1}{\varepsilon_0 t^\beta}\left(1-\frac{1}{t^\beta}
\right)
\endaligned
\end{equation}
uniformly in $\widehat{\Omega}_{p}$, and
\begin{equation}\label{3.68}
\aligned
|y-\xi'_i|\geq |\xi_i'-\xi_k'|-|y-\xi_k'|\geq
|\xi_i'-\xi_k'|-\frac{1}{\varepsilon_0 t^{2\beta}}
\geq
\frac{1}{\varepsilon_0 t^\beta}\left(1-\frac{1}{t^\beta}
\right)
\endaligned
\end{equation}
uniformly in $\widehat{\Omega}_{k1}$ with $k\neq i$.
By  (\ref{2.23}), (\ref{2.28}), (\ref{2.30}), (\ref{2.32}) and   (\ref{3.5}) we have that
in $\widehat{\Omega}_{i1}$,
$$
\aligned
\left[\frac{1}{\gamma_{i}^2}\frac{8}{\big(1+\big|\frac{y-\xi'_i}{\gamma_{i}}\big|^2\big)^2}-W\right]
\eta_{i2}Z_{ij}=&
\frac{1}{\gamma_{i}^2}\frac{8}{\big(1+\big|\frac{y-\xi'_i}{\gamma_{i}}\big|^2\big)^{5/2}}
\left[
O\left(\varepsilon_0t^{\beta}\left|\frac{y-\xi'_i}{\gamma_{i}}\right|\right)
+
O\left(\frac{\varepsilon_i^2\mu_i^2}{\gamma_i}\right)+
O\left(\frac{\varepsilon_0^2\mu_0^2t^{2\beta(1+\alpha)}}{\gamma_i}\right)\right.\\
&\left.
+\sum_{j=1,\,j\neq i}^m O\left(\frac{1}{\gamma_i}\varepsilon_j^2\mu_j^2t^{2\beta}\right)
\right],
\endaligned
$$
and in $\widehat{\Omega}_{p}$,  by (\ref{3.67}),
$$
\aligned
\left[\frac{1}{\gamma_{i}^2}\frac{8}{\big(1+\big|\frac{y-\xi'_i}{\gamma_{i}}\big|^2\big)^2}-W\right]
\eta_{i2}Z_{ij}=
\left[
O\left(\frac{\gamma_i^2}{|y-\xi'_i|^4}\right)+
\left(\frac{\varepsilon_0}{\rho_0v_0}\right)^2
O\left(\frac{8(1+\alpha)^2\big|\frac{\varepsilon_0 y-p}{\rho_0v_0}\big|^{2\alpha}}{\,\big(1+\big|\frac{\varepsilon_0 y-p}{\rho_0v_0}\big|^{2(1+\alpha)}\big)^2\,}\right)
\right]O\left(\frac1{|y-\xi'_i|}
\right),
\endaligned
$$
and in $\widehat{\Omega}_{k1}$, $k\neq i$, by (\ref{3.68}),
$$
\aligned
\left[\frac{1}{\gamma_{i}^2}\frac{8}{\big(1+\big|\frac{y-\xi'_i}{\gamma_{i}}\big|^2\big)^2}-W\right]
\eta_{i2}Z_{ij}=\left[
O\left(\frac{\gamma_i^2}{|y-\xi'_i|^4}\right)+
O\left(\frac{1}{\gamma_{k}^2}\frac{8}{\big(1+\big|\frac{y-\xi'_k}{\gamma_{k}}\big|^2\big)^2}\right)
\right]O\left(\frac1{|y-\xi'_i|}
\right),
\endaligned
$$
and in $\widehat{\Omega}_{2}$,
$$
\aligned
\left[\frac{1}{\gamma_{i}^2}\frac{8}{\big(1+\big|\frac{y-\xi'_i}{\gamma_{i}}\big|^2\big)^2}-W\right]
\eta_{i2}Z_{ij}=
O\left(\varepsilon_0^3\varepsilon_i^2\mu_i^2t^{10\beta}\right)
+O\left(\varepsilon_0^3 t^{2\beta(4m+5+2\alpha)} e^{-t\phi_1(\varepsilon_0 y)}
\right).
\endaligned
$$
Then
\begin{equation}\label{3.69}
\aligned
\left|\int_{\Omega_t}
\phi \mathcal{L}(\eta_{2i}Z_{1i})\right|
\leq C\varepsilon_0t^\beta \|\phi\|_{L^{\infty}(\Omega_t)}.
\endaligned
\end{equation}
On the other hand, since
$\|\eta_{i2}Z_{ij}\|_{L^{\infty}(\Omega_t)}\leq C\gamma_i^{-1}$, we know that
\begin{equation}\label{3.70}
\aligned
\int_{\Omega_t}h\eta_{i2}Z_{ij}=O\left(\frac{\|h\|_{*}}{\gamma_i}\right).
\endaligned
\end{equation}
Moreover,  if $k= i$, by (\ref{3.2}), (\ref{3.5}) and (\ref{3.6}),
\begin{equation}\label{3.71}
\aligned
\int_{\Omega_t}\chi_kZ_{kl}\eta_{k2}Z_{kj}=
\int_{\mathbb{R}^2}\chi(|z|)\mathcal{Z}_{l}(z) \mathcal{Z}_{j}(z)dz=C\delta_{lj},
\endaligned
\end{equation}
while if $k\neq i$, by (\ref{3.68}),
\begin{equation}\label{3.72}
\aligned
\int_{\Omega_t}\chi_kZ_{kl}\eta_{i2}Z_{ij}=O\left(\gamma_k\varepsilon_0 t^\beta\right).
\endaligned
\end{equation}
As a consequence, substituting estimates (\ref{3.69})-(\ref{3.72}) into  (\ref{3.66}),  we find
$$
\aligned
|c_{ij}|\leq C\left(\varepsilon_0t^\beta \|\phi\|_{L^{\infty}(\Omega_t)}+
\frac{1}{\gamma_i}\|h\|_{*}+
\sum\limits_{k\neq i}^m\sum\limits_{l=1}^2\gamma_k\varepsilon_0 t^\beta|c_{kl}|
\right),
\endaligned
$$
and then, by   (\ref{2.23}),
$$
\aligned
|c_{ij}|\leq C \left(\varepsilon_0t^\beta \|\phi\|_{L^{\infty}(\Omega_t)}+
\frac{1}{\gamma_i}\|h\|_{*}
\right).
\endaligned
$$
Combing this estimate with (\ref{3.65}), we conclude
\begin{equation}\label{3.73}
\aligned
|c_{ij}|\leq C\frac{1}{\gamma_i}\|h\|_{*},
\endaligned
\end{equation}
which proves (\ref{3.8}).

Now, we consider the Hilbert space
$$
\aligned
H_{\xi}=\left\{\phi\in H_0^1(\Omega_t)\left|\,
\int_{\Omega_t}\chi_iZ_{ij}\phi=0\,\,\,\,\,\,\forall\,\,\,i=1,\ldots,m,\,\,j=1,2\right.\right\}
\endaligned
$$
with the norm
$\|\phi\|_{H_\xi}=\|\nabla\phi\|_{L^{\infty}(\Omega_t)}$.
Equation (\ref{3.1}) is equivalent to find
$\phi\in H_\xi$, such that
$$
\aligned
\int_{\Omega_t}\nabla\phi\nabla\psi-\int_{\Omega_t}W\phi\psi
=\int_{\Omega_t}h\psi,\,\,\quad\,\,\forall
\,\,\psi\in H_\xi.
\endaligned
$$
By Fredholm's alternative this is equivalent to the uniqueness of solutions to this
problem, which in turn follows from  estimate (\ref{3.8}).
\end{proof}

The result of Proposition 3.1 implies that the unique solution $\phi=T(h)$
of (\ref{3.1}) defines a bounded linear map from the Banach space
$\mathcal{C}_*$ of all functions $h$ in $L^\infty$ for which $\|h\|_{*}<\infty$,
into $L^\infty$.

\vspace{1mm}
\vspace{1mm}
\vspace{1mm}

\noindent{\bf Lemma 3.2.}\,\,{\it
For any integer $m\geq1$, the
operator $T$ is differentiable  with
respect to the variables $\xi=(\xi_1,\ldots,\xi_m)$ in $\mathcal{O}_t$,
precisely for any $k=1,\ldots,m$ and $l=1,2$,
\begin{equation}\label{3.74}
\aligned
\|\partial_{\xi'_{kl}}T(h)\|_{L^{\infty}(\Omega_t)}\leq Ct^2\|h\|_{*}.
\endaligned
\end{equation}
}
\begin{proof}
Differentiating (\ref{3.1}) with respect to $\xi'_{kl}$, formally
$Z=\partial_{\xi'_{kl}}\phi$ should satisfy
$$
\aligned
\left\{
\aligned
&\mathcal{L}(Z)=\phi\,\partial_{\xi'_{kl}}W+
\sum\limits_{i=1}^m\sum\limits_{j=1}^2\left[c_{ij}\partial_{\xi'_{kl}}(\chi_iZ_{ij})+\widetilde{c}_{ij}\chi_iZ_{ij}\right]\,\,\quad
\ \textrm{in}\,\,\,\,\,\Omega_t,\\
&Z=0
\quad\qquad\qquad\qquad\quad
\,\quad\quad\,\,\,\,
\qquad\qquad\qquad\quad
\qquad\qquad\,\,
\textrm{on}
\,\,\,\po_t,\\[1mm]
&\int_{\Omega_t}\chi_iZ_{ij}Z=-\int_{\Omega_t}\phi\partial_{\xi'_{kl}}(\chi_iZ_{ij})
\qquad\quad\,
\forall\,\,\,i=1,\ldots,m,
\,\,j=1,2,
\endaligned
\right.
\endaligned
$$
where (still formally) $\widetilde{c}_{ij}=\partial_{\xi'_{kl}}(c_{ij})$.
Furthermore, if we consider the constants $b_{ij}$ defined as
$$
\aligned
b_{ij}\int_{\Omega_t}\chi^2_i|Z_{ij}|^2=\int_{\Omega_t}\phi\,\partial_{\xi'_{kl}}(\chi_iZ_{ij}),
\endaligned
$$
and set
$$
\aligned
\widetilde{Z}=Z+\sum_{i=1}^m\sum_{j=1}^2b_{ij}\chi_iZ_{ij},
\endaligned
$$
then we have
$$
\aligned
\left\{
\aligned
&\mathcal{L}(\widetilde{Z})=f+\sum\limits_{i=1}^m\sum\limits_{j=1}^2\widetilde{c}_{ij}\chi_iZ_{ij}
\quad\quad\,\,\,\,
\textrm{in}\,\,\,\,\,\,\,\Omega_t,\\
&\widetilde{Z}=0
\,\quad\qquad\qquad\qquad\qquad\qquad
\,\,\,\,\,\textrm{on}
\,\,\,\,\,\po_t,\\[1mm]
&\int_{\Omega_t}\chi_iZ_{ij}\widetilde{Z}=0\,\,\,\,\,
\,\,\forall\,\,i=1,\ldots,m,\,\,j=1,2,
\endaligned
\right.
\endaligned
$$
where
$$
\aligned
f=\phi\,\partial_{\xi'_{kl}}W
+\sum\limits_{i=1}^m\sum\limits_{j=1}^2
b_{ij}\mathcal{L}(\chi_iZ_{ij})
+
\sum\limits_{i=1}^m
\sum\limits_{j=1}^2c_{ij}\partial_{\xi'_{kl}}(\chi_iZ_{ij}).
\endaligned
$$
From Proposition 3.1 it follows that this equation  has a unique solution
$\widetilde{Z}$ and $\widetilde{c}_{ij}$, and hence
$\partial_{\xi'_{kl}}T(h)=T(f)-\sum_{i=1}^m\sum_{j=1}^2b_{ij}\chi_iZ_{ij}$
is well defined. Moreover, by (\ref{3.8}) we get
\begin{equation}\label{3.75}
\aligned
\|\partial_{\xi'_{kl}}T(h)\|_{L^{\infty}(\Omega_t)}\leq\|T(f)\|_{L^{\infty}(\Omega_t)}
+C\sum_{i=1}^m\sum_{j=1}^2\frac{1}{\gamma_i}|b_{ij}|
\leq C t\|f\|_{*}+C\sum_{i=1}^m\sum_{j=1}^2\frac{1}{\gamma_i}|b_{ij}|.
\endaligned
\end{equation}
Now, to prove estimate (\ref{3.74}), we
 first estimate $\partial_{\xi'_{kl}}W$.
Notice that $\partial_{\xi'_{kl}}W=W\partial_{\xi'_{kl}}V$. Obviously,
by (\ref{2.28}), (\ref{2.30}), (\ref{2.32}) and (\ref{3.7}) we  find
$\|W\|_{*}=O\left(1\right)$. On the other hand, similar to the proof of Lemma 2.1,
by (\ref{2.14})-(\ref{2.15}) we can  compute that
\begin{equation}\label{3.76}
\aligned
\partial_{\xi'_{kl}}H_0(\varepsilon_0y)=O\left(\varepsilon_0t^\beta\right)
\quad\qquad
\textrm{and}
\quad\qquad
\partial_{\xi'_{kl}}H_i(\varepsilon_0y)=O\left(\varepsilon_0t^\beta\right),
\quad\,i=1,\ldots,m,
\endaligned
\end{equation}
uniformly in $\overline{\Omega}_t$.
Furthermore,  by (\ref{2.4}),
(\ref{2.8}), (\ref{2.14}), (\ref{2.15}),
(\ref{2.23}), (\ref{3.2}) and (\ref{3.5})
we can directly check that
\begin{equation}\label{3.77}
\aligned
\partial_{\xi'_{kl}}V(y)=Z_{kl}(y)+O\left(\varepsilon_0t^\beta\right).
\endaligned
\end{equation}
This, together with the fact that $\frac{1}{\gamma_k}\leq C$ uniformly on $t$, immediately implies
\begin{equation}\label{3.78}
\aligned
\|\partial_{\xi'_{kl}}V\|_{L^{\infty}(\Omega_t)}=O\left(1\right)
\,\,\qquad\,\,\textrm{and}\,\,\qquad\,\,\|\partial_{\xi'_{kl}}W\|_{*}=O\left(1\right).
\endaligned
\end{equation}
Next, by definitions (\ref{3.5})-(\ref{3.6}),
a straightforward  computation gives
\begin{equation}\label{3.79}
\aligned
\|\partial_{\xi'_{kl}}(\chi_iZ_{ij})\|_{*}=\left\{
\aligned
&O\left(\varepsilon_0\gamma_it^\beta\right)
\,\quad\quad\,\,\textrm{if}\,\,\,\,i\neq k,\\[1mm]
&O\left(1\right)\,\,\quad
\,\qquad\quad\,\textrm{if}
\,\,\,\,i=k.
\endaligned
\right.
\endaligned
\end{equation}
Furthermore,
\begin{equation}\label{3.80}
\aligned
|b_{ij}|=\left\{
\aligned
&O\left(\varepsilon_0\gamma_it^\beta\right)\|\phi\|_{L^{\infty}(\Omega_t)}
\,\,\,\,\,\quad\,\textrm{if}\,\,\,\,i\neq k,\\[1mm]
&O\left(1\right)\|\phi\|_{L^{\infty}(\Omega_t)}\,\,\quad
\,\,\,\ \,\,\,\quad\,\,\,\textrm{if}
\,\,\,\,i=k.
\endaligned
\right.
\endaligned
\end{equation}
Finally, by  (\ref{3.8}),
(\ref{3.33}), (\ref{3.73}), (\ref{3.79})
and (\ref{3.80}), we obtain
\begin{equation}\label{}
\aligned
\|f\|_{*}\leq Ct\|h\|_{*}
\,\,\,\,\qquad\,\,\,\,
\textrm{and}
\,\,\,\,\qquad\,\,\,\,
|b_{ij}|\leq Ct\|h\|_{*}.
\endaligned
\end{equation}
Inserting these into (\ref{3.75}), we then prove  (\ref{3.74}).
\end{proof}

\vspace{1mm}
\vspace{1mm}
\vspace{1mm}

\section{The  nonlinear projected problem}
In this section we   solve  the
nonlinear projected problem: for any integer $m\geq1$
and
any points
$\xi=(\xi_1,\ldots,\xi_m)\in\mathcal{O}_t$,
we find a function $\phi$ and scalars $c_{ij}$,
$i=1,\ldots,m$, $j=1,2$,  such that
\begin{equation}\label{4.1}
\aligned
\left\{
\aligned
&\mathcal{L}(\phi)=-\Delta\phi-W\phi=E+N(\phi)+\sum\limits_{i=1}^m\sum\limits_{j=1}^2c_{ij}\chi_iZ_{ij}
\quad
\textrm{in}\,\,\,\,\,\Omega_t,\\
&\,\phi=0
\quad\quad\quad\quad\qquad
\qquad\qquad\qquad\qquad\,
\quad\quad\quad\quad\qquad
\textrm{on}
\,\,\,\po_t,\\[1mm]
&\int_{\Omega_t}\chi_iZ_{ij}\phi=0
\,\,\qquad\qquad\qquad\qquad\,\,\,
\forall\,\,\,i=1,\ldots,m,
\,\,j=1,2,
\endaligned
\right.
\endaligned
\end{equation}
where  $W$ is as in (\ref{2.28}), (\ref{2.30}) and (\ref{2.32}),
and $E$, $N(\phi)$ are given by
(\ref{2.21}) and (\ref{2.35}), respectively.

\vspace{1mm}
\vspace{1mm}
\vspace{1mm}
\vspace{1mm}

\noindent{\bf Proposition 4.1.}\,\,{\it
Let $m$ be a positive   integer.
Then there exist constants $t_m>1$ and $C>0$ such
that for any $t>t_m$ and any points
$\xi=(\xi_1,\ldots,\xi_m)\in\mathcal{O}_t$,
problem {\upshape(\ref{4.1})} admits
a unique solution
$\phi\in L^{\infty}(\Omega_t)$, and
scalars $c_{ij}\in\mathbb{R}$,
$i=1,\ldots,m$, $j=1,2$, such that
\begin{eqnarray}\label{4.2}
\quad\quad
\|\phi\|_{L^{\infty}(\Omega_t)}\leq C t
\max\big\{(\rho_0v_0)^{\min\{1,2(\alpha-\hat{\alpha})\}}t^{\beta},\,\varepsilon_0\gamma_1t^{\beta},\,\ldots,\,\varepsilon_0\gamma_mt^{\beta},\,
\|e^{-\frac{1}{2}t\phi_1}\|_{L^{\infty}(\Omega_t)}\big\}.
\end{eqnarray}
Furthermore, the map
$\xi'\mapsto\phi(\xi')\in C(\overline{\Omega}_t)$
is $C^1$, precisely for any $k=1,\ldots,m$ and  $l=1,2$,
\begin{eqnarray}\label{4.3}
\quad\quad
\|\partial_{\xi'_{kl}}\phi\|_{L^{\infty}(\Omega_t)}
\leq C t^{2}\max\big\{(\rho_0v_0)^{\min\{1,2(\alpha-\hat{\alpha})\}}t^{\beta},\,\varepsilon_0\gamma_1t^{\beta},\,\ldots,\,\varepsilon_0\gamma_mt^{\beta},\,
\|e^{-\frac{1}{2}t\phi_1}\|_{L^{\infty}(\Omega_t)}\big\},
\end{eqnarray}
where $\xi':=(\xi'_1,\ldots,\xi'_m)=(\frac1{\varepsilon_0}\xi_1,\ldots,\frac1{\varepsilon_0}\xi_m)$.
}

\vspace{1mm}
\vspace{1mm}

\begin{proof}
Let  $T$ be the  operator as defined in Proposition 3.1. Then $\phi$ solves (\ref{4.1})
if and only if
\begin{equation}\label{4.4}
\aligned
\phi=T(E+N(\phi))\equiv A(\phi).
\endaligned
\end{equation}
For a given number $\kappa>0$, let us consider the region
$$
\aligned
\mathcal{F}_\kappa=\left\{\phi\in C(\overline{\Omega}_t)\left|\,
\|\phi\|_{L^{\infty}(\Omega_t)}\leq \kappa t
\max\big\{(\rho_0v_0)^{\min\{1,2(\alpha-\hat{\alpha})\}}t^{\beta},\,\varepsilon_0\gamma_1t^{\beta},\,\ldots,\,\varepsilon_0\gamma_mt^{\beta},\,
\|e^{-\frac{1}{2}t\phi_1}\|_{L^{\infty}(\Omega_t)}\big\}\right.\right\}.
\endaligned
$$
Observe that, by (\ref{2.29}), (\ref{2.31}),  (\ref{2.33}) and (\ref{3.7}),
\begin{equation}\label{4.5}
\aligned
\|E\|_{*}\leq  C \max\big\{(\rho_0v_0)^{\min\{1,2(\alpha-\hat{\alpha})\}}t^{\beta},\,\varepsilon_0\gamma_1t^{\beta},\,\ldots,\,\varepsilon_0\gamma_mt^{\beta},\,
\|e^{-\frac{1}{2}t\phi_1}\|_{L^{\infty}(\Omega_t)}\big\}.
\endaligned
\end{equation}
Moreover, by definition (\ref{2.35}) of $N(\phi)$   and Lagrange's theorem
we have that for  $\phi$, $\phi_1$, $\phi_2\in\mathcal{F}_\kappa$,
$$
\aligned
\|N(\phi)\|_{*}\leq C\|W\|_{*}\|\phi\|^2_{L^{\infty}(\Omega_t)}\leq C\|\phi\|^2_{L^{\infty}(\Omega_t)},
\endaligned
$$
$$
\aligned
\|N(\phi_1)-N(\phi_2)\|_{*}
\leq C\big(\max\limits_{l=1,2}\|\phi_l\|_{L^{\infty}(\Omega_t)}\big)
\|\phi_1-\phi_2\|_{L^{\infty}(\Omega_t)}.
\endaligned
$$
where $C>0$ is independent of $\kappa$ and $t$.
Hence by (\ref{2.7}), (\ref{2.14}), (\ref{2.15}), (\ref{2.23})
and Proposition 3.1,
$$
\aligned
\|A(\phi)\|_{L^{\infty}(\Omega_t)}
\leq Ct \big(\|E\|_{*}+\|N(\phi)\|_{*}\big)
\leq C t\max\big\{&(\rho_0v_0)^{\min\{1,2(\alpha-\hat{\alpha})\}}t^{\beta},\,\varepsilon_0\gamma_1t^{\beta},\,\ldots,\,\varepsilon_0\gamma_mt^{\beta},\,
\|e^{-\frac{1}{2}t\phi_1}\|_{L^{\infty}(\Omega_t)}\big\},\\
\|A(\phi_1)-A(\phi_2)\|_{L^{\infty}(\Omega_t)}\leq C&t\|N(\phi_1)-N(\phi_2)\|_{*}
< \frac12\|\phi_1-\phi_2\|_{L^{\infty}(\Omega_t)}.
\endaligned
$$
This means  that
 for all $t$  large enough,
$A$ is a contraction on $\mathcal{F}_\kappa$ and thus a unique
fixed point of $A$ exists in the region.

We now analyze the differentiability of the map $\xi'\mapsto\phi$.
Assume for instance that the
partial derivative $\partial_{\xi'_{kl}}\phi$ exists.
Then,
formally
$$
\aligned
\partial_{\xi'_{kl}}\phi=
T\left(\partial_{\xi'_{kl}}E+\partial_{\xi'_{kl}}N(\phi)\right)
+\left(\partial_{\xi'_{kl}}T\right)\left(E+N(\phi)\right).
\endaligned
$$
By (\ref{3.74}) and (\ref{4.5}), we get
$$
\aligned
\left\|\left(\partial_{\xi'_{kl}}T\right)\left(E+N(\phi)\right)\right\|_{L^{\infty}(\Omega_t)}
&\leq Ct^2\left(
\|E\|_{*}+\|N(\phi)\|_{*}\right)\\
&\leq Ct^{2}\max\big\{(\rho_0v_0)^{\min\{1,2(\alpha-\hat{\alpha})\}}t^{\beta},\,\varepsilon_0\gamma_1t^{\beta},\,\ldots,\,\varepsilon_0\gamma_mt^{\beta},\,
\|e^{-\frac{1}{2}t\phi_1}\|_{L^{\infty}(\Omega_t)}\big\}.
\endaligned
$$
Observe that
$$
\aligned
\partial_{\xi'_{kl}}N(\phi)=\partial_{\xi'_{kl}}W\left(
e^{\phi}-\phi-1
\right)
+W\left(
e^\phi-1
\right)\partial_{\xi'_{kl}}\phi,
\endaligned
$$
so that, by (\ref{3.78}),
$$
\aligned
\|\partial_{\xi'_{kl}}N(\phi)\|_{*}
\leq C\|\phi\|_{L^{\infty}(\Omega_t)}\left(
\|\phi\|_{L^{\infty}(\Omega_t)}+\|\partial_{\xi'_{kl}}\phi\|_{L^{\infty}(\Omega_t)}\right).
\endaligned
$$
Also,  thanks to the  expansion of $\partial_{\xi'_{kl}}V$ in (\ref{3.77}) ,
by (\ref{2.22}), (\ref{2.28}), (\ref{2.30}) and (\ref{2.32}) we can directly check  that
\begin{equation}\label{4.6}
\aligned
\|\partial_{\xi'_{kl}}E\|_{*}
\leq  C \max\big\{(\rho_0v_0)^{\min\{1,2(\alpha-\hat{\alpha})\}}t^{\beta},\,\varepsilon_0\gamma_1t^{\beta},\,\ldots,\,\varepsilon_0\gamma_mt^{\beta},\,
\|e^{-\frac{1}{2}t\phi_1}\|_{L^{\infty}(\Omega_t)}\big\}.
\endaligned
\end{equation}
Hence by Proposition 3.1,  we then  prove
$$
\aligned
\|\partial_{\xi'_{kl}}\phi\|_{L^{\infty}(\Omega_t)}\leq
C t^{2}\max\big\{(\rho_0v_0)^{\min\{1,2(\alpha-\hat{\alpha})\}}t^{\beta},\,\varepsilon_0\gamma_1t^{\beta},\,\ldots,\,\varepsilon_0\gamma_mt^{\beta},\,
\|e^{-\frac{1}{2}t\phi_1}\|_{L^{\infty}(\Omega_t)}\big\}.
\endaligned
$$
The above computation can be  made rigorous by using the implicit function theorem
and the fixed point representation (\ref{4.4}) which guarantees $C^1$
regularity of $\xi'$.
\end{proof}

\section{The reduced problem: A maximization procedure }
In this section we study a maximization problem involving the variational reduction. Let
us consider the energy function $J_t$ associated to problem
(\ref{1.4}), namely
\begin{equation}\label{5.1}
\aligned
J_t(u)=\frac12\int_{\Omega}
|\nabla u|^2-\int_{\Omega}|x-p|^{2\alpha}k(x)e^{-t\phi_1}e^u,
\,\quad\,\, u\in H_0^1(\Omega).
\endaligned
\end{equation}
For any integer $m\geq1$,
we take its finite dimensional restriction
\begin{equation}\label{5.2}
\aligned
F_t(\xi)=J_t\big(U(\xi)+\tilde{\phi}(\xi)
\big)\,
\qquad\,\forall\,\,\,\xi=(\xi_1,\ldots,\xi_m)\in\overline{\mathcal{O}}_t,
\endaligned
\end{equation}
where $U(\xi)$ is our approximate solution defined in (\ref{2.8})
and $\tilde{\phi}(\xi)(x)=\phi(\frac{x}{\varepsilon_0},\frac{\xi}{\varepsilon_0})$,
$x\in\Omega$, with $\phi=\phi_{\xi'}$  the unique solution
to problem (\ref{4.1}) given by Proposition 4.1. Define
\begin{equation}\label{5.3}
\aligned
\mathcal{M}_m^t=\max\limits_{(\xi_1,\ldots,\xi_m)\in\overline{\mathcal{O}}_t}
F_t(\xi_1,\ldots,\xi_m).
\endaligned
\end{equation}
From
 the results obtained in Proposition 4.1 and the
definition of function $U(\xi)$ we have clearly that
for any integer $m\geq1$, the map
$F_t:\overline{\mathcal{O}}_t\rightarrow\mathbb{R}$
is of class $C^1$ and then this maximization problem
has a solution over $\overline{\mathcal{O}}_t$.

\vspace{1mm}
\vspace{1mm}
\vspace{1mm}
\vspace{1mm}

\noindent{\bf Proposition 5.1.}\,\,\,{\it
For any integer $m\geq1$ and any $t$ large enough, the maximization problem
\begin{eqnarray}\label{5.4}
\max\limits_{(\xi_1,\ldots,\xi_m)\in\overline{\mathcal{O}}_t}
F_t(\xi_1,\ldots,\xi_m)
\end{eqnarray}
has a solution $\xi_t=(\xi_{1,t},\ldots,\xi_{m,t})\in\mathcal{O}_t^o$, i.e., the interior of $\mathcal{O}_t$.
}

\begin{proof}
The proof of this result consists of  three steps which we state and prove next.

\vspace{1mm}

{\bf Step 1:}
With the choices  for  the parameters $\mu_0$ and $\mu_i$,
$i=1,\ldots,m$,
respectively given by {\upshape(\ref{2.12})} and {\upshape(\ref{2.13})},
let us prove that the following expansion   holds
\begin{eqnarray}\label{5.5}
J_t\big(U(\xi)\big)=
8\pi(1+\alpha)t
+8\pi t\sum_{i=1}^m
\phi_1(\xi_i)
+16\pi(2+\alpha)\sum_{i=1}^m\log|\xi_i-p|
+16\pi\sum_{i\neq j}^m
\log|\xi_i-\xi_j|+O\left(1\right)
\,\,\,
\end{eqnarray}
uniformly for all points $\xi=(\xi_1,\ldots,\xi_m)\in\mathcal{O}_t$ and for all $t$ large enough.

Observe first  that by (\ref{2.8}) and (\ref{2.9}),
\begin{eqnarray}\label{5.6}
\frac12\int_{\Omega}|\nabla U|^2
=-\frac12\int_{\Omega} U_0\Delta U_0-\sum_{j=1}^m\int_{\Omega} U_j\Delta U_0-\frac12\sum_{i,j=1}^m\int_{\Omega} U_j\Delta U_i
\qquad\qquad\qquad\qquad\quad\,\,\,\,\,
&&\nonumber\\
=-\frac12\int_{\Omega}\big(u_0+H_0\big)\Delta u_0
-\sum_{j=1}^m\int_{\Omega}\big(u_j+H_j\big)\Delta u_0
-\frac12\sum_{i,j=1}^m\int_{\Omega}\big(u_j+H_j\big)\Delta u_i.
&&
\end{eqnarray}
Let us analyze the behavior of the first term.
By (\ref{2.4}), (\ref{2.5}) and  (\ref{2.10}) we get
$$
\aligned
-\int_{\Omega}\big(u_0+H_0\big)\Delta u_0
=\int_{\Omega}\frac{8\varepsilon_0^2\mu_0^2(1+\alpha)^2|x-p|^{2\alpha}}{\,(\varepsilon_{0}^2\mu_{0}^2+|x-p|^{2(1+\alpha)})^2\,}
\left[\log
\frac{1}{(\varepsilon_{0}^2\mu_{0}^2+|x-p|^{2(1+\alpha)})^2}+
(1+\alpha)H(x,p)+O\left(\varepsilon_0^2\mu_0^2\right)
\right].
\endaligned
$$
Making the change of
variables $\rho_0v_0z=x-p$, we can derive that
$$
\aligned
-\int_{\Omega}\big(u_0+H_0\big)\Delta u_0
=\int_{\Omega_{\rho_0v_0}}\frac{8(1+\alpha)^2|z|^{2\alpha}}{(1+|z|^{2(1+\alpha)})^2}
\left[\log\frac{(\varepsilon_0\mu_0)^{-4}}{(1+|z|^{2(1+\alpha)})^2}
+(1+\alpha)H(p,p)+O\left(\rho_0v_0\big|z\big|\right)
+O\left(\varepsilon_0^2\mu_0^2\right)
\right],
\endaligned
$$
where $\Omega_{\rho_0v_0}=\frac{1}{\rho_0v_0}(\Omega-\{p\})$. Note that
$$
\aligned
\int_{\Omega_{\rho_0v_0}}\frac{8(1+\alpha)^2|z|^{2\alpha}}{(1+|z|^{2(1+\alpha)})^2}
=8\pi(1+\alpha)+O\left(\varepsilon_0^2\mu_0^2\right),
\endaligned
$$
and
$$
\aligned
\int_{\Omega_{\rho_0v_0}}\frac{8(1+\alpha)^2|z|^{2\alpha}}{(1+|z|^{2(1+\alpha)})^2}
\log\frac{1}{(1+|z|^{2(1+\alpha)})^2}=-16\pi(1+\alpha)+O(\varepsilon^2_0\mu^2_0).
\endaligned
$$
Then
\begin{equation}\label{5.7}
\aligned
-\int_{\Omega}\big(u_0+H_0\big)\Delta u_0
=8\pi(1+\alpha)\big[(1+\alpha)H(p,p)-2-4\log(\varepsilon_0\mu_0)\big]
+O\left(\rho_0v_0\right)
+O\left(\varepsilon_0^2\mu_0^2|\log(\varepsilon_0\mu_0)|\right).
\endaligned
\end{equation}
For the second term of (\ref{5.6}), by (\ref{2.4}), (\ref{2.5}), (\ref{2.11})
and  the change of
variables $\rho_0v_0z=x-p$ we have that for any $j=1,\ldots,m$,
\begin{eqnarray}\label{5.8}
-\int_{\Omega}\big(u_j+H_j\big)\Delta u_0
=\int_{\Omega}\frac{8\varepsilon_0^2\mu_0^2(1+\alpha)^2|x-p|^{2\alpha}}{\,(\varepsilon_{0}^2\mu_{0}^2+|x-p|^{2(1+\alpha)})^2\,}
\left[\log
\frac{1}{(\varepsilon_{j}^2\mu_{j}^2+|x-\xi_j|^2)^2}
+H(x,\xi_j)+O\left(\varepsilon_j^2\mu_j^2\right)
\right]dx
\quad\,\,\,
&&\nonumber\\
=\int_{\Omega_{\rho_0v_0}}\frac{8(1+\alpha)^2|z|^{2\alpha}}{(1+|z|^{2(1+\alpha)})^2}
\left[\log
\frac{1}{|p-\xi_j|^4}
+H(p,\xi_j)+O\left(\rho_0v_0t^\beta\big|z\big|\right)+O\left(\varepsilon_j^2\mu_j^2t^{2\beta}\right)
\right]dz
&&\nonumber\\[1mm]
=8\pi(1+\alpha)G(p,\xi_j)
+O\left(\rho_0v_0t^\beta\right)
+O\left(\varepsilon_0^2\mu_0^2t^\beta\right)+O\left(\varepsilon_j^2\mu_j^2t^{2\beta}\right).
\qquad\qquad\qquad\qquad\quad\,\,
&&
\end{eqnarray}
As for the last term of (\ref{5.6}), by (\ref{2.4}), (\ref{2.5}), (\ref{2.11})
and  the change of
variables $\varepsilon_i\mu_iz=x-\xi_i$ we observe that for any $i,\,j=1,\ldots,m$,
$$
\aligned
-\int_{\Omega}\big(u_j+H_j\big)\Delta u_i
=&\int_{\Omega}\frac{8\varepsilon_i^2\mu_i^2}{(\varepsilon_{i}^2\mu_{i}^2+|x-\xi_i|^2)^2}
\left[\log\frac{1}{(\varepsilon_{j}^2\mu_{j}^2+|x-\xi_j|^2)^2}+
H(x,\xi_j)+O(\varepsilon_j^2\mu_j^2)
\right]dx\\
=&\int_{\Omega_{\varepsilon_i\mu_i}}\frac{8}{(1+|z|^2)^2}
\left[\log\frac{1}{(\varepsilon_{j}^2\mu_{j}^2+|\xi_i-\xi_j+\varepsilon_i\mu_iz|^2)^2}
+H(\xi_i,\xi_j)+O\left(\varepsilon_i\mu_i|z|\right)
+O\left(\varepsilon_j^2\mu_j^2\right)
\right]dz,
\endaligned
$$
where $\Omega_{\varepsilon_i\mu_i}=\frac{1}{\varepsilon_i\mu_i}(\Omega-\{\xi_i\})$.
Then for all $i,\,j=1,\ldots,m$,
\begin{equation}\label{5.9}
\aligned
-\int_{\Omega}\big(u_j+H_j\big)\Delta u_i
=\left\{\aligned
&8\pi \big[H(\xi_i,\xi_i)-2-4\log(\varepsilon_i\mu_i)\big]
+O\left(\varepsilon_i\mu_i\right)
\quad\,\,\forall\,\,\,i=j,\\[2.5mm]
&8\pi
G(\xi_i,\xi_j)+O\left(\varepsilon_i\mu_it^\beta\right)+O\left(\varepsilon_j^2\mu_j^2t^{2\beta}\right)
\qquad\,
\forall\,\,\,i\neq j.
\endaligned\right.
\endaligned
\end{equation}
On the other hand, by (\ref{2.18}), (\ref{2.19}),  (\ref{2.20}) and  the change of variables
$x=\varepsilon_0y=e^{-\frac12t}y$, we obtain
$$
\aligned
\int_{\Omega}|x-p|^{2\alpha}k(x)e^{-t\phi_1}e^U=&
\int_{\Omega_t}|\varepsilon_0 y-p|^{2\alpha}k(\varepsilon_0 y)e^{-t\big[\phi_1(\varepsilon_0 y)-1\big]}e^{U(\varepsilon_0 y)-2t}dy
\\[1mm]
=&\int_{\Omega_t\setminus\left[\bigcup_{i=1}^mB_{\frac1{\varepsilon_0 t^{2\beta}}}(\xi'_i)\cup
B_{\frac1{\varepsilon_0 t^{2\beta}}}(p')\right]}
W dy
+\int_{B_{\frac1{\varepsilon_0 t^{2\beta}}}(p')}W dy
+\sum\limits_{i=1}^m\int_{B_{\frac1{\varepsilon_0 t^{2\beta}}}(\xi'_i)}W dy.
\endaligned
$$
By (\ref{2.28}), (\ref{2.30})  and (\ref{2.32}) we obtain
$$
\aligned
\int_{\Omega_t\setminus\left[\bigcup_{i=1}^mB_{\frac1{\varepsilon_0 t^{2\beta}}}(\xi'_i)\cup
B_{\frac1{\varepsilon_0 t^{2\beta}}}(p')\right]}
W dy=&
\int_{\Omega_t\setminus\left[\bigcup_{i=1}^mB_{\frac1{\varepsilon_0 t^{2\beta}}}(\xi'_i)\cup
B_{\frac1{\varepsilon_0 t^{2\beta}}}(p')\right]}
O\left(\frac{\varepsilon_0^2e^{-t\phi_1(\varepsilon_0 y)}}
{\,|\varepsilon_0 y-p|^{4+2\alpha}\,}\prod_{i=1}^m\frac{1}{\,|\varepsilon_0 y-\xi_i|^4\,}
\right)dy\\[2mm]
=&\,O\left(1\right),
\endaligned
$$
and
$$
\aligned
\int_{B_{\frac1{\varepsilon_0 t^{2\beta}}}(p')}Wdy=&
\int_{B_{\frac1{\varepsilon_0 t^{2\beta}}}(p')}\left(\frac{\varepsilon_0}{\rho_0v_0}\right)^2\frac{8(1+\alpha)^2\big|\frac{\varepsilon_0 y-p}{\rho_0v_0}\big|^{2\alpha}}{\,\big(1+\big|\frac{\varepsilon_0 y-p}{\rho_0v_0}\big|^{2(1+\alpha)}\big)^2\,}
\left[1+O\left(\varepsilon_0t^{\beta}|y-p'|\right)
+
o\left(1\right)
\right]dy\\[2mm]
=&\,8\pi(1+\alpha)+o(1),
\endaligned
$$
and for any $i=1,\ldots,m$,
$$
\aligned
\int_{B_{\frac1{\varepsilon_0 t^{2\beta}}}(\xi'_i)}Wdy=
\int_{B_{\frac1{\varepsilon_0 t^{2\beta}}}(\xi'_i)}\frac{1}{\gamma_{i}^2}\frac{8}{\big(1+\big|\frac{y-\xi'_i}{\gamma_{i}}\big|^2\big)^2}
\left[1+O\left(\varepsilon_0t^{\beta}|y-\xi'_i|\right)
+
o\left(1\right)
\right]dy
=8\pi+o(1).
\endaligned
$$
Then
\begin{equation}\label{5.10}
\aligned
\int_{\Omega}|x-p|^{2\alpha}k(x)e^{-t\phi_1}e^U=O\left(1\right).
\endaligned
\end{equation}
Hence  by (\ref{5.1}), (\ref{5.6})-(\ref{5.10}) we conclude that
$$
\aligned
J_t\left(U(\xi)\right)=-16\pi(1+\alpha)\log(\varepsilon_0\mu_0)-16\pi
\sum_{i=1}^m\log(\varepsilon_i\mu_i)
+8\pi\sum\limits_{i=1}^{m}\left[(1+\alpha)G(p,\xi_i)
+\sum_{j=i+1}^mG(\xi_j,\xi_i)\right]
+O(1),
\endaligned
$$
which, together with the definitions  of $\varepsilon_0$, $\varepsilon_i$ in (\ref{2.7})
and the choices of $\mu_0$, $\mu_i$ in  (\ref{2.12})-(\ref{2.13}),
implies that expansion (\ref{5.5}) holds.

\vspace{1mm}

{\bf Step 2:}
For any integer $m\geq1$ and any $t$ large enough,  let us claim that
the
following expansion holds
\begin{equation}\label{5.11}
\aligned
F_t(\xi)=J_t\big(U(\xi)\big)+o(1)
\endaligned
\end{equation}
uniformly on points $\xi=(\xi_1,\ldots,\xi_m)\in\mathcal{O}_t$.
Indeed, set
\begin{equation}\label{5.12}
\aligned
I_t(\omega)=\frac12\int_{\Omega_t}
|\nabla \omega|^2-\int_{\Omega_t}|\varepsilon_0y-p|^{2\alpha}q(y,t)e^\omega,
\,\quad\,\,\omega\in H_0^1(\Omega_t).
\endaligned
\end{equation}
By (\ref{2.7}) and (\ref{2.19}) we obtain
$$
\aligned
F_t(\xi)-
J_t\big(U(\xi)\big)=I_t\big(V(\xi')+\phi_{\xi'}
\big)-I_t\big(V(\xi')
\big).
\endaligned
$$
Using  $DI_t(V+\phi_{\xi'})[\phi_{\xi'}]=0$, a Taylor expansion and an integration by parts, we give
\begin{eqnarray*}
F_t(\xi)-
J_t\big(U(\xi)\big)=\int_0^1D^2I_t(V+\tau\phi_{\xi'})[\phi_{\xi'}]^2(1-\tau)d\tau
\,\,\qquad\qquad\,\,
\quad\qquad\qquad\,\,
&&\\
=\int_0^1\left\{
\int_{\Omega_t}\big[N(\phi_{\xi'})+E\big]\phi_{\xi'}+
W\big[1-e^{\tau\phi_{\xi'}}\big]\phi_{\xi'}^2
\right\}(1-\tau)d\tau.&&
\nonumber
\end{eqnarray*}
Thanks to $\|\phi_{\xi'}\|_{L^{\infty}(\Omega_t)}\leq C t\max\big\{(\rho_0v_0)^{\min\{1,2(\alpha-\hat{\alpha})\}}t^{\beta},\,\varepsilon_0\gamma_1t^{\beta},\,\ldots,\,\varepsilon_0\gamma_mt^{\beta},\,
\|e^{-\frac{1}{2}t\phi_1}\|_{L^{\infty}(\Omega_t)}\big\}$
and the estimates in Lemma $3.2$ and Proposition $4.1$,
 we have readily
$$
\aligned
F_t(\xi)-
J_t\left(U(\xi)\right)=O\left(t
\max\left\{\big(\rho_0v_0\big)^{\min\{2,4(\alpha-\hat{\alpha})\}}t^{2\beta},\,\big(\varepsilon_0\gamma_1t^{\beta}\big)^2,\,\ldots,\,
\big(\varepsilon_0\gamma_mt^{\beta}\big)^2,\,
\|e^{-t\phi_1}\|_{L^{\infty}(\Omega_t)}\right\}\right)
=o\left(1\right).
\endaligned
$$
The continuity in $\xi$ of the above  expression is inherited from that of $\phi_{\xi'}$  in the $L^\infty$ norm.

\vspace{1mm}

{\bf Step 3:} Proof of Proposition 5.1.
Let $\xi_t=(\xi_{1,t},\ldots,\xi_{m,t})$
be the maximizer of $F_t$ over $\overline{\mathcal{O}}_t$. We need to prove that $\xi_t$
belongs to the interior of $\mathcal{O}_t$.
First, we obtain a lower bound for $F_t$ over $\overline{\mathcal{O}}_t$.
Let us fix the point
$p$ as a strict local maximum point of $\phi_1$ in $\Omega$ and set
$$
\aligned
\xi^0_i=p+
\frac{1}{\sqrt{t}}\widehat{\xi}_i,
\endaligned
$$
where $\widehat{\xi}=(\widehat{\xi}_1,\ldots,\widehat{\xi}_m)$ is
a $m$-regular polygon in $\mathbb{R}^2$.
Clearly,
$\xi^0=(\xi^0_1,\ldots,\xi^0_m)\in\mathcal{O}_t$
because  $\beta>1$ and $\phi_1(\xi^0_i)=1+O(t^{-1})$.
By (\ref{5.5}) and (\ref{5.11})  we find
\begin{eqnarray}\label{5.13}
\max\limits_{\xi\in\overline{\mathcal{O}}_t}
F_t(\xi)\geq
8\pi(1+\alpha)t
+8\pi t\sum_{i=1}^m
\phi_1(\xi_i^0)
+16\pi(2+\alpha)\sum_{i=1}^m\log|\xi_i^0-p|
+16\pi\sum_{i\neq j}^m
\log|\xi_i^0-\xi_j^0|+O\left(1\right)
&&\nonumber\\
\geq
8\pi(m+1+\alpha) t
-8\pi m(m+1+\alpha)
\log t+O(1).
\qquad\qquad\qquad\qquad\qquad
\qquad\qquad\qquad\,\,\,\,
&&
\end{eqnarray}
Next, we suppose $\xi_t=(\xi_{1,t},\ldots,\xi_{m,t})\in\partial\mathcal{O}_t$.
Then there  exist three possibilities:\\
C1. \,\,There exists an $i_0$ such that
$\phi_1(\xi_{i_0,t})=1-\frac{1}{\sqrt{t}}$;\\
C2. \,\,There exist indices $i_0$, $j_0$,
$i_0\neq j_0$ such that
$|\xi_{i_0,t}-\xi_{j_0,t}|=t^{-\beta}$;\\
C3. \,\,There exists an  $i_0$
such that
$|\xi_{i_0,t}-p|=t^{-\beta}$.\\
For the first case,   we have
\begin{eqnarray}\label{5.14}
\max\limits_{\xi\in\overline{\mathcal{O}}_t}
F_t(\xi)=F_t(\xi_t)
\leq
8\pi(1+\alpha)t
+8\pi t\left[(m-1)+1-\frac1{\sqrt{t}}\right]+O\big(\log t\big),
\end{eqnarray}
which contradicts to (\ref{5.13}).
For the second case, we have
\begin{equation}\label{5.15}
\aligned
\max\limits_{\xi\in\overline{\mathcal{O}}_t}
F_t(\xi)=F_t(\xi_t)
\leq
8\pi(m+1+\alpha) t -16\pi\beta
\log t
+O\big(1\big).
\endaligned
\end{equation}
For the last case, we have
\begin{equation}\label{5.16}
\aligned
\max\limits_{\xi\in\overline{\mathcal{O}}_t}
F_t(\xi)=F_t(\xi_t)
\leq
8\pi(m+1+\alpha) t -16\pi(2+\alpha)\beta
\log t
+O\big(1\big).
\endaligned
\end{equation}
Combining (\ref{5.15})-(\ref{5.16}) with (\ref{5.13}), we give
\begin{equation}\label{5.17}
\aligned
16\pi(2+\alpha)\beta\log t
+O\big(1\big)\leq
8\pi m(m+1+\alpha)
\log t+O(1),
\endaligned
\end{equation}
which  is impossible by the choice of $\beta$ in (\ref{2.3}).
\end{proof}

\section{Proof of Theorem 1.1}
\noindent {\bf Proof of Theorem 1.1.}
According to Proposition 4.1, we have that for
any integer $m\geq1$, any points
$\xi=(\xi_1,\ldots,\xi_m)\in\mathcal{O}_t$ and any $t$ large enough,
there exists a function
$\phi_{\xi'}$    such that
$$
\aligned
-\Delta\big(V(\xi')+\phi_{\xi'}\big)-|\varepsilon_0y-p|^{2\alpha}q(y,t)e^{V(\xi')+\phi_{\xi'}}
=\sum\limits_{i=1}^m\sum\limits_{j=1}^2c_{ij}(\xi')\chi_iZ_{ij},
\qquad
\int_{\Omega_t}\chi_iZ_{ij}\phi_{\xi'}=0
\,\,\,\,\,
\endaligned
$$
for some coefficients $c_{ij}(\xi')$,
$i=1,\ldots,m$, $j=1,2$.
%$$
%\aligned
%\left\{
%\aligned
%&-\Delta\big(V(\xi')+\phi_{\xi'}\big)-|\varepsilon_0y-p|^{2\alpha}q(y,t)\exp\big\{V(\xi')+\phi_{\xi'}\big\}
%=\sum\limits_{i=1}^m\sum\limits_{j=1}^2c_{ij}(\xi')\chi_iZ_{ij}
%\quad
%\textrm{in}\,\,\,\,\,\Omega_t,\\
%&\,\phi_{\xi'}=0
%\qquad\qquad\qquad\qquad\qquad
%\qquad\qquad\qquad\qquad\qquad
%\qquad\quad\,\,
%\qquad\qquad\qquad\qquad
%\textrm{on}
%\,\,\,\po_t,\\[1mm]
%&\int_{\Omega_t}\chi_iZ_{ij}\phi_{\xi'}=0
%\qquad\qquad\qquad\qquad\qquad
%\qquad\qquad\qquad\qquad\qquad
%\forall\,\,\,i=1,\ldots,m,
%\,\,\,j=1,2,
%\endaligned
%\right.
%\endaligned
%$$
Therefore, in order to
construct a solution to problem
(\ref{2.17}) and hence to the original problem  (\ref{1.4}),
we need to adjust  $\xi$ in $\mathcal{O}_t$ such that
the above coefficients $c_{ij}(\xi')$  satisfy
\begin{equation}\label{6.1}
\aligned
c_{ij}(\xi')=0\,\,\quad\,\,\textrm{for all}\,\,\,i=1,\ldots,m,\,\,j=1,2.
\endaligned
\end{equation}
On the other hand,  from Proposition 5.1, there is
a $\xi_t=(\xi_{1,t},\ldots,\xi_{m,t})\in\mathcal{O}_t^o$
that achieves the maximum for the maximization problem in Proposition 5.1.
Let
$\omega_t=V(\xi'_t)+\phi_{\xi'_t}$.
Then we have
\begin{equation}\label{6.2}
\aligned
\partial_{\xi_{kl}}F_t(\xi_t)=0\,\,\quad\,\,\textrm{for all}\,\,\,k=1,\ldots,m,\,\,l=1, 2.
\endaligned
\end{equation}
Notice that by (\ref{5.1}), (\ref{5.2}) and (\ref{5.12}),
$$
\aligned
\partial_{\xi_{kl}}F_t(\xi_t)=&\,\partial_{\xi_{kl}}J_t\big(U(\xi_t)+\tilde{\phi}(\xi_t)
\big)
=\,\frac1{\varepsilon_0}
\partial_{\xi'_{kl}}I_t\big(V(\xi_t')+\phi_{\xi_t'}
\big)\\
=&\,\frac1{\varepsilon_0}\left\{
\int_{\Omega_t}\nabla\omega_t\nabla\left[\partial_{\xi'_{kl}}V(\xi'_t)+\partial_{\xi'_{kl}}\phi_{\xi'_t}
\right]-\int_{\Omega_t}|\varepsilon_0y-p|^{2\alpha}q(y,t)e^{\omega_t}\left[\partial_{\xi'_{kl}}V(\xi'_t)+\partial_{\xi'_{kl}}\phi_{\xi'_t}
\right]
\right\}.
\endaligned
$$
Then for all $k=1,\ldots,m$
and $l=1, 2$,
$$
\aligned
\sum\limits_{i=1}^m\sum\limits_{j=1}^2c_{ij}(\xi'_t)\int_{\Omega_t}\chi_iZ_{ij}
\left[\partial_{\xi'_{kl}}V(\xi'_t)+\partial_{\xi'_{kl}}\phi_{\xi'_t}
\right]=0.
\endaligned
$$
Since
$\partial_{\xi'_{kl}}V(\xi'_t)(y)=Z_{kl}(y)+O\big(\varepsilon_0t^\beta\big)$
and $\|\partial_{\xi'_{kl}}\phi_{\xi'_t}\|_{L^{\infty}(\Omega_t)}
\leq Ct^{2}\max\big\{(\rho_0v_0)^{\min\{1,2(\alpha-\hat{\alpha})\}}t^{\beta},
\varepsilon_0\gamma_1t^{\beta},\ldots,\varepsilon_0\gamma_mt^{\beta},
\\
\|e^{-\frac{1}{2}t\phi_1}\|_{L^{\infty}(\Omega_t)}\big\}$,
by (\ref{2.7}), (\ref{2.14}), (\ref{2.15}) and (\ref{2.23})
we get the validity of a system of equations of the form
\begin{eqnarray}\label{6.3}
\sum\limits_{i=1}^m\sum\limits_{j=1}^2c_{ij}(\xi_t')\int_{\Omega_t}\chi_iZ_{ij}\big[
Z_{kl}(y)+o\left(1\right)
\big]=0,\,\ \ \ \ \,\,\,\,\,k=1,\ldots,m,\,\,l=1,2.
\end{eqnarray}
Note that
$$
\aligned
\int_{\Omega_t}\chi_iZ_{ij}
Z_{kl}(y)=\left\{
\aligned
&\int_{\mathbb{R}^2}\chi(|z|)\mathcal{Z}_{j}(z)\mathcal{Z}_{l}(z)dz=C\delta_{jl}
\,\,\ \,\,\textrm{if}\,\,\,\,i=k,\\[1mm]
&O\left(\varepsilon_0\gamma_i t^\beta\right)
\,\quad\quad\qquad\,\quad\qquad\qquad\,\,
\textrm{if}\,\,\,\,i\neq k.
\endaligned
\right.
\endaligned
$$
Hence the  coefficient matrix of  system (\ref{6.3}) is strictly diagonal dominant and then
$c_{ij}(\xi_t')=0$ for all $i=1,\ldots,m$, $j=1,2$.
As a consequence, we obtain a solution $u_t$ to problem (\ref{1.4}) of the form
$U(\xi_t)+\tilde{\phi}(\xi_t)$ with
the qualitative properties  as predicted in
Theorem 1.1.\qquad\qquad\qquad\qquad\qquad\qquad\qquad
\qquad\qquad\qquad\qquad\qquad\qquad\qquad
\qquad\quad$\square$

\end{document}